# LARGE DEVIATIONS FOR A CLASS OF NONHOMOGENEOUS MARKOV CHAINS[1]

By Zach Dietz and Sunder Sethuraman

*Tulane University and Iowa State University*

Large deviation results are given for a class of perturbed nonhomogeneous Markov chains on finite state space which formally includes some stochastic optimization algorithms. Specifically, let $\{P_n\}$ be a sequence of transition matrices on a finite state space which converge to a limit transition matrix $P$. Let $\{X_n\}$ be the associated nonhomogeneous Markov chain where $P_n$ controls movement from time $n-1$ to $n$. The main statements are a large deviation principle and bounds for additive functionals of the nonhomogeneous process under some regularity conditions. In particular, when $P$ is reducible, three regimes that depend on the decay of certain "connection" $P_n$ probabilities are identified. Roughly, if the decay is too slow, too fast or in an intermediate range, the large deviation behavior is trivial, the same as the time-homogeneous chain run with $P$ or nontrivial and involving the decay rates. Examples of anomalous behaviors are also given when the approach $P_n \to P$ is irregular. Results in the intermediate regime apply to geometrically fast running optimizations, and to some issues in glassy physics.

**1. Introduction.** The purpose of this paper is to provide some large deviation bounds and principles for a class of nonhomogeneous Markov chains related to some popular stochastic optimization algorithms such as Metropolis and simulated annealing schemes. In a broad sense, these algorithms are stochastic perturbations of steepest descent or "greedy" procedures to find the global minimum of a function $H$ and are in the form of nonhomogeneous Markov chains whose connecting transition kernels converge to a limit kernel associated with steepest descent.

Received December 2003; revised February 2004.
[1]Supported in part by NSF Grant DMS-00-03811.
*AMS 2000 subject classifications.* Primary 60J10; secondary 60F10.
*Key words and phrases.* Large deviations, nonhomogeneous, Markov, optimization, geometric cooling, glassy models.

For instance, in the Metropolis algorithm on finite state space $\Sigma$, the transition kernel connecting times $n-1$ and $n$ is given by

$$(1.1) \quad P_n(i,j) = \begin{cases} g(i,j)\exp\{-\beta_n(H(j)-H(i))_+\}, & \text{for } j \neq i, \\ 1 - \sum_{l \neq i} P_n(i,l), & \text{for } j = i, \end{cases}$$

where $g$ is an irreducible transition function and $\beta_n$ represents an inverse temperature parameter which diverges, $\beta_n \to \infty$. Here, the limit kernel $P = \lim_n P_n$ corresponds to steepest descent in that jumps from $i$ to $j$ when $H(j) > H(i)$ are not allowed.

These types of schemes are intensively used in image analysis [35], neural networks [4], statistical physics of glassy systems and combinatorial optimizations [26]. More general tutorials include [8, 16, 17] and [32].

Virtually all previous large deviations work with respect to optimization chains has been through Freidlin–Wentzell-type methods [14]. This approach is to consider a sequence of time-homogeneous Markov chains, parametrized by temperature, which approaches the steepest descent chain as the temperature cools, and then to transfer "short time" large deviation estimates to a single related system in which temperature varies with time. For instance, with respect to the Metropolis algorithm, by studying the sequence of time-homogeneous chains $\{X_\cdot^\beta : \beta \geq 0\}$ where $\beta_n \equiv \beta$ and $\beta \uparrow \infty$, estimates can be made on the nonhomogeneous chain where $\beta_n$ varies. Although this approach has had much success, especially related to statistical physics metastability questions, it seems that only large deviation bounds are recovered for the position of the nonhomogeneous process rather than large deviation principles (LDPs) (see [2, 5, 6, 8, 9] and references therein). It would be then natural to ask about LDPs for empirical averages which are more regular objects than the positions.

In a different, more general vein, LDPs have been shown for independent nonidentically distributed variables whose Cesaro empirical averages converge [29], and also for some types of Gibbs measures, which include nonhomogeneous chains whose connecting transition kernels are positive entrywise and converge in Cesaro mean to a positive limit matrix [31].

Other work in the literature treats an intermediate case of nonhomogeneity, namely Markov chains whose transition kernels are chosen at random from a time-homogeneous process. The results here are then to prove an LDP for almost all realized nonhomogeneous Markov chains chosen in this fashion [20, 30]. Also, we note that an LDP has been shown for a class of near irreducible time-homogeneous processes that satisfy some mixing conditions [1].

In this context, we develop here an LDP in natural scale $n$ with explicit rate function for the empirical averages of nonhomogeneous Markov chains on finite state spaces whose transition kernels converge to the general limit matrix which allows for reducibility, a key concern in optimization schemes. We note the methods used here differ from Freidlin–Wentzell-type arguments in that they focus on the nonhomogeneous process itself rather than homogeneous approximations. The specific techniques used are constructive and involve various "surgeries" of path realizations and some coarse graining.

Let $\Sigma = \{1, 2, \ldots, \mathfrak{r}\}$ be a finite set of points. Let $P_n = \{p_n(i,j) : i, j \in \Sigma\}$ be a sequence of $\mathfrak{r} \times \mathfrak{r}$ stochastic matrices for $n \geq 1$ and let $\pi$ be a distribution on $\Sigma$. Let now $\mathbb{P}_\pi = \mathbb{P}_\pi^{\{P_n\}}$ be the (nonhomogeneous) Markov measure on the sequence space $\Sigma^\infty$ with Borel sets $\mathcal{B}(\Sigma^\infty)$ that correspond to initial distribution $\pi$ and transition kernels $\{P_n\}$. That is, with respect to the coordinate process $X_0, X_1, \ldots$, we have the Markov property

$$\mathbb{P}_\pi(X_{n+1} = j | X_0, X_1, \ldots, X_{n-1}, X_n = i) = p_{n+1}(i, j)$$

for all $i, j \in \Sigma$ and $n \geq 0$. We see then that $P_{n+1}$ controls "transitions" between times $n$ and $n+1$.

We now specify the class of nonhomogeneous processes focused on in this article. Let $\pi$ be a distribution and let $P = \{p(i,j)\}$ be a stochastic matrix on $\Sigma$. Define the collection

$$\mathbb{A}(P) = \{\mathbb{P}_\pi^{\{P_n\}} : P_n \to P\},$$

where the convergence $P_n \to P$ is elementwise, that is, $\lim_{n \to \infty} p_n(i,j) = p(i,j)$ for all $i, j \in \Sigma$. The collection $\mathbb{A}$ can be thought of as perturbations of the time-homogeneous Markov chain run with $P$ and is a natural class in which to explore how nonhomogeneity enters into the large deviation picture.

We also remark that this class has been studied in connection with other types of problems such as ergodicity [19], laws of large numbers [34, 35] and fluctuations [18]. See also [24] and [15] for some laws of large numbers and fluctuation results for generalized annealing algorithms and Markov chains with rare transitions.

Let now $f : \Sigma \to \mathbb{R}^d$ be a $(d \geq 1)$-dimensional function. Let also $\mathbb{P}_\pi \in \mathbb{A}(P)$ be a $P$-perturbed nonhomogeneous Markov measure. In terms of the coordinate process, define the additive sum $Z_n = Z_n(f)$ for $n \geq 1$ by

$$Z_n = \frac{1}{n} \sum_{i=1}^n f(X_i).$$

The specific goal of this paper is to understand the large deviation behavior of the induced distributions of $\{Z_n : n \geq 1\}$ with respect to $\mathbb{P}_\pi$ in scale $n$. That is, we search for a rate function $\mathbb{J}$ so that for Borel sets $B \subset \mathbb{R}^d$,

$$-\inf_{z \in B^o} \mathbb{J}(z) \leq \liminf \frac{1}{n} \log \mathbb{P}_\pi(Z_n \in B)$$
$$\leq \limsup \frac{1}{n} \log \mathbb{P}_\pi(Z_n \in B) \leq -\inf_{z \in \overline{B}} \mathbb{J}(z).$$

An immediate question which comes to mind is whether these large deviations for the nonhomogeneous chain, if they exist, differ from the deviations with respect to the time-homogeneous chain run with $P$. The general answer found in our work is "Yes" and "No," and as might be suspected depends on the rate of convergence $P_n \to P$ and the structure of the limit matrix $P$.

More specifically, when $P$ is irreducible, it turns out that the large deviation of behavior of $\{Z_n\}$ under $\mathbb{P}_\pi$ is the same as that under the time-homogeneous chain associated with $P$ and independent of the rate of convergence of $P_n$ to $P$. (Note that [31] covers the case $P$ is positive entrywise and [29] covers the case when each $P_n$ has identical rows.)

Perhaps the more interesting case is when the target matrix $P$ is reducible. Indeed, this is the case with stochastic optimization algorithms where $H$ has several local minima, for example, with respect to the Metropolis process, the local minima sets of $H$ do not communicate in the limit steepest descent chain. In this situation, the large deviations of $\{Z_n\}$ depend on both the type of reducibilities of $P$ and the decay rate, with respect to $P_n$, of certain "connection probabilities" between $P$-irreducible sets, and fall into three categories. Namely, when the decay is fast, or superexponential, the large deviation behavior is the same as for the time-homogeneous Markov chain run under $P$; when the speed is slow, or subexponential, a trivial large deviation behavior is obtained; finally, when the speed is intermediate, or when the connection probabilities are on the order $e^{-Cn}$, a nontrivial behavior is found which differs from stationarity.

We remark now, in terms of applications, the intermediate processes are important in situations such as (i) fast annealing simulations, and (ii) models of glass formation.

(i) In Metropolis-type procedures, classic convergence theorems mandate that the temperatures satisfy $\beta_n \leq C_H \log n$ with respect to a known constant $C_H$ for the process to converge to the global minima set of $H$ (cf. [8] and [17]):

$$\lim_{n \to \infty} \mathbb{P}_\pi(X_n \in \text{ global minima set of } H) = 1.$$

However, with only finite time and resources, the optimal logarithmic speed is too slow to yield good results. In fact, in violation of classic results, exponentially fast schemes where $\beta_n \sim n$ are often used for which the process may actually converge to a nonglobal but local minimum of $H$. Whereas connecting probabilities between local minima sets are on the order of $\exp(-C\beta_n)$, these chains fit naturally in the intermediate framework mentioned above (cf. discussion after Corollary 3.1). Although there are some good error bounds for these geometrically cooling experiments in finite time [7], it seems the structure of the associated dynamics is not that well understood (cf. [35], Section 6.2).

(ii) In the manufacture of glass, a hot, fired material is quickly quenched into a substance which is not quite solid or liquid. The interpretation is that under rapid cooling the constructed glass is caught in a local energy optimum associated with some spatial disorder—not the regularly structured global one associated with a solid—from which over much longer time scales it may move to other states [22]. Such glassy systems are intensively studied in the literature. Two rough concerns can be identified: What are the typical glass landscapes which specify the local optima and what are the dynamics of the quick quenching phase and beyond? Much discussion is focused on the first concern [27], but even in systems where statics are quantified, dynamical questions remain open [23], Part IV, and [26]. However, with respect to metastability, as mentioned earlier, much work has been accomplished (cf. [3, 11] and [33] and references therein). Less work has been done though when certain time inhomogeneities are severe, say on exponential scale $e^{-Cn}$, in the context of Metropolis models in the intermediate regime.

At this point, we observe, as alluded to above in the two examples, that (from Borel–Cantelli arguments) the typical large scale picture of general intermediate speed nonhomogeneous Markov chains is to get trapped in one of the irreducible sets that correspond to the limit $P$ (e.g., the local $H$-minima sets in the Metropolis scheme). In this sense, the large deviations rate function $\mathbb{J}$, found with respect to averages $\{Z_n\}$, is relevant to understanding how atypical deviations arise, namely how the process average can "survive" for long times, that is, how $Z_n \sim z$ for large $n$ when $z$ is *not* a $P$-irreducible set average. More specifically, when $P$ corresponds to $K \geq 2$ irreducible sets $\{C_{\zeta_j}\}$, we show that $\mathbb{J}$ is an optimization between two types of costs and is in the form

$$\mathbb{J}(z) = \min_{\sigma \in \mathbb{S}} \inf_{\mathbf{v} \in \Omega} \inf_{\mathbf{x} \in D(\mathbf{v},z)} -\sum_{i=1}^{K-1}\left(\sum_{j=1}^{i} v_j\right)\mathcal{U}(\zeta_{\sigma(i)},\zeta_{\sigma(i+1)}) + \sum_{i=1}^{K} v_i \mathbb{I}_{\zeta_{\sigma(i)}}(x_i).$$

Here $\mathbb{I}_{\zeta_j}$ is the rate function for the $P$ time-homogeneous chain restricted to $C_{\zeta_j}$ and represents a "resting" cost of moving within $C_{\zeta_j}$, and $\mathcal{U}(\zeta_j, \zeta_k)$ is a large deviation "routing" cost of traveling between $C_{\zeta_j}$ and $C_{\zeta_k}$. Also, $\mathbb{S}$ and $\Omega$ are the sets of permutations and probabilities on $\{1, 2, \ldots, K\}$, respectively, and $D(\mathbf{v}, z)$ is the set of vectors $\mathbf{x}$ such that $\sum_{j=1}^{K} v_j x_j = z$. The intuition then is that $Z_n$ optimally deviates to $z$ by visiting sets $\{C_{\zeta_j}\}$ finitely many times, in a certain order $\sigma$ with time proportions $\mathbf{v}$, so that the average $z$ is maintained, and resting and routing costs are minimized.

Our main theorem (Theorem 3.3) is that under some natural regularity conditions on the approach $P_n \to P$, the average $Z_n$ satisfies an LDP with rate function $\mathbb{J}$. When $\mathcal{U} \equiv -\infty$ or $\mathcal{U} \equiv 0$, that is, when connection probabilities vanish too fast or too slow, the rate $\mathbb{J}$ reduces to the rate function for the time-homogeneous chain run under $P$ or a trivial rate. When the connections are exponential, $-\infty < \mathcal{U} < 0$ and $\mathbb{J}$ nontrivially incorporates the convergence exponents (Corollary 3.1). Some comments on the Metropolis algorithm are made at the end of Section 3. When the approach is irregular, large deviation bounds (Theorems 3.1 and 3.2) and examples (Section 12) of anomalous behaviors are also given.

Finally, it is natural to ask about the large deviations on scales $\alpha_n$ different from scale $n$, that is, the liminf and limsup limits of $(1/\alpha_n)\log \mathbb{P}_\pi(Z_n \in B)$. The metaresult should be, if the typical system behavior is to be absorbed into certain sets, the analogous large deviation (LD) behavior holds in scale $\alpha_n$ with revised resting and routing costs reflecting the scale. In fact, with respect to the Metropolis model, by the methods in this article, large deviation bounds and principles in scale $\beta_n$ can be derived as long as $\liminf \beta_n/n^\theta = \infty$ for some $\theta > 0$. In principle, similar results should hold when $\beta_n > C \log n$ and $C > 1$, although this is not pursued here. On the other hand, large deviation principles in scale $\beta_n \leq \log n$ are of a completely different category, because in this case there is no local minima absorption (see, however, [9] for LD bounds with respect to metastability concerns).

**2. Preliminaries.** We now recall and develop some definitions and notation before arriving at the main theorems. Throughout, we use the convention that $\pm\infty \cdot 0 = 0$ and $\log 0 = -\infty$.

2.1. *Rate functions and extended LDP.* Let $\mathbb{I}:\mathbb{R}^d \to \mathbb{R} \cup \{\infty\}$ be an extended real-valued function. We say that $\mathbb{I}$ is an *extended rate function* if $\mathbb{I}$ is lower semicontinuous and, further, that $\mathbb{I}$ is a *good extended rate function* if, in addition, the level sets of $\mathbb{I}$, namely $\{x : \mathbb{I}(x) \leq a\}$ for $a \in \mathbb{R}$, are compact. This definition extends the usual notion of rate function where negative values are not allowed (cf. [10], Section 1.2). Namely, we say $\mathbb{I}$ is a (*good*) *rate function* if $\mathbb{I}:\mathbb{R}^d \to [0,\infty]$ is a (good) extended rate function.

We denote $Q_\mathbb{I} \subset \mathbb{R}^d$ as the domain of finiteness, $Q_\mathbb{I} = \{x \in \mathbb{R}^d : \mathbb{I}(x) < \infty\}$. We also recall the standard notation for $B \subset \mathbb{R}^d$ that $\mathbb{I}(B) = \inf_{x \in B} \mathbb{I}(x)$.

Let now $\{\mu_n : n \geq 1\}$ be a sequence of nonnegative measures with respect to Borel sets on $\mathbb{R}^d$. We say that $\{\mu_n\}$ satisfies a large deviation principle with (extended) rate function $\mathbb{I}$ if, for all Borel sets $B \subset \mathbb{R}^d$, we have

$$(2.1) \quad -\inf_{x \in B^o} \mathbb{I}(x) \leq \liminf \frac{1}{n} \log \mu_n(B) \leq \limsup \frac{1}{n} \log \mu_n(B) \leq -\inf_{x \in \overline{B}} \mathbb{I}(x).$$

2.2. *Nonnegative matrices.* Let $U = \{u(i,j)\}$ be a matrix on $\Sigma$ and let $C \subset \Sigma$ be a subset of states. Define $U_C = \{u(i,j) : i, j \in C\}$ as the corresponding submatrix. We say that $U_C$ is *nonnegative*, denoted $U_C \geq 0$, if all entries are nonnegative. Analogously, $U_C$ is *positive*, denoted $U_C > 0$, if its entries are all positive. We say a nonnegative matrix $U_C$ is *stochastic* if all rows add to 1, $\sum_{j \in C} u(i,j) = 1$ for all $i \in C$; of course, $U_C$ is *substochastic* when $\sum_{j \in C} u(i,j) \leq 1$ for all $i \in C$. Also, we say $U_C$ is *primitive* if there is an integer $k \geq 1$ such that $U_C^k > 0$ is positive. In addition, we say $U_C$ is *irreducible* if, for any $i, j \in C$, there is a finite path $i = x_0, x_1, \ldots, x_n = j$ in $C$ with positive weight, $U_C(x_0, x_1) \cdots U_C(x_{n-1}, x_n) > 0$. The *period* of a state $i \in C$ is defined as $d_C(i) = \text{g.c.d}\{n \geq 1 : U_C^n(i,i) > 0\}$. When $U_C$ is irreducible, all states in $C$ have the same period $d_C$. When $d_C = 1$, we say $U_C$ is *aperiodic*. Finally, note that $U_C$ is primitive $\Leftrightarrow U_C$ is irreducible and aperiodic $\Leftrightarrow (U_C)^{\mathfrak{r}} > 0$.

2.3. *Construction CON.* We now construct a sequence of nonnegative Markov-like measures. Let $U_k = \{u_k(i,j)\}$ for $1 \leq k \leq n$ be a sequence of $\mathfrak{r} \times \mathfrak{r}$ nonnegative matrices. Let also $\pi$ be a measure on $\Sigma$. Then define the nonnegative measure $\mathbb{U}_\pi$ on $\Sigma^n$ for $n \geq 1$, where $\mathbb{U}_\pi(X_0 \in B) = \pi(B)$ and

$$\mathbb{U}_\pi(\mathbf{X}_n \in B) = \sum_{x_0 \in \Sigma} \sum_{\mathbf{x}_n \in B} \pi(x_0) \prod_{i=1}^n u_i(x_{i-1}, x_i),$$

where $\mathbf{X}_n = \langle X_1, \ldots, X_n \rangle$ is the coordinate process up to time $n$. Let also $\mathbf{X}_i^j = \langle X_i, \ldots, X_j \rangle$ for $0 \leq i \leq j$ be the observations between times $i$ and $j$, and denote, for $0 \leq k \leq m \leq l$,

$$\mathbb{U}_{(k,\pi)}(\mathbf{X}_m^l \in B) = \mathbb{U}_\pi'(\mathbf{X}_{m-k}^{l-k} \in B),$$

where $\mathbb{U}_\pi'$ is made with respect to $U_i' = U_{i+k}$ for $i \geq 1$. When $\pi$ is the point mass $\delta_x$ for $x \in \Sigma$, we denote $\mathbb{U}_{(k,\delta_x)} = \mathbb{U}_{(k,x)}$ for simplicity.

The measure $\mathbb{U}_\pi$ shares the Markov property:

$$\mathbb{U}_\pi(\mathbf{X}_k \in A, \mathbf{X}_{k+1}^n \in B) = \sum_{x_0 \in \Sigma} \sum_{\mathbf{x}_k \in A} \sum_{\mathbf{x}_{k+1}^n \in B} \pi(x_0) \prod_{i=1}^n u_i(x_{i-1}, x_i)$$

(2.2)

$$= \sum_{\mathbf{x}_k \in A} \mathbb{U}_\pi(\mathbf{X}_k = \mathbf{x}_k) \mathbb{U}_{(k,x_k)}(\mathbf{X}_{k+1}^n \in B).$$

2.4. *LDP for homogeneous nonnegative processes.* Let $U$ be a nonnegative matrix on $\Sigma$. Let also $C \subset \Sigma$ and let $f : \Sigma \to \mathbb{R}^d$ be a subset and function on the state space.

For $\lambda \in \mathbb{R}^d$, define the "tilted" matrix $\Pi_{C,\lambda,f,U} = \Pi_{C,\lambda}$ by

$$\Pi_{C,\lambda} = \{u(i,j) e^{\langle \lambda, f(j) \rangle} : i, j \in C\}.$$

Suppose now that $C$ is such that $U_C$ is irreducible. Then $\Pi_{C,\lambda}$ is irreducible for all $\lambda$ and $f$, and we may define

(2.3) $\rho(C, \lambda) = \rho(C, \lambda; f, U)$     as the Perron–Frobenius eigenvalue of $\Pi_{C,\lambda}$

(cf. [10], Theorem 3.1.1, or [28]). Define also the extended function $\mathbb{I}_C = \mathbb{I}_{C,f,U} : \mathbb{R}^d \to \mathbb{R} \cup \{\infty\}$ by

$$\mathbb{I}_{C,f,U}(x) = \sup_{\lambda \in \mathbb{R}^d} \{\langle \lambda, x \rangle - \log \rho(C, \lambda)\}$$

and let $Q_C = Q_{\mathbb{I}_C}$ be its domain of finiteness.

Let now $\pi$ be a distribution on $\Sigma$ and let $\mathbb{U}_\pi$ be made from CON with $U_k = U$ for all $k \geq 1$. We call such a measure $\mathbb{U}_\pi$ a *homogeneous nonnegative process*. Also, for $x_0 \in C$, define the measures on $\mathbb{R}^d$ for $n \geq 2$ by

$$\mu_n(B) = \mathbb{U}_{x_0}(Z_n(f) \in B, \mathbf{X}_n \in C^n).$$

Define also for $1 \leq k \leq l$ that $Z_k^l = Z_k^l(f) = (1/l - k + 1) \sum_{i=k}^l f(X_i)$. Note, as $|\Sigma| < \infty$, that $f$ is bounded, $\|f\| = \max_{1 \leq i \leq d} \|f_i\|_{L^\infty} < \infty$ and so $Z_k^l$ varies within the closed cube $\mathbb{K} = \overline{B}_{cu}(0, \|f\|)$ of width $2\|f\|$ about the origin.

The following proposition is proved in the Appendix.

PROPOSITION 2.1.   *The function $\mathbb{I}_C$ and domain $Q_C$ satisfy the following criteria*:

1. *Domain $Q_C$ is a nonempty convex compact subset of the cube $\mathbb{K}$.*
2. *Function $\mathbb{I}_C$ is a good extended rate function. In fact, when $U_C$ is substochastic, $\mathbb{I}_C$ is a good rate function.*
3. *Function $\mathbb{I}_C$ is convex on $\mathbb{R}^d$ and strictly convex on the relative interior of $Q_C$. Also, when restricted to $Q_C$, $\mathbb{I}_C$ is uniformly continuous and hence bounded on $Q_C$.*
4. *Measure $\{\mu_n\}$ satisfies an LDP (2.1) with extended rate function $\mathbb{I}_C$.*

2.5. *Upper block form.*  For a stochastic matrix $P = \{p(i,j)\}$ on $\Sigma$, we now recall the upper block form. By reordering $\Sigma$ if necessary, the matrix $P$ may be put in the form

$$(2.4) \qquad P = \begin{bmatrix} U(0,0) & U(0,1) & \cdots & \cdots & U(0,M_0) \\ 0 & S(1) & 0 & \cdots & 0 \\ \vdots & 0 & \ddots & & \vdots \\ \vdots & & \ddots & \ddots & \vdots \\ 0 & \cdots & \cdots & 0 & S(M_0) \end{bmatrix},$$

where $1 \leq M_0 \leq \mathfrak{r}$ and $S(1), \ldots, S(M_0)$ are stochastic irreducible submatrices that correspond to disjoint subsets of recurrent states—denoted as *stochastic sets*—and submatrices $U(0,0), \ldots, U(0,M_0)$ correspond to transient states when they exist.

When there are transient states, the square block $U(0,0)$ itself may be decomposed as (cf. [28], Section 1.2)

$$U(0,0) = \begin{bmatrix} R(1) & V(1,2) & \cdots & \cdots & V(1,N_0) \\ 0 & R(2) & V(2,3) & \cdots & V(2,N_0) \\ \vdots & 0 & \ddots & & \vdots \\ \vdots & & \ddots & \ddots & \vdots \\ 0 & \cdots & \cdots & 0 & R(N_0) \end{bmatrix},$$

where $1 \leq N_0 \leq \mathfrak{r} - 1$ and $R(i)$ is either the $1 \times 1$ zero matrix or an irreducible submatrix that corresponds to a subset of transient states for $1 \leq i \leq N_0$. We call the $R(i) = [0]$ matrices and corresponding states *degenerate transient*, and the irreducible $R(i)$ and associated states *nondegenerate transient*, since returns to these states are, respectively, impossible and possible under the time-homogeneous chain run with $P$.

Define the number of degenerate transient submatrices as

$$N = \begin{cases} 0, & \text{when no transient states in } P, \\ |\{1 \leq i \leq N_0 : R(i) = [0]\}|, & \text{otherwise.} \end{cases}$$

Also let the number of nondegenerate and stochastic submatrices be

$$M = \begin{cases} M_0, & \text{when no transient states in } P, \\ (N_0 - N) + M_0, & \text{otherwise.} \end{cases}$$

It will be useful to rewrite the upper block form by inserting the form for $U(0,0)$ into (2.4). To this end, when there are transient states, let $P(i) = R(i)$ for $1 \leq i \leq N_0$ and let $P(i) = S(i - N_0)$ for $N_0 + 1 \leq i \leq N_0 + M_0$. When all states are recurrent, let $P(i) = S(i)$ for $1 \leq i \leq M_0$. Also, in the following discussion, let $T(i,j)$ for $i < j$ denote the appropriate "connecting" submatrix $U(\cdot,\cdot)$ or $V(\cdot,\cdot)$. We remark that $T(i,j)$ is a matrix of zeroes for $N_0 + 1 \leq i < j \leq N + M$.

We have now the canonical decomposition

$$P = \begin{bmatrix} P(1) & T(1,2) & \cdots & \cdots & T(1,N+M) \\ 0 & P(2) & T(2,3) & \cdots & T(2,N+M) \\ \vdots & 0 & \ddots & & \vdots \\ \vdots & & \ddots & \ddots & \vdots \\ 0 & \cdots & \cdots & 0 & P(N+M) \end{bmatrix}.$$

Let now $C_i = C_i(P) \subset \Sigma$ be the subset which corresponds to $P(i)$ so that $P_{C_i} = P(i) = \{p(x,y) : x, y \in C_i\}$ for $1 \leq i \leq N + M$. Define also the sets $\mathcal{D} = \mathcal{D}(P)$, $\mathcal{N} = \mathcal{N}(P)$, $\mathcal{M} = \mathcal{M}(P)$ and $\mathcal{G} = \mathcal{G}(P)$ by

$\mathcal{D} = \{i : P(i) \text{ degenerate transient}\},$

$\mathcal{N} = \{i : P(i) \text{ nondegenerate transient}\},$

$$\mathcal{M} = \{i : P(i) \text{ stochastic}\},$$
$$\mathcal{G} = \mathcal{N} \cup \mathcal{M} \qquad (= \{i : P(i) \text{ nondegenerate transient or stochastic}\}).$$

To link with previous notation, note that $N = |\mathcal{D}|$ and $M = |\mathcal{G}|$.

It will be convenient to enumerate the elements of $\mathcal{G}$ as $\mathcal{G} = \{\zeta_1, \zeta_2, \ldots, \zeta_M\}$. Whereas $P(i)$ is (sub)stochastic and irreducible for $i \in \mathcal{G}$, we may denote, with respect to $f : \Sigma \to \mathbb{R}^d$, the rate function $\mathbb{I}_i = \mathbb{I}_{C_i, f, P}$ and its domain of finiteness $Q_i = Q_{C_i}$. In addition, let

$$(2.5) \qquad p_{\min} = \min\{p(x, y) : p(x, y) \neq 0, x, y \in C_i, i \in \mathcal{G}\}$$

be the minimum positive transition probability in the irreducible submatrices of $P$.

Consider now a sequence of transition matrices $\{P_n\}$, where $P_n = \{p_n(i, j)\}$ converges to $P$. With respect to the sets $\{C_i(P) : 1 \leq i \leq N + M\}$ above for the matrix $P$, the $n$th step matrix $P_n$ can be put in the form

$$P_n = \begin{bmatrix} P_n(1) & T_n(1,2) & \cdots & & \cdots & T_n(1, N+M) \\ T_n(2,1) & P_n(2) & T_n(2,3) & & \cdots & T_n(2, N+M) \\ \vdots & T_n(3,2) & & & & \vdots \\ \vdots & & \ddots & & & \vdots \\ T_n(N+M, 1) & \cdots & & \cdots & T_n(N+M, N+M-1) & P_n(N+M) \end{bmatrix},$$

where $P_n(i) = (P_n)_{C_i} \to P(i)$ for $1 \leq i \leq N + M$, $T_n(i, j)$ governs $P_n$ transitions from $C_i$ to $C_j$, and $T_n(i, j) \to T(i, j)$ for $i < j$ and vanishes otherwise. As a warning, we note that the form above for $P_n$ is NOT the canonical decomposition of $P_n$.

2.6. *Routing costs and deviations.* Let $\mathbb{S}_M$ and $\Omega_M$ be the set of permutations and the collection of probability vectors on $\{1, 2, \ldots, M\}$,

$$\Omega_M = \left\{ \mathbf{v} \in \mathbb{R}^M : \sum_{i=1}^{M} v_i = 1, 0 \leq v_i \leq 1 \text{ for } 1 \leq i \leq M \right\}.$$

For $\mathbf{v} \in \Omega_M$ and $z \in \mathbb{R}^d$, define the set of convex combinations

$$(2.6) \qquad D(M, \mathbf{v}, z) = \left\{ \mathbf{x} = \langle x_1, \ldots, x_M \rangle \in (\mathbb{R}^d)^M : \sum_{i=1}^{M} v_i x_i = z \right\}.$$

Let also $U = \{u(i, j) : 1 \leq i, j \leq M\}$ be a matrix of extended nonpositive real numbers. For a permutation $\sigma \in \mathbb{S}_M$, $\mathbf{v} \in \Omega_M$, $\mathbf{x} \in (\mathbb{R}^d)^M$ and $z \in \mathbb{R}^d$, define the extended functions

$$C_{\mathbf{v}, U}(\sigma, \mathbf{x}) = \begin{cases} -\sum_{i=1}^{M-1} \left( \sum_{j=1}^{i} v_j \right) u(\zeta_{\sigma(i)}, \zeta_{\sigma(i+1)}) + \sum_{i=1}^{M} v_i \mathbb{I}_{\zeta_{\sigma(i)}}(x_i), & \text{for } M \geq 2, \\ \mathbb{I}_1(x_1), & \text{for } M = 1, \end{cases}$$

and

$$\mathbb{J}_U(z) = \inf_{\mathbf{v}\in\Omega_M} \inf_{\mathbf{x}\in D(M,\mathbf{v},z)} \min_{\sigma\in S_M} C_{\mathbf{v},U}(\sigma,\mathbf{x}).$$

It will be shown that $\mathbb{J}_U$ is a good rate function (Proposition 4.1). Moreover, it will turn out, for well chosen routing cost matrices $U$, that $\mathbb{J}_U(z)$ measures various upper and lower large deviation rates of the additive sums $\{Z_n(f)\}$. Note that $\mathbb{J}_U$ is defined in terms of $\{\varsigma_i\} = \mathcal{G}$ and depends on $i \in \mathcal{D}$ only possibly through the routing cost $U$, which makes sense since it would be too expensive to rest on degenerate transient states in any positive time proportion. Also, we observe when $M = 1$, that is, when any transient states with respect to $P$ do not allow returns, and $P$ corresponds to exactly one irreducible stochastic block, the function $\mathbb{J}_U(z) = \mathbb{I}_1(z)$ is independent of $U$.

2.7. *Upper and lower cost matrices.* With respect to a $\mathbb{P}_\pi \in \mathbb{A}(P)$, we now specify certain relevant upper and lower costs $U$ when $N + M \geq 2$. Define, for distinct $1 \leq i, j \leq N + M$,

(2.7) $$t(n,(i,j)) = \max_{\substack{x\in C_i \\ y\in C_j}} p_n(x,y)$$

and the extended nonpositive numbers

$$\upsilon(i,j) = \limsup_{n\to\infty} \frac{1}{n}\log t(n,(i,j)) \quad \text{and} \quad \tau(i,j) = \liminf_{n\to\infty} \frac{1}{n}\log t(n,(i,j)).$$

Also, for $0 \leq k \leq N+M-2$, let $l_0 = i$, $l_{k+1} = j$ and let $L_k = \langle l_0, l_1, \ldots, l_k, l_{k+1}\rangle$ be a $(k+2)$-tuple of distinct indices. Now define the upper cost

(2.8) $$\mathcal{U}_0(i,j) = \max_{0\leq k\leq N+M-2} \max_{L_k} \sum_{s=0}^{k} \upsilon(l_s, l_{s+1})$$

and the lower cost

$$\mathcal{T}_0(i,j) = \max_{0\leq k\leq N+M-2} \max_{L_k} \sum_{s=0}^{k} \tau(l_s, l_{s+1}).$$

We remark briefly that $\mathcal{U}_0(i,j)$ and $\mathcal{T}_0(i,j)$ represent, respectively, maximal and minimal asymptotic travel costs of moving from $C_i$ to $C_j$ in $k+1 \leq N+M-1$ steps by visiting sets $\{C_i\}$ in the order $L_k$.

A more subtle lower cost $\mathcal{T}_1$ is the following. Let $0 \leq k \leq N+M-2$, $l_0 = i$, $l_{k+1} = j$ and $L_k$ be as before. Let also

(2.9) $$1 \leq q_0, \quad q_{k+1} \leq \mathfrak{r} \quad \text{and} \quad \text{when } k > 1 \quad \text{and} \quad 1 \leq s \leq k,$$

$$\text{let } 1 \leq q_s \leq \mathfrak{r}+1$$

and call $Q_k = \langle q_0, \ldots, q_{k+1}\rangle$. Let $\mathbf{x}^0 = \langle x_1^0, \ldots, x_{q_0}^0\rangle$ and $\mathbf{x}^{k+1} = \langle x_1^{k+1}, \ldots, x_{q_{k+1}}^{k+1}\rangle$ be vectors with components in $C_i$ and $C_j$, respectively, and when $k \geq 1$, let $\mathbf{x}^i = \langle x_1^i, \ldots, x_{q_i}^i\rangle$ be a vector with elements in $C_{l_i}$ for $1 \leq i \leq k$. Denote also the $(k+2)$-tuple $V_k = \langle \mathbf{x}^0, \mathbf{x}^1, \ldots, \mathbf{x}^{k+1}\rangle$.

For distinct $i, j \in \mathcal{G}$, and $y \in C_i$ and $z \in C_j$, define

$$\underline{\gamma}^1(n,y,z) = \max_{0 \leq k \leq N+M-2} \max_{L_k} \max_{Q_k} \max_{V_k} \mathbb{P}_{(n-1,y)}(\mathbf{X}_n^{n+r(k+1)} = \langle \mathbf{x}^0, \ldots, \mathbf{x}^{k+1}, z\rangle),$$

where the concatenated vector $\langle \mathbf{x}^0, \ldots, \mathbf{x}^{k+1}, z\rangle = \langle x_1^0, \ldots, x_{q_{k+1}}^{k+1}, z\rangle$ is of length at most $E_0(N, M) + 1$. Here, $E_0(N, M) = (\mathfrak{r}+1)(M-2) + N + 2\mathfrak{r}$ and $r(u) = \sum_{l=0}^u q_l$ for $0 \leq u \leq k+1$.

Also define

(2.10) $$\underline{\gamma}^1(n, (i,j)) = \inf_{y \in C_i, z \in C_j} \underline{\gamma}^1(n, y, z).$$

Finally, define

$$\mathcal{T}_1(i,j) = \liminf \frac{1}{n} \log \underline{\gamma}^1(n, (i,j)).$$

We now interpret the objects $\underline{\gamma}^1(n,y,z)$, $\underline{\gamma}^1(n,(i,j))$ and $\mathcal{T}_1(i,j)$. As with the routing cost $\mathcal{T}_0$, $L_k$ is an ordered list of sets to visit on the way from point $y$ to point $z$. More specifically here, $Q_k$ lists the $O(\mathfrak{r})$ number of steps taken in each visited set and $V_k$ details on which states this travel is made. Here, $\mathfrak{r}$ is chosen since all movement in a given irreducible $C_i \subset \Sigma$ is possible in at most $\mathfrak{r} = |\Sigma|$ steps. Then $\underline{\gamma}^1(n,y,z)$ is the largest probability of movement from $y$ to $z$ within the constraints of $O(\mathfrak{r})$ travel among distinct sets. Also, $\underline{\gamma}^1(n,(i,j))$ is the smallest such chance of moving from $C_i$ to $C_j$, and $\mathcal{T}_1(i,j)$ is the asymptotic exponential rate of this quantity.

**3. Results.** We now come to the main results for processes $\mathbb{P}_\pi \in \mathbb{A}(P)$. After a general upper bound and some lower bounds which depend on natural assumptions, we present an LDP which follows from these bounds. Some remarks on the Metropolis scheme and on the format of the article are made at the end of this section.

The upper bound statement is the following.

THEOREM 3.1. *With respect to good rate function $\mathbb{J}_{\mathcal{U}_0}$ and Borel $\Gamma \subset \mathbb{R}^d$, we have*

$$\limsup_{n \to \infty} \frac{1}{n} \log \mathbb{P}_\pi(Z_n \in \overline{\Gamma}) \leq -\inf_{z \in \overline{\Gamma}} \mathbb{J}_{\mathcal{U}_0}(z).$$

We now label conditions and assumptions to give LD lower bounds.

*Sufficient initial ergodicity.* To avoid degenerate cases, we introduce an initial ergodicity condition for $\mathbb{P}_\pi$ so that all information about $P$ is relevant. A typical situation to avoid is when $P_n = P$ for $n \geq m$, and distribution $\pi P_1 \cdots P_m$ locks the process evolution into a strict $P$-irreducible subset of $\Sigma$. To avoid lengthy technicalities and to be concrete, we impose the following assumption on the chains considered in this article. Let $n_0 = n_0(\{P_n\}) \geq 1$ be the first index $m$ so that for all $s,t \in C_i$ and $i \in \mathcal{G}$ when $p(s,t) > 0$ we have $p_n(s,t) > 0$ for $n \geq m$. Such an $n_0 < \infty$ exists since $P_n \to P$.

CONDITION SIE. There is an $n_1 \geq n_0 - 1$ such that
$$\mathbb{P}_\pi(X_{n_1} \in C_i) > 0 \qquad \text{for all } i \in \mathcal{G}.$$

A simpler condition which implies Condition SIE is the following.

CONDITION SIE-1. Let $n_0 = 1$ and let $\pi(C_i) > 0$ for all $i \in \mathcal{G}$.

We say that a distribution $\pi$ is SIE-1 *positive* if $\pi(C_i) > 0$ for all $i \in \mathcal{G}$. A trivial condition for SIE-1 positivity is when $\pi$ is positive [e.g., when $\pi(x) > 0$ for all $x \in \Sigma$].

*Assumptions* A, B *and* C. We now state three assumptions on the regularity of the asymptotic approach $P_n \to P$.

ASSUMPTION A. Suppose $\upsilon(i,j) = \tau(i,j)$ for all distinct $1 \leq i,j \leq N + M$.

ASSUMPTION B. Suppose for all distinct $1 \leq i,j \leq N + M$ there exists an element $a = a(i,j) \in C_i$ and a sequence $\{b_n = b_n(i,j)\} \subset C_j$ such that
$$\tau(i,j) = \lim_{n \to \infty} \frac{1}{n} \log p_n(a, b_n).$$
In other words, $\tau(i,j)$ is achieved on a fixed departing point $a \in C_i$.

ASSUMPTION C. Define $P^*(i) = \{p^*(s,t) : s,t \in C_i\}$ by
$$p^*(s,t) = \begin{cases} p(s,t), & \text{when } p(s,t) > 0, \\ 1, & \text{when } \liminf(1/n)\log p_n(s,t) = 0 \text{ and } p(s,t) = 0, \\ 0, & \text{otherwise.} \end{cases}$$
Suppose that $P^*(i)$ is primitive for $i \in \mathcal{G}$.

In words, Assumption A specifies that the maximal connection probabilities in the $(1/n)\log$ sense have limits. Assumption B states that $\tau(i,j)$ can be achieved in a systematic manner. Assumption C ensures there is "primitivity" in the system and covers the case when $P$ is periodic but the approach $P_n$ is slow enough to give a sense of primitivity. We now list some easy sufficient conditions to verify these assumptions.

PROPOSITION 3.1.

LIM. *Assumptions* A *and* B *hold if, for distinct* $1 \leq i, j \leq N + M$ *and each pair* $x \in C_i$ *and* $y \in C_j$, *we have* $\lim_{n \to \infty}(1/n) \log p_n(x,y)$ *exists.*

PRM. *Assumption* C *holds when* $\{P(i) : i \in \mathcal{G}\}$ *are primitive.*

We now come to lower bound statements for the process that obeys Condition SIE, the first of which holds in general and the second of which holds under Assumption B or C.

THEOREM 3.2. *Let* $\mathbb{P}_\pi$ *satisfy Condition* SIE.

(i) *Then with respect to good rate function* $\mathbb{J}_{\mathcal{T}_1}$ *and Borel* $\Gamma \subset \mathbb{R}^d$, *we have*

$$- \inf_{z \in \Gamma^o} \mathbb{J}_{\mathcal{T}_1}(z) \leq \liminf_{n \to \infty} \frac{1}{n} \log \mathbb{P}_\pi(Z_n \in \Gamma^o).$$

(ii) *In addition, when either Assumption* B *or* C *holds, we have with respect to good rate function* $\mathbb{J}_{\mathcal{T}_0}$ *that*

$$- \inf_{z \in \Gamma^o} \mathbb{J}_{\mathcal{T}_0}(z) \leq \liminf_{n \to \infty} \frac{1}{n} \log \mathbb{P}_\pi(Z_n \in \Gamma^o).$$

We note in the case $M = 1$ (i.e., when $P$ possesses exactly one irreducible recurrent stochastic set and possibly some degenerate transient states) that Theorems 3.1 and 3.2 already give an LDP with rate function $\mathbb{J}_{\mathcal{T}_1} = \mathbb{J}_{\mathcal{U}_0} = \mathbb{I}_1$. In particular, in this case, the large deviation behavior under $\mathbb{P}_\pi$ is independent of the approach $P_n \to P$.

However, in the general situation when $M \geq 2$, the lower and upper bounds may be different. In fact, there are nonhomogeneous processes $\mathbb{P}_\pi$ for which the lower and upper rate function bounds in Theorems 3.1 and 3.2(i) differ and are achieved so that the result is sharp in a certain sense (e.g., the example in Section 12.2).

Also, we remark that the two lower bounds in Theorem 3.2 may differ when there is some periodicity in the system and the maximal connection weight sequence is not regular. In this case, the process may not be allowed to visit freely various states because certain cyclic patterns may be in force. Therefore, the asymptotic routing costs in this general case should be larger than under Assumption B or C when some regularity is imposed on connection probabilities or when a form of primitivity is present; hence, the use of $\mathcal{T}_1$ instead of $\mathcal{T}_0$ in the lower estimates. See Section 12.3 for an explicit process where lower bounds do not respect $\mathcal{T}_0$.

It is natural now to ask when the lower and upper bounds match in the previous results so that a large deviation principle holds. For $z \in \mathbb{R}^d$, let

$$\mathbb{J}(z) = \mathbb{J}_{\mathcal{U}_0}(z).$$

Under Assumption A, costs $\mathcal{T}_0 = \mathcal{U}_0$ and so the following is a direct corollary of Theorems 3.1 and 3.2.

THEOREM 3.3. *Suppose $\mathbb{P}_\pi$ satisfies Condition SIE and Assumption A, and also either Assumption B or C. Then, with respect to good rate function $\mathbb{J}$ and Borel sets $\Gamma \subset \mathbb{R}^d$, we have the LDP*

$$-\inf_{z \in \Gamma^o} \mathbb{J}(z) \leq \liminf_{n \to \infty} \frac{1}{n} \log \mathbb{P}_\pi(Z_n \in \Gamma^o)$$

$$\leq \limsup_{n \to \infty} \frac{1}{n} \log \mathbb{P}_\pi(Z_n \in \overline{\Gamma})$$

$$\leq -\inf_{z \in \overline{\Gamma}} \mathbb{J}(z).$$

Hence, by Proposition 3.1, when all limits exist (LIM; in particular, e.g., in the time-homogeneous case, $P_n \equiv P$) or when Assumption A holds and there is no periodicity (PRM), the LDP is available. Also note that by taking $f(x) = \langle \mathbb{1}_1(x), \mathbb{1}_2(x), \ldots, \mathbb{1}_{\mathfrak{r}}(x) \rangle$, Theorem 3.3 gives the LDP for the empirical measure and so is a form of Sanov's theorem for these nonhomogeneous chains.

We remark that it may be tempting to think Assumption A by itself may be sufficient for an LDP, but it turns out there are processes which satisfy Condition SIE and Assumption A but neither B nor C for which the LDP cannot hold (e.g., the example in Section 12.3). On the other hand, we note that Assumption A is not even necessary for an LDP, for instance, with respect to chains where $P_n$ alternates between two alternatives (cf. Section 12.1). So although Theorem 3.3 is broad in a sense, more work is required to identify necessary and sufficient conditions for an LDP.

We now comment on the three types of LD behaviors mentioned in the Introduction which follow from Theorem 3.3. These are (1) homogeneous, (2) trivial and (3) intermediate behaviors for which easy sufficient (but not necessary) conditions are given below.

COROLLARY 3.1. *Let Condition SIE, and Assumption A and either Assumption B or C hold. Let also $N + M \geq 2$.*

1. *Suppose $v(i,j) = -\infty$ when $\limsup t(n, (i,j)) = 0$ for distinct $1 \leq i, j \leq N + M$. Then $\mathbb{J}$ is also the rate function for the time-homogeneous chain run under $P$ (because the routing costs are the same as if $P_n \equiv P$).*
2. *Suppose $|\mathcal{M}| \geq 2$ and $\mathcal{U}_0(i,j) = 0$ for all distinct $i, j \in \mathcal{M}$. Then $\mathbb{J}$ vanishes on the convex hull of $\bigcup_{i \in \mathcal{M}} \{z : \mathbb{I}_i(z) = 0\}$ and so is in a sense trivial.*
3. *Suppose $|\mathcal{M}| \geq 2$ and $\mathcal{U}_0(i,j) \in (-\infty, 0)$ for all distinct $i, j \in \mathcal{M}$. Then $\mathbb{J}$ differs from the rate function for the time-homogeneous chain run with $P$ and also involves nontrivially the convergence speed of $P_n$ to $P$ in terms of routing costs.*

We now briefly comment on application to the Metropolis algorithm. Note

$$1 - \sum_{j \neq i} p_n(i,j) = g(i,i) + \sum_{j \neq i} g(i,j)[1 - \exp(-\beta_n(H(j) - H(i))_+)].$$

Also, as $\beta_n \to \infty$, we have $\lim_{n \to \infty} g(i,j) \exp\{-\beta_n(H(j) - H(i))_+\} = g(i,j)\mathbb{1}_{[H(j) \leq H(i)]}$. Therefore, the limit matrix $P$ is formed in terms of entries

$$\lim_n p_n(i,j) = \begin{cases} g(i,j)\mathbb{1}_{[H(j) \leq H(i)]}, & \text{if } i \neq j, \\ g(i,i) + \sum_{j \neq i} g(i,j)\mathbb{1}_{[H(j) > H(i)]}, & \text{if } i = j. \end{cases}$$

We now decompose $P$ into components $\mathcal{D}$, $\mathcal{N}$ and $\mathcal{M}$. First, note that a state $x \in \Sigma$ belongs to the "level" set

$$C_x = \{x\} \cup \left\{ y : \exists \text{ path } x = x_0, \ldots, x_n = y, \text{ where } \prod_{i=0}^{n-1} g(x_i, x_{i+1}) > 0 \right.$$

$$\left. \text{and } H(x_i) = H(x) \text{ for } 1 \leq i \leq n \right\},$$

which corresponds to one of three types, $\mathcal{D}$, $\mathcal{N}$ or $\mathcal{M}$.

In particular, $C_x$ is a stochastic set that corresponds to $\mathcal{M}$ exactly when $H(x) = \min\{H(y) : g(x,y) > 0\}$ is a local minimum. Also, $C_x$ is a nondegenerate transient set exactly when $H(x)$ is not a local minimum and either $g(x,x) > 0$ or $g(x,y) > 0$, where $H(y) = H(x)$. Additionally, $C_x$ is a degenerate singleton exactly when $H(x)$ is not a local minimum, $g(x,x) = 0$, and when $g(x,y) > 0$ we have $H(y) \neq H(x)$.

We now discuss the rate of convergence $P_n \to P$. Observe for distinct $1 \leq i, j \leq N + M$, and $x \in C_i$ and $y \in C_j$ that

$$\limsup \frac{1}{n} \log p_n(x,y) = \begin{cases} -(H(y) - H(x))_+ \limsup(\beta_n/n), & \text{if } g(x,y) > 0, \\ -\infty, & \text{if } g(x,y) = 0, \end{cases}$$

with analogous expressions for $\liminf(1/n) \log p_n(x,y)$. Hence LIM holds when $\beta = \lim \beta_n/n$ exists. Also, we remark that when $g(x,x) > 0$ for $x \in \Sigma$, there are no degenerate transient states, so all $P$ submatrices are primitive and PRM holds. In addition, given irreducibility of $g$, Condition SIE is satisfied with respect to any initial distribution $\pi$.

Therefore, by Corollary 3.1, as routing costs are computed with respect to different level sets, the three types of LD behavior follow when the limit $\beta$ exists and there is more than one local minimum. Namely, trivial, intermediate or homogeneous behaviors occur when $\beta = 0$, $\beta \in (0, \infty)$ or $\beta = \infty$.

Finally, we give a concrete example with respect to a simple geometrically cooling Metropolis chain where $\beta = 1$. Let $H$ be defined on $\Sigma = \{1, 2, \ldots, 9\}$ in terms of its graph (Figure 1) and let $f(x) = H(x)$, so that $Z_n$ is the average $H$ value seen by the chain. Typically, for large $n$, these values $Z_n$ will be near an $H$-local minimum average.

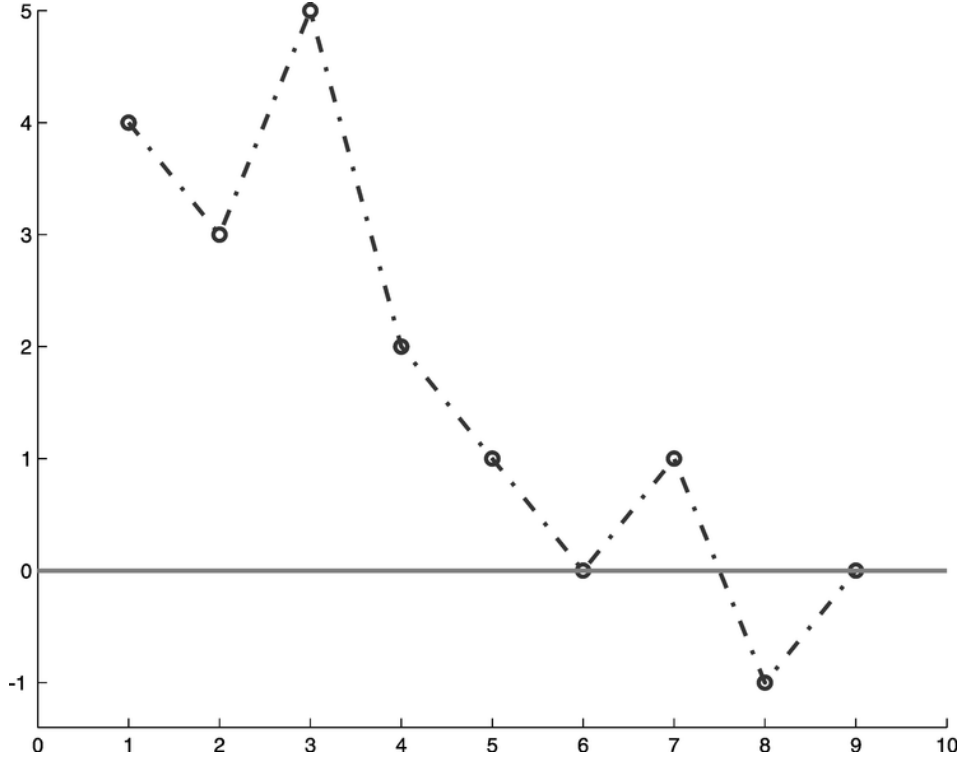

Fig. 1. *Graph of $H$.*

Let the kernel $g$ be a random walk so that $g(i, i+1) = 1/2$ for $i = 2, 6, 7, 8$, $g(i+1, i) = 1/2$ for $i = 1, 2, 6, 7$, and $g(1,2) = 1$, $g(9,8) = 1$, $g(3,4) = 1/2$, $g(4,3) = (1-a)/2$, $g(4,4) = a$, $g(4,5) = (1-a)/2$, $g(5,4) = (1-b)/2$, $g(5,5) = b$ and $g(5,6) = (1-b)/2$ with $0 < a, b < 1$. Then states $\{2\}, \{6\}, \{8\}$ are distinct local minima, $\{4\}, \{5\}$ are nondegenerate transient singletons and the remaining states are degenerate transient.

The routing costs satisfy, for distinct sets,

$$\mathcal{U}_0(\{i\}, \{j\}) = \begin{cases} -\sum_{l=i}^{j-1} (H(l) - H(l+1))_+, & \text{for } i < j, \\ -\sum_{l=j}^{i-1} (H(l+1) - H(l))_+, & \text{for } i > j. \end{cases}$$

Also, the rate functions that correspond to local minima $2, 6$ and $8$ are degenerate, and equal $\infty \cdot \mathbb{1}_{H(2)}(y)$, $\infty \cdot \mathbb{1}_{H(6)}$ and $\infty \cdot \mathbb{1}_{H(8)}$, respectively.

For the nondegenerate transient states 4 and 5, we have

$$I_{\{4\}}(y) = \begin{cases} -\log \dfrac{1+a}{2}, & \text{for } y = H(4), \\ \infty, & \text{otherwise,} \end{cases}$$

and

$$I_{\{5\}}(y) = \begin{cases} -\log \dfrac{1+b}{2}, & \text{for } y = H(5), \\ \infty, & \text{otherwise.} \end{cases}$$

When $-\log(1+a)/2 = 1/3$ and $-\log(1+b)/2 = 2/3$, we compute, by analyzing the not-too-large number of possibilities, the nonconvex rate function

$$\mathbb{J}(z) = \begin{cases} \infty, & \text{for } z < -1 \text{ and } z > 3, \\ 4z/9 + 4/9, & \text{for } -1 \leq z \leq -2/11, \\ -2z, & \text{for } -2/11 \leq z \leq 0, \\ z/6, & \text{for } 0 \leq z \leq 2, \\ 5z/3 - 3, & \text{for } 2 \leq z \leq 12/5, \\ -5z/3 + 5, & \text{for } 12/5 \leq z \leq 3. \end{cases}$$

Not surprisingly, $\mathbb{J}$ vanishes at local minima and is largest near $z \sim 2^+$ (excluding infinite costs), with exact value $z = 12/5$ found from computation. The $\mathbb{J}$ calculation (see Figure 2) also gives optimal scenarios under which $Z_n \sim z$; these include, for $-1 \leq z \leq -2/11$ that the average $Z_n$ is a convex combination of rest stays initially on $\{4\}$ and then at $\{8\}$; for $-2/11 \leq z \leq 0$, at $\{8\}$, then $\{6\}$; for $0 \leq z \leq 2$, at $\{4\}$, then $\{6\}$; for $2 \leq z \leq 12/5$, at $\{2\}$, then $\{4\}$; for $12/5 \leq z \leq 3$, at $\{6\}$, then $\{2\}$.

We now discuss the plan of the paper. In the next section, we outline the proof structure of the main theorems. After supplying proofs of stated results in the outline in Sections 5–11, we give the three examples in Section 12 commented on earlier. Finally, in the Appendix some technical proofs are collected.

**4. Outline of the proofs of the main theorems.** Consider a process $\mathbb{P}_\pi \in \mathbb{A}(P)$ and a function $f : \Sigma \to \mathbb{R}^d$. We first observe that $\mathbb{J}_{\mathcal{U}_0}$, $\mathbb{J}_{\mathcal{T}_0}$ and $\mathbb{J}_{\mathcal{T}_1}$ are all good rate functions from the following proposition, which is proved in the Appendix.

PROPOSITION 4.1. *For a nonpositive cost $U$, the function $\mathbb{J}_U$ is a good rate function and the domain of finiteness $Q_{\mathbb{J}_U} \subset \mathbb{K}$.*

In the following discussion, we say that the path $\mathbf{X}_n$ enters or visits a subset $C \subset \Sigma$ when $X_i \in C$ for some $1 \leq i \leq n$. We now outline the proofs of Theorems 3.1 and 3.2.

4.1. *Upper bounds*: *proof of Theorem* 3.1. The proof follows by first a *surgery of paths* estimate, then a *homogeneous rest cost* comparison, a *coarse graining cost* estimate and finally a limit relationship on a perturbed rate function. Let $\Gamma \subset \mathbb{R}^d$ be a Borel set.

*Surgery of paths estimate.* The first step is to overestimate $\mathbb{P}_\pi$ by another measure $\hat{\mu}_{\pi,\varepsilon_1,\varepsilon_2}$ which allows more movement in terms of parameters $\varepsilon_1, \varepsilon_2 > 0$. However, we restrict the process to those paths which make at most one long sojourn to each of the sets $\{C_i : 1 \leq i \leq N+M\}$, but connect among them in short visits.

Before getting to the first bound, the following technical monotonicity lemma, proved in the Appendix, is needed.

LEMMA 4.1. *Let $\delta \in [0,1]$ and let $\{t_n\} \subset [0,1]$ be a sequence which converges to $\delta$. Then there exists a sequence $\{\hat{t}_n\} \subset (0,1]$ such that* (i) $t_n \leq \hat{t}_n$, (ii) $\hat{t}_n \downarrow \delta$ *monotonically and* (iii) *the limit $\lim (1/n) \log \hat{t}_n$ exists and equals*

$$\lim_{n \to \infty} \frac{1}{n} \log \hat{t}_n = \limsup_{n \to \infty} \frac{1}{n} \log t_n.$$

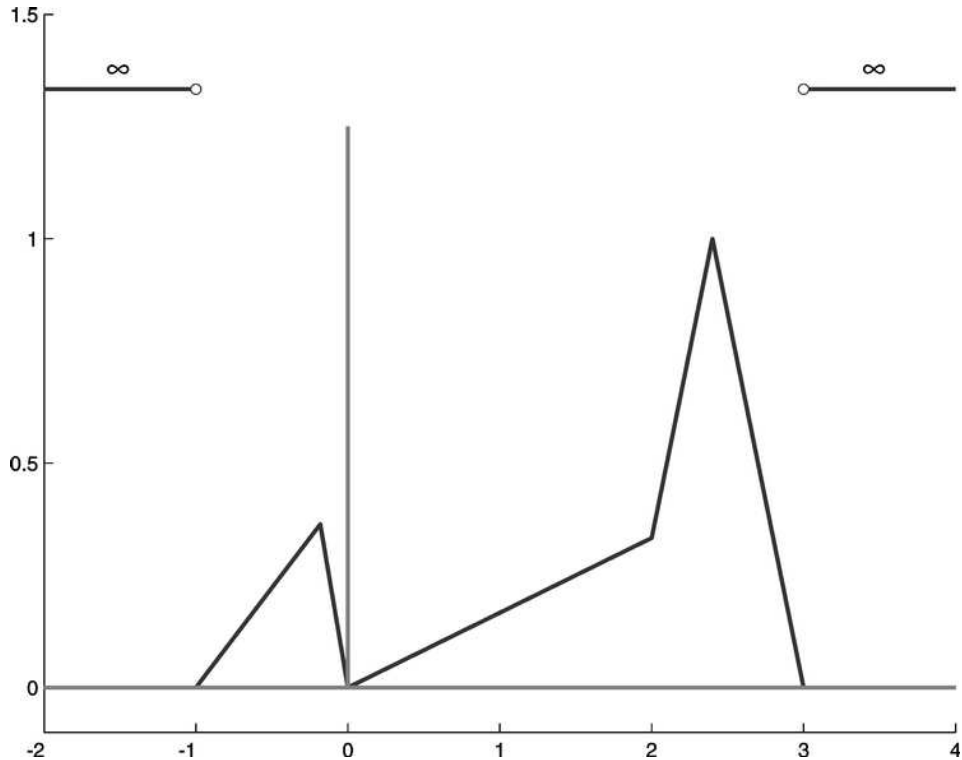

FIG. 2. *Graph of* $\mathbb{J}$.

Recall now the definition of $t(n,(i,j))$ [cf. (2.7)] and let

$\{\hat{t}(n,(i,j))\}$ be the sequence made from $\{t(n,(i,j))\}$ and Lemma 4.1.

Also, for distinct $1 \leq i, j \leq N + M$, as in the definition of $\mathcal{U}_0(i,j)$ [cf. (2.8)], let $0 \leq k \leq N + M - 2$, let $l_0 = i$ and $l_{k+1} = j$, and let $L_k = \langle l_0, l_1, \ldots, l_k, l_{k+1} \rangle$ be composed of distinct indices. Then define

$$\gamma(n,(i,j)) = \max_{0 \leq k \leq N+M-2} \max_{L_k} \prod_{s=0}^{k} \hat{t}(n+s,(l_s, l_{s+1})).$$

The term $\gamma(n,(i,j))$ bounds the largest possible transition probability between sets $C_i$ and $C_j$ in at most $N + M - 1$ steps.

We now create a certain sequence of positive transition matrices. For general $P$ and approaching sequence $\{P_n\}$, the submatrices $P(i)$ and $P_n(i)$ for $1 \leq i \leq N + M$ need not be positive. It will be helpful, however, to majorize them as follows. Let $\varepsilon \geq 0$, and let $P(i, \varepsilon) = \{p(s, t; \varepsilon) : s, t \in C_i\}$ and $P_n(i, \varepsilon) = \{p_n(s, t; \varepsilon) : s, t \in C_i\}$, where

$$p(s, t; \varepsilon) = \max\{p(s,t), \varepsilon\} \quad \text{and} \quad p_n(s, t; \varepsilon) = \max\{p_n(s,t), \varepsilon\}.$$

Define now $\widehat{\mathcal{P}}_{n,\varepsilon_1,\varepsilon_2} = \{\hat{p}_{n,\varepsilon_1,\varepsilon_2}(s,t)\}$ by

$$\hat{p}_{n,\varepsilon_1,\varepsilon_2}(s,t) = \begin{cases} \gamma(n,(i,j)), & \text{for } s \in C_i,\ t \in C_j \\ & \text{and distinct } 1 \leq i,\ j \leq N+M, \\ p_n(s,t;\varepsilon_2), & \text{for } s, t \in C_i \text{ and } i \in \mathcal{G}, \\ p_n(s,t;\varepsilon_1), & \text{for } s, t \in C_i \text{ and } i \in \mathcal{D}, \end{cases}$$

when $n \geq 2$; for $n = 1$, let $\widehat{\mathcal{P}}_{1,\varepsilon_1,\varepsilon_2}$ be the unit constant matrix, $\hat{p}_{1,\varepsilon_1,\varepsilon_2}(s,t) \equiv 1$. Form also through CON the measure $\hat{\mu}_{\pi,\varepsilon_1,\varepsilon_2}$ with respect to initial distribution $\pi$ and transition matrices $\{\widehat{\mathcal{P}}_{n,\varepsilon_1,\varepsilon_2}\}$.

PROPOSITION 4.2. *For $\varepsilon_1, \varepsilon_2 > 0$, the following upper bound holds:*

$$\limsup \frac{1}{n} \log \mathbb{P}_\pi(Z_n \in \Gamma)$$

$$\leq \limsup \frac{1}{n} \log \hat{\mu}_{\pi,\varepsilon_1,\varepsilon_2}(Z_n \in \Gamma, \mathbf{X}_n \text{ enters each } C_i \text{ at most once}).$$

The proof of this proposition is found in Section 5.

*Homogeneous rest cost comparison.* Next, we compare measure $\hat{\mu}_{\pi,\varepsilon_1,\varepsilon_2}$ with a measure $\bar{\mu}_{\pi,\varepsilon_1,\varepsilon_2}$, which replaces nonhomogeneous transitions within sets $C_i$ by limiting homogeneous transition weights.

Define, for $\varepsilon_1, \varepsilon_2 \geq 0$, $\overline{\mathcal{P}}_{n,\varepsilon_1,\varepsilon_2} = \{\bar{p}_{n,\varepsilon_1,\varepsilon_2}(s,t)\}$ by

$$\bar{p}_{n,\varepsilon_1,\varepsilon_2}(s,t) = \begin{cases} \gamma(n,(i,j)), & \text{for } s \in C_i,\ t \in C_j \\ & \text{and distinct } 1 \leq i,j \leq N+M, \\ p(s,t;\varepsilon_2), & \text{for } s,t \in C_i \text{ and } i \in \mathcal{G}, \\ \varepsilon_1, & \text{for } s,t \in C_i \text{ and } i \in \mathcal{D}, \end{cases}$$

when $n \geq 2$ and $\overline{\mathcal{P}}_{1,\varepsilon_1,\varepsilon_2} = \widehat{\mathcal{P}}_{1,\varepsilon_1,\varepsilon_2}$. Let now $\bar{\mu}_{\pi,\varepsilon_1,\varepsilon_2}$ be formed from CON and matrices $\{\overline{\mathcal{P}}_{n,\varepsilon_1,\varepsilon_2}\}$ and $\pi$.

PROPOSITION 4.3. *For $\varepsilon_1, \varepsilon_2 > 0$, we have*

$$\limsup \frac{1}{n} \log \hat{\mu}_{\pi,\varepsilon_1,\varepsilon_2}(Z_n \in \Gamma, \mathbf{X}_n \text{ enters each } C_i \text{ at most once})$$
$$\leq \limsup \frac{1}{n} \log \bar{\mu}_{\pi,\varepsilon_1,\varepsilon_2}(Z_n \in \Gamma, \mathbf{X}_n \text{ enters each } C_i \text{ at most once}).$$

The proof of this proposition is found in Section 7.

*Coarse graining estimate.* The next step is to further bound the right-hand side in Proposition 4.3 through a detailed decomposition of visit times and locations in terms of an $\varepsilon_1, \varepsilon_2$-perturbed rate $\mathbb{J}_{\mathcal{U}_0,\varepsilon_1,\varepsilon_2}$.

Observe for $1 \leq i \leq N+M$ that the submatrix $(\overline{P}_{n,\varepsilon_1,\varepsilon_2})_{C_i} = P(i,\varepsilon_1,\varepsilon_2)$ is independent of $n$ and

$$P(i, \varepsilon_1, \varepsilon_2) = \begin{cases} (\varepsilon_1), & \text{for } i \in \mathcal{D}, \\ P(i, \varepsilon_2), & \text{for } i \in \mathcal{G}. \end{cases}$$

Denote the extended rate function $\mathbb{I}_{i,\varepsilon_1,\varepsilon_2} = \mathbb{I}_{C_i,f,P(i,\varepsilon_1,\varepsilon_2)}$ and associated domain of finiteness $Q_{i,\varepsilon_1,\varepsilon_2} = Q_{C_i,f,P(i,\varepsilon_1,\varepsilon_2)}$. In fact, explicitly when $i \in \mathcal{D}$,

$$(4.1) \quad \mathbb{I}_{i,\varepsilon_1,\varepsilon_2}(x) = \begin{cases} -\log(\varepsilon_1), & \text{for } x = f(m_i), \text{ where } C_i = \{m_i\}, \\ \infty, & \text{otherwise}, \end{cases}$$

and $\mathbb{I}_{i,\varepsilon_1,\varepsilon_2}(x) = \mathbb{I}_{i,\varepsilon_2}(x) = \mathbb{I}_{C_i,f,P(i,\varepsilon_2)}$ when $i \in \mathcal{G}$.

Recall now the object $C_{\mathbf{v},U}$ near (2.6), and define for $\mathbf{v} \in \Omega_{N+M}$, $\mathbf{x} \in (\mathbb{R}^d)^{N+M}$, $\sigma \in \mathbb{S}_{N+M}$ and matrix $U = \{u(i,j) : 1 \leq i, j \leq N+M\}$, the function

$$C_{\mathbf{v},U,\varepsilon_1,\varepsilon_2}(\sigma,\mathbf{x}) = -\sum_{i=1}^{N+M-1} \left(\sum_{j=1}^{i} v_j\right) u(\sigma(i),\sigma(i+1)) + \sum_{i=1}^{N+M} v_i \mathbb{I}_{\sigma(i),\varepsilon_1,\varepsilon_2}(x_i)$$

when $N + M \geq 2$ and $C_{\mathbf{v},U,\varepsilon_1,\varepsilon_2}(\sigma,\mathbf{x}) = \mathbb{I}_{1,\varepsilon_1,\varepsilon_2}(x_1)$ when $N+M = 1$. Define also, for $z \in \mathbb{R}^d$,

$$(4.2) \quad \mathbb{J}_{U,\varepsilon_1,\varepsilon_2}(z) = \inf_{\mathbf{v} \in \Omega_{N+M}} \inf_{\mathbf{x} \in D(N+M,\mathbf{v},z)} \min_{\sigma \in S_{N+M}} C_{\mathbf{v},U,\varepsilon_1,\varepsilon_2}(\sigma,\mathbf{x}).$$

We comment that when $N = 0$ and all $P(i) > 0$ for $i \in \mathcal{G}$, that $\mathbb{J}_{U,\varepsilon_1,\varepsilon_2} = \mathbb{J}_U$ for all $\varepsilon_1, \varepsilon_2$ small, so the following result already gives the desired upper bound.

PROPOSITION 4.4. *For $\varepsilon_1, \varepsilon_2 > 0$, we have*

$$\limsup_{n \to \infty} \frac{1}{n} \log \bar{\mu}_{\pi, \varepsilon_1, \varepsilon_2}(Z_n \in \overline{\Gamma}, \mathbf{X}_n \text{ enters each } C_i \text{ at most once})$$
$$\leq -\mathbb{J}_{\mathcal{U}_0, \varepsilon_1, \varepsilon_2}(\overline{\Gamma} \cap \mathbb{K}).$$

The proof of the proposition is given in Section 8.

*Limit estimate on $\mathbb{J}_{\mathcal{U}_0, \varepsilon_1, \varepsilon_2}$.* The last step is to analyze $\mathbb{J}_{\mathcal{U}_0, \varepsilon_1, \varepsilon_2}$ as $\varepsilon_1, \varepsilon_2 \downarrow 0$ in the following proposition, which is proved in Section 10.

PROPOSITION 4.5. *We have*

$$\limsup_{\varepsilon_2 \downarrow 0} \limsup_{\varepsilon_1 \downarrow 0} -\mathbb{J}_{\mathcal{U}_0, \varepsilon_1, \varepsilon_2}(\overline{\Gamma} \cap \mathbb{K}) \leq -\mathbb{J}_{\mathcal{U}_0}(\overline{\Gamma}).$$

Now, putting together the results above gives Theorem 3.1.

4.2. *Lower bounds: proof of Theorem 3.2.* The argument is similar in structure to the upper bound. To prove part (i), a reduction is first made with respect to initial ergodicity, which can be skipped if one is willing to assume that $\mathbb{P}_\pi$ satisfies the stronger Condition SIE-1 rather than just Condition SIE. Then a surgery of paths estimate, a homogeneous rest cost comparison and finally a coarse graining cost estimate are given. Last, having proved part (i), the second lower bound part (ii) is argued.

Let $\Gamma \subset \mathbb{R}^d$ be a Borel set. If $\Gamma^o = \varnothing$, the bound is trivial. Otherwise, let $x_0 \in \Gamma^o$ and $\Gamma_1 = B(x_0, a) \subset \Gamma^o$ be an open ball of radius $a > 0$.

*SIE estimate.* The following estimate shows that under Condition SIE, the first few transition kernels do not contribute effectively to lower bounds and, in particular, Condition SIE may be replaced with Condition SIE-1. When $\mathbb{P}_\pi$ satisfies Condition SIE, let $P'_n = P_{n+n_1}$ for $n \geq 1$, and let $\eta(l) = \mathbb{P}_\pi(X_{n_1} = l)$ for $l \in \Sigma$. Let also $\mathbb{P}'_\eta$ be constructed with respect to $\{P'_n\}$ and distribution $\eta$. Clearly, we have $n_0(\{P'_n\}) = 1$ and $\mathbb{P}'_\eta$ satisfies Condition SIE-1.

PROPOSITION 4.6. *Let $\Gamma_2 = B(x_0, a/2)$ and suppose $\mathbb{P}_\pi$ satisfies Condition* SIE. *Then we have*

$$\liminf \frac{1}{n} \log \mathbb{P}_\pi(Z_n \in \Gamma_1) \geq \liminf \frac{1}{n} \log \mathbb{P}'_\eta(Z_n \in \Gamma_2).$$

PROOF. Note that

$$\{Z_n \in B(x_0, a)\} \supset \left\{\frac{n - n_1}{n} Z^n_{n_1+1} \in B\left(x_0, a - \frac{c_1}{n}\right)\right\},$$

where $c_1 = n_1 \|f\|$. Then

$$\begin{aligned}
\mathbb{P}_\pi(Z_n \in \Gamma_1) &\geq \mathbb{P}_\pi(((n-n_1)/n)Z^n_{n_1+1} \in B(x_0, a - c_1/n)) \\
&= \sum_{l \in \Sigma} \eta(l) \mathbb{P}_{(n_1,l)}(((n-n_1)/n)Z^n_{n_1+1} \in B(x_0, a - c_1/n)) \\
&= \mathbb{P}'_\eta(Z_{n-n_1} \in (n/(n-n_1))B(x_0, a - c_1/n)).
\end{aligned}$$

The proposition now follows by simple calculations. □

In view of the last proposition, with regard to the standard lower bound methods, we may just as well assume that $\mathbb{P}_\pi$ satisfies Condition SIE-1 if Condition SIE already holds.

*Surgery of paths estimate.* We underestimate $\mathbb{P}_\pi$ by another measure $\check{\mu}_{\pi,\varepsilon_1,\varepsilon_2}$ whose connection transitions correspond to $\mathcal{T}_1$. Slightly different from the surgery for the upper bound, the paths focused on here are those which make at most one long visit to sets $\{C_i : i \in \mathcal{G}\}$, but travel between them in short trips through all $\{C_i : 1 \leq i \leq N+M\}$.

Let $E(N, M) = (M-1)E_0(N, M)$ and recall the connecting weight $\underline{\gamma}^1(n, (i, j))$ for distinct $i, j \in \mathcal{G}$ [cf. (2.10)]. Define

$$\check{\gamma}^0(n, (i,j)) = \min_{0 \leq k \leq E(N,M)} \underline{\gamma}^1(n+k, (i,j)),$$

which picks the smallest weight in a traveling frame.

Define also $\check{\mathcal{P}}_n = \{\check{p}_n(s,t)\}$ for $n \geq 1$ by

$$\check{p}_n(s,t) = \begin{cases} \check{\gamma}^0(n, (i,j)), & \text{for all } s \in C_i, t \in C_j \\ & \text{and distinct } i, j \in \mathcal{G}, \\ p_n(s,t), & \text{for } s, t \in C_i \text{ and } i \in \mathcal{G} \\ & \text{or } s \in C_i, t \in C_j \text{ when } i \text{ or } j \in \mathcal{D}. \end{cases}$$

Let $\check{\mu}_\pi$ be made through CON with $\{\check{\mathcal{P}}_n\}$ and $\pi$.

In addition, for convenience, let

$$G_n = \{\mathbf{X}_n \text{ enters only } \{C_i : i \in \mathcal{G}\} \text{ with at most one visit to each set}\}.$$

PROPOSITION 4.7. *Let $\Gamma_3 = B(x_0, a/4)$ and suppose $\mathbb{P}_\pi$ satisfies Condition* SIE-1. *Then*

$$\liminf \frac{1}{n} \log \mathbb{P}_\pi(Z_n(f) \in \Gamma_2) \geq \liminf \frac{1}{n} \log \check{\mu}_\pi(Z_n(f) \in \Gamma_3, G_n).$$

The proof is given in Section 6.

*Homogeneous rest cost comparision.* As before, we compare $\check{\mu}_\pi$ with a measure $\underline{\mu}_\pi$, which replaces nonhomogeneous transitions within sets $C_i$ with limiting homogeneous transition weights.

Define $\underline{\mathcal{P}}_n = \{\underline{p}_n(s,t)\}$ for $n \geq 1$ by

$$\underline{p}_n(s,t) = \begin{cases} \check{\gamma}^0(n,(i,j)), & \text{for all } s \in C_i,\, t \in C_j \text{ and distinct } i,j \in \mathcal{G}, \\ p(s,t), & \text{for } s,t \in C_i \text{ and } i \in \mathcal{G}, \\ 0, & \text{otherwise.} \end{cases}$$

Correspondingly, define $\underline{\mu}_\pi$ through CON with $\{\underline{\mathcal{P}}_n\}$ and initial distribution $\pi$.

PROPOSITION 4.8. *Suppose $n_0(\{P_n\}) = 1$. Then we have*

$$\liminf \frac{1}{n} \log \check{\mu}_\pi(Z_n \in \Gamma_3, G_n) \geq \liminf \frac{1}{n} \log \underline{\mu}_\pi(Z_n \in \Gamma_3, G_n).$$

The proof is given in Section 7.

*Coarse graining estimate.* Again, we bound the right-hand side above through a decomposition of visit times and locations.

PROPOSITION 4.9. *Let $\pi$ be SIE-1-positive. Then*

$$\liminf \frac{1}{n} \log \underline{\mu}_\pi(Z_n \in \Gamma_3, G_n) \geq -\mathbb{J}_{\mathcal{T}_1}(\Gamma_3).$$

The proof is given in Section 9.

Finally, whereas $x_0 \in \Gamma^o$ is arbitrary, we have that

$$\liminf_{n \to \infty} \frac{1}{n} \log \mathbb{P}_\pi(Z_n \in \Gamma^o) \geq - \inf_{z \in \Gamma^o} \mathbb{J}_{\mathcal{T}_1}(z)$$

and so part (i) is proved.

PROOF OF THEOREM 3.2(ii). The following cost bound, proved in Section 11, is the key step.

PROPOSITION 4.10. *We have under Assumptions B or C that $\mathcal{T}_1 \geq \mathcal{T}_0$ and so $\mathbb{J}_{\mathcal{T}_1} \leq \mathbb{J}_{\mathcal{T}_0}$.*

Therefore, given the lower bound in part (i), the second part follows directly. □

**5. Path surgery upper bound.** The strategy of Proposition 4.2 is to compare the probability of a path which moves many times between sets with that of a respective rearranged path with fewer sojourns. To make estimates we need a few more definitions.

Let $\underline{t}(n)$ be the largest entry which connects upward with respect to the ordering of the sets $\{C_i\}$ in the canonical decomposition of $P$:

$$\underline{t}(n) = \max_{1 \leq j < i \leq N+M} \hat{t}(n,(i,j)).$$

Observe that as movement up the tree is impossible in the limit or, more precisely, as $T_n(i,j)$ vanishes for $1 \leq j < i \leq N+M$, we have $\underline{t}(n) \to 0$ as $n \to \infty$.

Define also for $\varepsilon_1, \varepsilon_2 \geq 0$, the matrix $\widetilde{P}_{n,\varepsilon_1,\varepsilon_2} = \{\tilde{p}_{n,\varepsilon_1,\varepsilon_2}(s,t)\}$ by

$$\tilde{p}_{n,\varepsilon_1,\varepsilon_2}(s,t) = \begin{cases} \hat{t}(n,(i,j)), & \text{for } s \in C_i,\ t \in C_j \\ & \text{and distinct } 1 \leq i,\ j \leq N+M, \\ p_n(s,t;\varepsilon_2), & \text{for } s,t \in C_i \text{ and } i \in \mathcal{G}, \\ p_n(s,t;\varepsilon_1), & \text{for } s,t \in C_i \text{ and } i \in \mathcal{D}, \end{cases}$$

for $n \geq 1$. Form now through CON the measure $\nu_{\pi,\varepsilon_1,\varepsilon_2}$ with respect to initial distribution $\pi$ and transition matrices $\{\widetilde{P}_{n,\varepsilon_1,\varepsilon_2}\}$.

Let also $\tilde{p} = \min\{\varepsilon_1, \varepsilon_2\}$ and observe that $\tilde{p}$ is less than the minimum transition probability within subblocks:

$$\tilde{p} \leq \min_{1 \leq l \leq N+M} \min_{s,t \in C_l} \tilde{p}_{n,\varepsilon_1,\varepsilon_2}(s,t).$$

We now describe a procedure to cut paths into resting and traveling parts, which then are rearranged through a rearrangement map. Let $\mathbf{x}_n = \langle x_1, \ldots, x_n \rangle \in \Sigma^n$ be a path of length $n \geq 2$. We say that $\mathbf{x}_n$ possesses a "switch" at time $1 \leq i \leq n-1$ if $x_i \in C_j$ and $x_{i+1} \in C_k$ for $j \neq k$. For a path $\mathbf{x}_n$ which switches $l \geq 1$ times, let $g_k(\mathbf{x}_n)$ be the time of the $k$th switch, where $1 \leq k \leq l$. Set also $g_0(\mathbf{x}_n) = 0$ and $g_{l+1}(\mathbf{x}_n) = n$.

Define now, for $1 \leq k \leq l$, the path segments between switch times: $J_k(\mathbf{x}_n) = \langle x_{g_{k-1}(\mathbf{x}_n)+1}, \ldots, x_{g_k(\mathbf{x}_n)} \rangle$, and the remainder $J_{l+1}(\mathbf{x}_n) = \langle x_{g_l(\mathbf{x}_n)+1}, \ldots, x_n \rangle$. Define also that $J_{k,2}(\mathbf{x}_n) = \langle x_{g_{k-1}(\mathbf{x}_n)+2}, \ldots, x_{g_k(\mathbf{x}_n)} \rangle$ when $g_k(\mathbf{x}_n) \geq g_{k-1}(\mathbf{x}_n) + 2$.

In addition, let $C_{i_k}$ be the subset in which path $J_k$ lies for $1 \leq k \leq l+1$ and let $\mathcal{C}_l = \mathcal{C}_l(\mathbf{x}_n) = \langle C_{i_1}, \ldots, C_{i_{l+1}} \rangle$ be the sequence of subsets visited, given in the order of visitation. Also, let $\|\mathcal{C}_l\|$ be the number of distinct elements in $\mathcal{C}_l$. We say $\mathbf{x}_n$ has no repeat visits if the sequence $\mathcal{C}_l$ contains no repetitions.

For $0 \leq k \leq n-1$ and $1 \leq j \leq N+M$, define the sets

$$A_n(k) = \{\mathbf{x}_n : \mathbf{x}_n \text{ switches } k \text{ times}\}$$

and

$$A'_n(j) = \{\mathbf{x}_n : \mathbf{x}_n \text{ switches } j \text{ times, with no repeat visits}\}.$$

When there are at least two sets, $N + M \geq 2$, we define the map
$$\sigma_l : A_n(l) \to \bigcup_{j=1}^{\min\{N+M-1,l\}} A'_n(j)$$
for $l \geq 1$, in the following steps.

1. Let $\mathbf{x}_n \in A_n(l)$. Let $s_{\|\mathcal{C}_l\|} = l+1$ and $s_{\|\mathcal{C}_l\|-1} = l$. Inductively define, for $k < \|\mathcal{C}_l\|$,
$$s_k = \max\left\{j : C_{i_j} \notin \left\{C_{i_{s_{k+1}}}, C_{i_{s_{k+2}}}, \ldots, C_{i_{s_{\|\mathcal{C}_l\|}}}\right\}\right\}.$$
In words, $C_{i_{s_{\|\mathcal{C}_l\|}}}, \ldots, C_{i_{s_1}}$ are the $\|\mathcal{C}_l\|$ distinct subsets visited in reverse order starting from the last state of $\mathbf{x}_n$.

2. For $1 \leq k \leq \|\mathcal{C}_l\|$, let $J_{\alpha_1^k}, \ldots, J_{\alpha_{d_k}^k}$, where $\alpha_1^k < \cdots < \alpha_{d_k}^k = s_k$, be the $d_k \geq 1$ paths which lie in $C_{i_{s_k}}$.

3. Define
$$\sigma_l(\mathbf{x}_n) = \left\langle J_{\alpha_1^1}, \ldots, J_{\alpha_{d_1}^1}, \ldots, J_{\alpha_1^{\|\mathcal{C}_l\|}}, \ldots, J_{\alpha_{d_{\|\mathcal{C}_l\|}}^{\|\mathcal{C}_l\|}} \right\rangle.$$

In words, $\sigma_l$ rearranges the paths that correspond to distinct subsets so that the reverse visiting order is preserved. We comment that the last path $J_{\alpha_{d_{l+1}}^{l+1}}$ is preserved under $\sigma_l$ and that $\sigma_1$ is the identity map.

EXAMPLE 1. Suppose $N + M = 8$ and $\mathbf{x}_n \in A_n(25)$, where
$$\mathcal{C}_{25} = \langle C_8, C_6, C_8, C_7, C_5, C_7, C_6, C_5, C_6, C_4, C_2, C_4,$$
$$C_3, C_1, C_3, C_1, C_2, C_1, C_6, C_7, C_5, C_4, C_2, C_5, C_2, C_4 \rangle.$$
Here, $\|\mathcal{C}\| = 8$, $s_1 = 3$, $s_2 = 15$, $s_3 = 18$, $s_4 = 19$, $s_5 = 20$, $s_6 = 24$, $s_7 = 25$ and $s_8 = 26$. Then
$$\langle C_{i_{s_1}}, C_{i_{s_2}}, C_{i_{s_3}}, C_{i_{s_4}}, C_{i_{s_5}}, C_{i_{s_6}}, C_{i_{s_7}}, C_{i_{s_8}} \rangle = \langle C_8, C_3, C_1, C_6, C_7, C_5, C_2, C_4 \rangle$$
and
$$\sigma_{25}(\mathbf{x}_n) = \langle J_1, J_3, J_{13}, J_{15}, J_{14}, J_{16}, J_{18}, J_2, J_7, J_9,$$
$$J_{19}, J_4, J_6, J_{20}, J_5, J_8, J_{21}, J_{24} J_{11}, J_{17} J_{23}, J_{25}, J_{10}, J_{12}, J_{22}, J_{26} \rangle.$$

Finally, we recall at this point useful versions of the "union of events" bound.

LEMMA 5.1. *Let $N \geq 1$ and let $\{a_n^i : i, n \geq 1\}$ be an array of nonnegative numbers. We have then*
$$\limsup_{n \to \infty} \frac{1}{n} \log \sum_{i=1}^N a_n^i = \max_{1 \leq i \leq N} \limsup_{n \to \infty} \frac{1}{n} \log a_n^i$$

*and*

$$\liminf_{n\to\infty} \frac{1}{n}\log \sum_{i=1}^{N} a_n^i = \liminf \max_{1\le i\le N} \frac{1}{n}\log a_n^i \ge \max_{1\le i\le N} \liminf_{n\to\infty} \frac{1}{n}\log a_n^i.$$

*In addition, let $\alpha \ge 1$ be an integer and let $\{\beta(n)\}$ be a sequence where $\beta(n) \le n^\alpha$ for $n \ge 1$. Then*

$$\limsup_{n\to\infty} \frac{1}{n}\log \sum_{i=1}^{\beta(n)} a_n^i = \limsup_{n\to\infty} \max_{1\le i\le \beta(n)} \frac{1}{n}\log a_n^i$$

*with the same equality when $\liminf$ replaces $\limsup$.*

See [10], Lemma 1.2.15, for the "lim sup" proof. The other statements follow similarly.

PROOF OF PROPOSITION 4.2. As $P_n \le \widetilde{P}_{n,\varepsilon_1,\varepsilon_2}$ elementwise, we have

$$\mathbb{P}_\pi(Z_n \in \overline{\Gamma}) \le \nu_{\pi,\varepsilon_1,\varepsilon_2}(Z_n \in \overline{\Gamma}).$$

Now consider the case $N + M = 1$ when $P$ corresponds to one irreducible set $C_1 = \Sigma$. Trivially in this case $\mathbf{X}_n$ does not leave $C_1$, so more than one switch is impossible. Therefore, the upper bound statement holds immediately.

We now assume that $N + M \ge 2$. By Lemma 5.1,

$$(5.1) \quad \limsup \frac{1}{n}\log \mathbb{P}_\pi(Z_n \in \Gamma) \le \max_{\substack{x_0 \in \Sigma \\ \pi(x_0) > 0}} \limsup \frac{1}{n}\log \nu_{x_0,\varepsilon_1,\varepsilon_2}(Z_n \in \overline{\Gamma}).$$

Hence, it suffices to focus on $\nu_{x_0,\varepsilon_1,\varepsilon_2}$ for a given $x_0 \in \Sigma$ such that $\pi(x_0) > 0$.

The main idea exploited now is that for a realization $\mathbf{X}_n$ which switches between sets $\{C_i\}$ many times there will be guaranteed a large number of these switches "up the tree" between sets $C_i$ and $C_j$ for $i > j$ whose chance is small, and so such paths are unlikely. For notational simplicity, we now suppress $\varepsilon_1$ and $\varepsilon_2$ subscripts.

STEP 1. Decompose according to the number of switches:

$$(5.2) \qquad \nu_{x_0}(Z_n \in \Gamma) = \sum_{i=0}^{n-1} \nu_{x_0}(Z_n \in \Gamma, A_n(i)).$$

STEP 2. Let $l \ge 1$ and let $\mathbf{x}_n \in \{Z_n \in \Gamma\} \cap A_n(l)$. Let also $\mathbf{y}_n \in \sigma_l^{-1}(\sigma_l(\mathbf{x}_n))$, that is, $\mathbf{y}_n$ is a path with $l$ switches which rearranges to $\sigma_l(\mathbf{x}_n)$. As $\mathbf{y}_n =$ $\langle J_1(\mathbf{y}_n), \ldots, J_{l+1}(\mathbf{y}_n)\rangle$, where $J_k(\mathbf{y}_n)$ is a path in $C_{i_k}$ for $1 \le k \le l+1$, we have

$$\nu_{x_0}(\mathbf{X}_n = \mathbf{y}_n)$$
(5.3)
$$= \nu_{x_0}(\mathbf{X}_1^{g_1} = J_1)$$
$$\times \prod_{\substack{1 \le k \le l \\ g_k \ge g_{k-1}+2}} \nu_{(g_k+1, y_{g_k+1})}(\mathbf{X}_{g_k+2}^{g_{k+1}} = J_{k+1,2}) \prod_{k=1}^{l} \hat{t}(g_k+1, (i_k, i_{k+1})),$$

where $g_k = g_k(\mathbf{y}_n)$ and $J_{k+1,2} = J_{k+1,2}(\mathbf{y}_n)$ (defined above) are shortened for clarity.

We now bound the right-hand side of (5.3) by

(5.4)
$$(\tilde{p}_1(x_0, y_1)/1)\hat{\mu}_{x_0}(\mathbf{X}_n = \sigma_l(\mathbf{x}_n)) \prod_{k=1}^{l} \hat{t}(g_k(\mathbf{y}_n)+1, (i_k, i_{k+1}))$$
$$\times \prod_{k=1}^{\|\mathcal{C}_l\|-1} \gamma^{-1}(g_k(\sigma_l(\mathbf{x}_n))+1, (i_{s_k}, i_{s_{k+1}})) \cdot (1/\tilde{p})^{l-(\|\mathcal{C}_l\|-1)}.$$

The bound (5.4) is explained by first recalling that in $\sigma_l(\mathbf{x}_n)$ there are $\|\mathcal{C}_l\| - 1$ connections between different sets $\{C_i\}$. Equation (5.3) is then multiplied and divided by corresponding connection probabilities with respect to $\hat{\mu}_{x_0}$ to give the $\prod \gamma^{-1}(\cdots)$ term. Second, the prefactor $(\tilde{p}_1(x_0, y_1)/1) \le 1$ arises in connecting $x_0$ to the first state of $\sigma_l(\mathbf{x}_n)$ with respect to $\hat{\mu}_{x_0}$ and noting the constant form of $\widehat{\mathcal{P}}_1$. Third, in forming $\sigma_l(\mathbf{x}_n)$ from $\mathbf{y}_n$, with respect to $\nu_{x_0}$, $l - \|\mathcal{C}_l\| + 1$ connections between different sets are replaced by corresponding internal transition probabilities and divided by them. These $l - \|\mathcal{C}_l\| + 1$ divisors are then underestimated by the product of $\tilde{p}$'s.

STEP 3.　We now bound further the product terms in (5.4). Consider the subproduct

(5.5)
$$\prod_{k=s_r}^{s_{r+1}-1} \hat{t}(g_k(\mathbf{y}_n)+1, (i_k, i_{k+1}))$$

whose factors correspond to transitions between sets in subsequence $\langle C_{i_{s_r}}, \ldots, C_{i_{s_{r+1}}}\rangle$ for $1 \le r \le \|\mathcal{C}_l\| - 1$. From this subsequence, we derive a smaller subsequence in the following algorithm.

1. Let $\beta_1^r$ be the smallest index $s_r + 1 \le q \le s_{r+1}$ such that $C_{i_q} = C_{i_{s_{r+1}}}$.
2. If $\beta_1^r > s_r + 1$, let $\beta_2^r$ be the smallest index $s_r + 1 \le q \le \beta_1^r - 1$ such that $C_{i_q} = C_{i_{\beta_1^r-1}}$. Otherwise, stop.

3. Continue iteratively: If $\beta_m^r > s_r + 1$, let $\beta_{m+1}^r$ be the smallest index $s_r + 1 \le q \le \beta_m^r - 1$ such that $C_{i_q} = C_{i_{\beta_m^r - 1}}$. Otherwise, stop. Recalling the definition of $s_r$, there are at most $\|\mathcal{C}_l\| - r$ distinct sets in the sequence $\langle C_{i_{s_r}}, \ldots, C_{i_{s_{r+1}}} \rangle$. The above process finishes in $n(r) \le \|\mathcal{C}_l\| - r$ steps to find $\beta_{n(r)}^r = s_r + 1$.

EXAMPLE 2. With respect to the path $\mathbf{x}_n$ in Example 1, we consider the algorithm for $r = 1$. We saw that $s_1 = 3$ and $s_2 = 15$, and

$$\langle C_{i_{s_1}}, C_{i_{s_1+1}}, \ldots, C_{i_{s_2}} \rangle = \langle C_8, C_7, C_5, C_7, C_6, C_5, C_6, C_4, C_2, C_4, C_3, C_1, C_3 \rangle.$$

Here, there are $n(1) = 4$ distinct sets and $\beta_1^1 = s_1 + 10$ is the smallest index so that $C_{i_q} = C_3$. Similarly, $\beta_2^1 = s_1 + 7$ is smallest, where $C_{i_q} = C_{i_{s_1+9}} = C_4$. Also, $\beta_3^1 = s_1 + 4$ and $\beta_4^1 = s_1 + 1$.

By construction, the terms

$$\hat{t}(g_{s_r}(\mathbf{y}_n) + 1, (i_{s_r}, i_{\beta_{n(r)}^r})),$$

$$\hat{t}(g_{\beta_{n(r)}^r}(\mathbf{y}_n) + 1, (i_{\beta_{n(r)}^r}, i_{\beta_{n(r)-1}^r})), \ldots, \hat{t}(g_{\beta_2^r}(\mathbf{y}_n) + 1, (i_{\beta_2^r}, i_{\beta_1^r}))$$

all appear as factors in (5.5). Also, by monotonicity of $\hat{t}(n, (i, j))$,

$$\hat{t}(g_{s_r}(\mathbf{y}_n) + 1, (i_{s_r}, i_{\beta_{n(r)}^r})) \prod_{k=1}^{n(r)-1} \hat{t}(g_{\beta_{k+1}^r}(\mathbf{y}_n) + 1, (i_{\beta_{k+1}^r}, i_{\beta_k^r}))$$

(5.6)

$$\le \hat{t}(g_{s_r}(\mathbf{y}_n) + 1, (i_{s_r}, i_{\beta_{n(r)}^r})) \prod_{k=1}^{n(r)-1} \hat{t}(g_{s_r}(\mathbf{y}_n) + n(r) - k + 1, (i_{\beta_{k+1}^r}, i_{\beta_k^r})).$$

Also, by construction, the $r$th switch time between sets $C_{i_{s_r}}$ and $C_{i_{s_{r+1}}}$ in the rearranged path $\sigma_l(\mathbf{x}_n)$ is less than the last time to switch to $C_{i_{s_{r+1}}}$ in path $\mathbf{y}_n$:

$$g_r(\sigma_l(\mathbf{x}_n)) \le g_{s_r}(\mathbf{y}_n).$$

So, by monotonicity again, the right-hand side of (5.6) is bounded above by $\gamma(g_r(\sigma_l(\mathbf{x}_n)) + 1, (i_{s_r}, i_{s_{r+1}}))$. Also, in particular, it will be convenient to note the gross bound, because $\hat{t}(n, (i, j)) \le 1$ applies to those terms in (5.5) not covered by (5.6), that $\prod_{k=s_r}^{s_{r+1}-1} \hat{t}(g_k(\mathbf{y}_n) + 1, (i_k, i_{k+1})) \le \gamma(g_r(\sigma_l(\mathbf{x}_n)) + 1, (i_{s_r}, i_{s_{r+1}}))$ and so

$$(5.7) \quad \prod_{k=1}^{l} \hat{t}(g_k(\mathbf{y}_n) + 1, (i_k, i_{k+1})) \le \prod_{k=1}^{\|\mathcal{C}_l\|-1} \gamma(g_k(\sigma_l(\mathbf{x}_n)) + 1, (i_{s_k}, i_{s_{k+1}})).$$

STEP 4.  We now consider cases when $l$ is large and small. Suppose first that $l$ is small, namely $l \leq \|\mathcal{C}_l\|(\|\mathcal{C}_l\| - 1)/2 + N + M - 1$. Then we have the bound, noting (5.3), (5.4) and (5.7), that

$$(5.8) \qquad \nu_{x_0}(\mathbf{X}_n = \mathbf{y}_n) \leq (1/\tilde{p})^l \hat{\mu}_{x_0}(\mathbf{X}_n = \sigma_l(\mathbf{x}_n)).$$

Suppose now that $l$ is large, that is, $l > \|\mathcal{C}_l\|(\|\mathcal{C}_l\| - 1)/2 + N + M - 1$. Whereas the chain can only make at most $N + M - 1$ consecutive downward switches (i.e., from sets $C_i$ to $C_j$ for $i < j$), in $q > N + M - 1$ switches there will be at least $[q/(N + M - 1)]$ upward switches from sets $C_i$ to $C_j$ for $i > j$.

Whereas $n(k) \leq \|\mathcal{C}_l\| - k$ and so $\sum_{k=1}^{\|\mathcal{C}_l\|-1} n(k) \leq \|\mathcal{C}_l\|(\|\mathcal{C}_l\| - 1)/2$, we see carefully in Step 3 that we take at most $\|\mathcal{C}_l\|(\|\mathcal{C}_l\| - 1)/2$ factors from $\prod_{k=1}^{l} \hat{t}(g_k(\mathbf{y}_n) + 1, (i_k, i_{k+1}))$ whose product is then dominated by $\prod_{k=1}^{\|C_l\|-1} \gamma_k(\sigma_l(\mathbf{x}_n) + 1, (i_{s_k}, i_{s_{k+1}}))$. Hence, remaining in the original product are at least $l - \|\mathcal{C}_l\|(\|\mathcal{C}_l\| - 1)/2$ uncommitted factors of which at least

$$\underline{l} = \lfloor (l - \|\mathcal{C}_l\|(\|\mathcal{C}_l\| - 1)/2)/(N + M - 1) \rfloor$$

correspond to upward transitions.

Then, using monotonicity of $\underline{t}(n)$, we have

$$\prod_{k=1}^{l} \hat{t}(g_k(\mathbf{y}_n) + 1, (i_k, i_{k+1})) \leq \prod_{k=1}^{\|C_l\|-1} \gamma_k(\sigma_l(\mathbf{x}_n) + 1, (i_{s_k}, i_{s_{k+1}})) \prod_{j=1}^{\underline{l}} \underline{t}(j).$$

Furthermore, noting (5.3) and (5.4), we have, for $l$ large,

$$(5.9) \qquad \nu_{x_0}(\mathbf{X}_n = \mathbf{y}_n) \leq (1/\tilde{p})^l \left[\prod_{j=1}^{\underline{l}} \underline{t}(j)\right] \hat{\mu}_{x_0}(\mathbf{X}_n = \sigma_l(\mathbf{x}_n)).$$

STEP 5.  We now estimate the size of the set $\sigma_l^{-1}(\sigma_l(\mathbf{x}_n))$. Observe that the ordering of states within the $l+1$ subpaths in $\sigma_l(\mathbf{x}_n)$ is preserved among the paths $\sigma_l^{-1}(\sigma_l(\mathbf{x}_n))$ with $l$ switches. Then, to overestimate $|\sigma_l^{-1}(\sigma_l(\mathbf{x}_n))|$, we need only to specify the sequence in which the pairwise distinct sets $C_{j_1} \neq C_{j_2} \neq \cdots \neq C_{j_{l+1}}$ are visited and how long each visit takes, since once the ordering of the sets and switch times are fixed, the arrangement within the $l+1$ subpaths is determined.

A simple overcount of this procedure yields that

$$|\sigma_l^{-1}(\sigma_l(\mathbf{x}_n))| \leq \binom{n}{l} M^{l+1}.$$

Therefore, from (5.8) and (5.9) we have that

$$\nu_{x_0}(\mathbf{X}_n \in \sigma_l^{-1}(\sigma_l(\mathbf{x}_n)))$$

(5.10)
$$\leq \begin{cases} \binom{n}{l} M^{l+1} \tilde{p}^{-l} \hat{\mu}_{x_0}(\mathbf{X}_n = \sigma_l(\mathbf{x}_n)), & \text{for } l \text{ small,} \\ \binom{n}{l} M^{l+1} \tilde{p}^{-l} \left[ \prod_{i=1}^{l} \underline{t}(i) \right] \hat{\mu}_{x_0}(\mathbf{X}_n = \sigma_l(\mathbf{x}_n)), & \text{for } l \text{ large.} \end{cases}$$

STEP 6. By Stirling's formula,

$$\frac{1}{n} \log \binom{n}{l} = o(1) - \frac{l}{n} \log\left(\frac{l}{n}\right) - \frac{n-l}{n} \log\left(\frac{n-l}{n}\right).$$

With this estimate, we now analyze the factor $\binom{n}{l} M^l \prod_{i=1}^{l} \underline{t}(i)$ in (5.10). We consider cases when $l = o(n)$ and when $l \leq n$ is otherwise.

*Case* 1. When $l = l_n = o(n)$, then $\log \binom{n}{l_n}/n \to 0$. Also, $M^{l_n} = e^{o(n)}$, $\tilde{p}^{-l_n} = e^{o(n)}$ and $\prod_{i=1}^{l_n} \underline{t}(i) = e^{o(n)}$.

*Case* 2. When $l = l_n$ satisfies $\limsup l_n/n \geq \varepsilon$ for some $0 < \varepsilon \leq 1$, let $n'$ be a maximal subsequence. Then $(\log \binom{n'}{l_{n'}})/n' = O(1)$, $(\log M^{l_{n'}})/n' \leq 1 + \log M$ and $(\log \tilde{p}^{-l_{n'}})/n' \leq 1 + \log \tilde{p}^{-1}$, but, as $\underline{t}(n') \downarrow 0$ and $\limsup \underline{l}'_n/n' \geq \varepsilon/(N+M-1)$, we have $\log[\prod_{i=1}^{l_{n'}} \underline{t}(i)]/n' \to -\infty$ as $n' \to \infty$.

Therefore, with respect to a $C_n = e^{o(n)}$, independent of $l \geq 1$ and the path, we have from (5.10) that

$$\nu_{x_0}(\mathbf{X}_n \in \sigma_l^{-1}(\sigma_l(\mathbf{x}_n))) \leq C_n \hat{\mu}_{x_0}(\mathbf{X}_n = \sigma_l(\mathbf{x}_n)).$$

STEP 7. Let $l \geq 1$, and let $\overline{A}_n(l) = \bigcup_{j=1}^{\min\{N+M-1,l\}} A'_n(j)$. Let also $\widehat{A}_n(l) = \sigma_l(\{Z_n \in \Gamma, A_n(l)\})$ and $\widetilde{A}_n(l) = \{Z_n \in \Gamma, \overline{A}_n(l)\}$. Whereas the average $Z_n$ is independent of the order of observations $\{X_1, \ldots, X_n\}$,

$$\widehat{A}_n(l) \subset \widetilde{A}_n(l) \quad \text{and} \quad \{Z_n \in \Gamma, A_n(l)\} = \sigma_l^{-1} \sigma_l(Z_n \in \Gamma, A_n(l)).$$

Then we can write

$$\nu_{x_0}(Z_n \in \Gamma, A_n(l)) = \nu_{x_0}(\sigma_l^{-1}(\sigma_l(Z_n \in \Gamma, A_n(l))))$$

$$= \nu_{x_0}\left(\mathbf{X}_n \in \bigcup_{\mathbf{x}_n \in \widehat{A}_n} \sigma_l^{-1}(\mathbf{x}_n)\right)$$

$$\leq \sum_{\mathbf{x}_n \in \widehat{A}_n} \nu_{x_0}(\mathbf{X}_n \in \sigma_l^{-1}(\mathbf{x}_n))$$

$$\leq e^{o(n)} \sum_{\mathbf{x}_n \in \widehat{A}_n} \hat{\mu}_{x_0}(\mathbf{X}_n = \mathbf{x}_n)$$

$$\leq e^{o(n)} \sum_{\mathbf{x}_n \in \widetilde{A}_n} \hat{\mu}_{x_0}(\mathbf{X}_n = \mathbf{x}_n)$$

$$= e^{o(n)} \hat{\mu}_{x_0}(Z_n \in \Gamma, \overline{A}_n(l)).$$

STEP 8. Whereas $\bigcup_{l \geq 1} \overline{A}_n(l) \cup A_n(0) \subset \{\mathbf{X}_n \text{ enters each } C_i \text{ at most once}\}$, we have

$$\sum_{i=0}^{n-1} \nu_{x_0}(Z_n \in \Gamma, A_n(i))$$

(5.11)

$$\leq ((1 + (n-1)e^{o(n)}) \hat{\mu}_{x_0}(Z_n \in \Gamma, \mathbf{X}_n \text{ enters each } C_i \text{ at most once}).$$

Then, noting (5.1), (5.2) and (5.11), we have

$$\limsup \frac{1}{n} \log \mathbb{P}_\pi(Z_n \in \Gamma)$$

$$\leq \max_{\substack{x_0 \in \Sigma \\ \pi(x_0) > 0}} \limsup \frac{1}{n} \log \hat{\mu}_{x_0}(Z_n \in \Gamma, \mathbf{X}_n \text{ enters each } C_i \text{ at most once}).$$

Applying Lemma 5.1 completes the proof. □

**6. Path surgery lower bound.** The lower bound strategy is informed by the upper bound result. Namely, given the rearranged paths focused on in the upper bound surgery, we can more or less restrict to them and gain lower bounds.

PROOF OF PROPOSITION 4.7. When $N + M = 1$, $P$ is irreducible, $C_1 = \Sigma$ and $\mathcal{D} = \varnothing$. Then $\check{\mathcal{P}}_n = P_n$ for all $n \geq 1$ and so $\mathbb{P}_\pi(Z_n \in \Gamma) = \check{\mu}_\pi(Z_n \in \Gamma)$. Also, as in the upper bound, $\mathbf{X}_n$ does not switch in this case. Hence, the lower bound holds trivially.

We now assume that $N + M \geq 2$. Consider the subset $B \subset \Sigma^n$ formed via the following procedure.

1. Tor $1 \leq m \leq N + M$, let $J_1, J_2, \ldots, J_m$ be subpaths that belong, respectively, to distinct sets $C_{i_1}, C_{i_2}, \ldots, C_{i_m}$, where $\{i_1, \ldots, i_m\} \subset \mathcal{G}$. Let $j_i = |J_i|$, $J_i = \langle y_1^i, \ldots, y_{j_i}^i \rangle$ and $J_{i,2} = \langle y_2^i, \ldots, y_{j_i}^i \rangle$ when $|j_i| \geq 2$ for $1 \leq i \leq m$. We impose now that the lengths satisfy $\sum_{i=1}^m j_i = n - E(N, M)$.
2. When $m \geq 2$, we connect subpaths $J_s$ and $J_{s+1}$ for $s = 1, \ldots, m-1$ as follows. Let $0 \leq k \leq N + M - 2$ be the number of sets entered in the connection and let $L_k$ with $i = i_s$ and $j = i_{s+1}$, $Q_k$ and $V_k = \{\mathbf{x}^{s,0}, \ldots, \mathbf{x}^{s,k+1}\}$ be as near (2.9). Denote $\mathbf{w}^s = \langle \mathbf{x}^{s,0}, \ldots, \mathbf{x}^{s,k+1} \rangle$ and $k_s = |\mathbf{w}^s|$. Also, denote $b(s) = j_s + \sum_{i=1}^{s-1}(j_i + k_i)$. Let now $\mathbf{w}^s$ be such that

$$\mathbb{P}_{(b(s),y_2^s)}(\mathbf{X}_{b(s)+1}^{b(s)+k_s+1} = \langle \mathbf{w}^s, y_1^{s+1} \rangle) = \underline{\gamma}^1(b(s)+1, y_{j_s}^s, y_1^{s+1}).$$

Then, in particular, as $\sum_{i=1}^{s-1} k_i \leq E(N, M)$, we have

$$\mathbb{P}_{(b(s),y_2^s)}(\mathbf{X}_{b(s)+1}^{b(s)+k_s+1} = \langle \mathbf{w}^s, y_1^{s+1} \rangle) \geq \underline{\gamma}^1(b(s)+1, (i_s, i_{s+1}))$$

$$\geq \check{\gamma}^0\left(\sum_{i=1}^{s} j_i + 1, (i_s, i_{s+1})\right).$$

3. For $m \geq 2$, as $\sum_{i=1}^{m-1} k_i \leq E(N, M)$, the length of the concatenation satisfies

$$L = |\langle J_1, \mathbf{w}^1, J_2, \ldots, \mathbf{w}^{m-1}, J_m \rangle|$$

$$= n - E(N, M) + \sum_{i=1}^{m-1} k_i \leq n.$$

When $m = 1$, the length $L = |\langle J_1 \rangle| = n - E(N, M)$.

If now $L < n$, we then augment the last subpath $J_m$ by $n - L \leq E(N, M)$ states in $C_{i_m}$. Specifically, define

$$J'_m = \begin{cases} J_m, & \text{if } L = n, \\ \langle J_m, x_1^m, \ldots, x_{n-L}^m \rangle, & \text{if } L < n, \end{cases}$$

where $\langle y_{j_m}^m, x_1^m, \ldots, x_{n-L}^m \rangle$ is a sequence of $n - L + 1$ elements in $C_{i_m}$ with positive weight. Let also $J'_{m,2} = J_{m,2}$ when $L = n$ and $J'_{m,2} = \langle J_{m,2}, x_1^m, \ldots, x_{n-L}^m \rangle$ otherwise.

Now let

$$\mathbf{x}_n = \begin{cases} \langle J_1, \mathbf{w}^1, \ldots, \mathbf{w}^{m-1}, J'_m \rangle, & \text{when } m \geq 2, \\ \langle J'_1 \rangle, & \text{when } m = 1. \end{cases}$$

Finally, we define $B$ as the set of all such sequences $\mathbf{x}_n$ possible.

Now write

$$\mathbb{P}_\pi(Z_n \in \Gamma_2)$$
$$\geq \mathbb{P}_\pi(Z_n \in \Gamma_2, \mathbf{X}_n \in B)$$
$$= \sum_{\mathbf{x}_n \in \{Z_n \in \Gamma_2\} \cap B} \mathbb{P}_\pi(\mathbf{X}_{j_1} = J_1) \underline{\gamma}(j_1 + 1, y_{j_1}^1, y_1^2)$$
$$\times \mathbb{P}_{(j_1+k_1+1, y_1^2)}(\mathbf{X}_{j_1+k_1+2}^{j_1+k_1+j_2} = J_{2,2})$$

(6.1)

$$\times \cdots \times \underline{\gamma}(b(m-1)+1, y_{j_{m-1}}^{m-1}, y_1^m)$$

$$\times \mathbb{P}_{(\sum_1^{m-1}(j_i+k_i)+1, y_1^m)}(\mathbf{X}_{\sum_1^{m-1}(j_i+k_i)+2}^{\sum_1^{m-1}(j_i+k_i)+n-L} = J'_{m,2})$$

$$\geq c(L)\check{\mu}_\pi(Z_{n-E(N,M)} \in \Gamma_{2,n}, \mathbf{X}_{n-E(N,M)} \text{ only enters } \{C_i : i \in \mathcal{G}\}$$

with at most one visit to each set),

where

$$c(L) = \begin{cases} \mathbb{P}_{(b(m), y_{j_m}^m)}(\mathbf{X}_{b(m)+1}^{b(m)+n-L} = \langle x_1^m, \ldots, x_{n-L}^m \rangle), & \text{when } n > L, \\ 1, & \text{when } n = L, \end{cases}$$

and

$$\Gamma_{2,n} = \frac{n}{n - E(N,M)} B\left(x_0, \frac{a}{2} - \frac{E(N,M)\|f\|}{n}\right).$$

In the last step, we rewrote $\mathbb{P}_\pi$ in terms of the measure $\check{\mu}_\pi$ by collapsing together the subpaths $\{J_i\}$. At the same time, since the collapsed path $\langle J_1, \ldots, J_m \rangle$ is of length $n - E(N,M)$, we correct the set $\Gamma_2$ to $\Gamma_{2,n}$.

We now estimate the prefactor $c(L)$. With respect to the minimum probability $p_{\min}$ [cf. (2.5)] and $n > L$ large, as $P_n \to P$, we can certainly bound

$$\mathbb{P}_{(b(m), y_{j_m}^m)}(\mathbf{X}_{b(m)+1}^{b(m)+n-L} = \langle x_1^m, \ldots, x_{n-L}^m \rangle) \geq p_{\min}^{E(N,M)}/2.$$

Therefore, $\lim(\log c(L))/n = 0$.

Hence, the proposition follows by taking lim inf in (6.1) and simple estimates. $\square$

**7. Homogeneous "rest cost" replacement.** We replace certain a priori nonhomogeneous "resting" weights with homogeneous ones for both upper and lower bound estimates.

PROOFS OF PROPOSITIONS 4.3 AND 4.8. The proofs follow as direct corollaries of the more general Proposition 7.1 below. $\square$

PROPOSITION 7.1. *Let $\{B_n\} \subset \mathbb{R}^d$ be a sequence of Borel sets.*
Upper bound. *For $\varepsilon_1, \varepsilon_2 > 0$, we have*

$$\limsup \frac{1}{n} \log \hat{\mu}_{\pi, \varepsilon_1, \varepsilon_2}(\mathbf{X}_n \in B_n) \leq \limsup \frac{1}{n} \log \bar{\mu}_{\pi, \varepsilon_1, \varepsilon_2}(\mathbf{X}_n \in B_n).$$

Lower bound. *Suppose $n_0(\{P_n\}) = 1$. Then we have*

$$\liminf \frac{1}{n} \log \check{\mu}_\pi(\mathbf{X}_n \in B_n) \geq \liminf \frac{1}{n} \log \underline{\mu}_\pi(\mathbf{X}_n \in B_n).$$

PROOF. We prove the lower bound part, because the upper bound estimate follows analogously and more simply. Let $G = \{(s,t) : p(s,t) > 0$ where $s, t \in C_i$ for $i \in \mathcal{G}\}$. As $P_n \to P$, the state space is finite and, by assumption $n_0 = 1$, there exists $\alpha > 0$ and a sequence $\alpha \leq m(k) \uparrow 1$ such that $m(k) \leq p_k(s,t)/p(s,t)$ for all $(s,t) \in G$ and $k \geq 1$. Write now that

$\check{\mu}_\pi(\mathbf{X}_n \in B_n)$

$$= \sum_{x_0 \in \Sigma} \sum_{\mathbf{x}_n \in B_n} \pi(x_0) \prod_{i=1}^n \check{p}_i(x_{i-1}, x_i)$$

$$= \sum_{x_0 \in \Sigma} \sum_{\mathbf{x}_n \in B_n} \pi(x_0)$$

$$\times \prod_{(x_{i-1},x_i) \in G^c} \check{p}_i(x_{i-1}, x_i) \prod_{(x_{i-1},x_i) \in G} \frac{\check{p}_i(x_{i-1}, x_i)}{p(x_{i-1}, x_i)} p(x_{i-1}, x_i)$$

$$\geq \sum_{x_0 \in \Sigma} \sum_{\mathbf{x}_n \in B_n} \pi(x_0)$$

$$\times \prod_{(x_{i-1},x_i) \in G^c} \underline{p}_i(x_{i-1}, x_i) \prod_{(x_{i-1},x_i) \in G} \frac{p_i(x_{i-1}, x_i)}{p(x_{i-1}, x_i)} \underline{p}_i(x_{i-1}, x_i)$$

$$\geq \left[\prod_{i=1}^n m(i)\right] \sum_{x_0 \in \Sigma} \sum_{\mathbf{x}_n \in B_n} \pi(x_0)$$

$$\times \prod_{(x_{i-1},x_i) \in G^c} \underline{p}_i(x_{i-1}, x_i) \prod_{(x_{i-1},x_i) \in G} \underline{p}_i(x_{i-1}, x_i)$$

$$= \left[\prod_{i=1}^n m(i)\right] \underline{\mu}_\pi(\mathbf{X}_n \in B_n).$$

Indeed, for the first bound, we note, if $(x_{i-1}, x_i) \notin G$, that $\check{p}_i(x_{i-1}, x_i) = \underline{p}_i(x_{i-1}, x_i)$ when $(x_{i-1}, x_i)$ connects distinct sets in $\mathcal{G}$, and $\check{p}_i(x_{i-1}, x_i) \geq 0 = \underline{p}_i(x_{i-1}, x_i)$ otherwise. The second bound follows by monotonicity of $\{m(i)\}$.

Then the proposition lower bound follows as $(\sum_1^n \log m(i))/n \to 0$. $\square$

**8. Upper coarse graining bounds.** The plan is to optimize over a coarse graining of the possible locations $Z_n$ visits in $\mathbb{K}$ and associated visit times. Some additional definitions which build on those in Section 5 are required in this effort.

Define, for $1 \leq H \leq N + M$ and $\underline{i}_H = \langle i_1, \ldots, i_H \rangle$ composed of distinct indices in $\{1, \ldots, N + M\}$, that

$C(\underline{i}_H) = \{\mathbf{X}_n$ starts in $C_{i_1}$ and enters successively $C_{i_2}, \ldots, C_{i_H}\}$.

Also, let $k_0 = 0$, $k_H = n$ and, when $H \geq 2$, let $1 \leq k_1 < \cdots < k_{H-1} \leq n-1$, and denote $\mathbf{k}_H = \langle k_0, \ldots, k_H \rangle$ and

$$S_n(\mathbf{k}_H) = \{\mathbf{X}_n \text{ switches at times } k_1, k_2, \ldots, k_{H-1}\}.$$

Let also

$$\mathbf{v}_{\mathbf{k}_H} = \langle k_1/n, (k_2 - k_1)/n, \ldots, (n - k_{H-1})/n \rangle.$$

We now specify a certain cube decomposition. For $\mathbf{v} \in \Omega_H$ and $z \in \mathbb{R}^d$, recall the set $D(H, \mathbf{v}, z)$ [cf. (2.6)] and let $D(H, \mathbf{v}, B) = \bigcup_{z \in B} D(H, \mathbf{v}, z)$ for sets $B \subset \mathbb{R}^d$.

Let now $F_1$ be the regular partition of $\mathbb{K}$ into $2^d$ closed cubes, $\{\Delta_s^1 : 1 \leq s \leq 2^d\}$, whose interiors nonintersect and $\bigcup_s \Delta_s^1 = \mathbb{K}$. For $n \geq 2$, let also $F_n$ be the regular refinement of $F_{n-1}$ into $2^{n-1}(2^d)$ closed cubes, $\{\Delta_s^n : 1 \leq s \leq 2^{n-1}(2^d)\}$, where also $\bigcup_s \Delta_s^n = \mathbb{K}$. Observe also that the $(2^{n-1}(2^d))^H$ subcubes formed from $F_n$, $\{\Delta(n, \mathbf{s}) = \Delta_{s_1}^n \times \cdots \times \Delta_{s_H}^n : 1 \leq s_i \leq 2^{n-1}(2^d)\}$, refine $\mathbb{K}^H$ as well.

For $B \subset \mathbb{K}$ and $j \geq 1$, define

$$D_j(H, \mathbf{v}, B) = \bigcup \{\Delta(j, \mathbf{s}) : \Delta(j, \mathbf{s}) \cap D(H, \mathbf{v}, B) \neq \varnothing\}$$

be the nonempty union of all subcubes with respect to $j$th partition which intersect $D(H, \mathbf{v}, B)$. Let also

$$F(H, n, \mathbf{v}, B) = \{\mathbf{s} : \Delta(n, \mathbf{s}) \subset D_n(H, \mathbf{v}, B)\}.$$

For $\alpha > 0$, let $m_\alpha$ be the first partition level $m$ so that, for each $1 \leq l \leq N + M$, $|\mathbb{I}_{l,\varepsilon_1,\varepsilon_2}(x) - \mathbb{I}_{l,\varepsilon_1,\varepsilon_2}(y)| \leq \alpha$ when $|x - y| \leq \text{diam}(\Delta(m, \cdot))$ and $x, y \in Q_{l,\varepsilon_1,\varepsilon_2}$.

We also need the following technical lemmas, which can be skipped on first reading.

LEMMA 8.1. *For distinct $i, j \in \mathcal{G}$, we have*

$$\mathcal{U}_0(i, j) = \limsup \frac{1}{n} \log \gamma(n, (i, j)).$$

PROOF. Write the left-hand side as

$$\limsup \frac{1}{n} \log \gamma(n, (i, j))$$

$$= \limsup \max_{0 \leq k \leq M-2} \max_{L_k} \sum_{s=0}^{k} \frac{1}{n} \log \hat{t}(n + s, (l_s, l_{s+1}))$$

$$= \max_{0 \leq k \leq M-2} \max_{L_k} \sum_{s=0}^{k} \lim \frac{1}{n} \log \hat{t}(n, (l_s, l_{s+1}))$$

$$= \max_{0 \leq k \leq M-2} \max_{L_k} \sum_{s=0}^{k} v(l_s, l_{s+1}) = \mathcal{U}_0(i, j),$$

where the second and third lines follow since the limit $\lim \log \hat{t}(n,(k,l))/n = v(k,l)$ holds from Lemma 4.1. □

In the next result, let $1 \leq H \leq N + M$ and let $\Gamma \subset \mathbb{K}$ be a closed set. Let also $\mathbb{I}_l^\theta = \min\{\mathbb{I}_{l,\varepsilon_1,\varepsilon_2}, \theta\}$ for $\theta \geq 1$ and $1 \leq l \leq N + M$.

LEMMA 8.2. *Let $\mathbf{v}^n \in \Omega_H$ be a convergent sequence, $\lim_n \mathbf{v}^n = \mathbf{v} \in \Omega_H$. Then, for any $\underline{i}_H$, we have*

$$\limsup_{\theta \uparrow \infty} \limsup_{m \uparrow \infty} \limsup_{n \to \infty} \inf_{\mathbf{x} \in D_m(H, \mathbf{v}^n, \Gamma)} \sum_{j=1}^H v_j^n \mathbb{I}_{i_j}^\theta(x_j)$$

$$\geq \inf_{\mathbf{x} \in D(H, \mathbf{v}, \Gamma) \cap \mathbb{K}} \sum_{j=1}^H v_j \mathbb{I}_{i_j, \varepsilon_1, \varepsilon_2}(x_j).$$

PROOF. Whereas $D_m(H, \mathbf{v}^n, \Gamma) \subset \mathbb{K}^H$ and $\mathbb{K}^H$ is compact, we can find a convergent sequence $\mathbf{x}^{m,n_k} \in D_m(H, \mathbf{v}^{n_k}, \Gamma) \to \mathbf{x}^m \in \mathbb{K}^H$ so that by lower semicontinuity of $\{\mathbb{I}_l^\theta\}$,

$$\limsup_{n \to \infty} \inf_{\mathbf{x} \in D_m(H, \mathbf{v}^n, \Gamma)} \sum_{j=1}^H v_j^n \mathbb{I}_{i_j}^\theta(x_j) = \lim_{k \to \infty} \sum_{j=1}^H v_j^{n_k} \mathbb{I}_{i_j}^\theta(x_j^{m,n_k})$$

$$\geq \sum_{j=1}^H v_j \mathbb{I}_{i_j}^\theta(x_j^m).$$

Now, out of $\{\mathbf{x}^m\} \subset \mathbb{K}^H$, let $\mathbf{x}^{m_j} \to \mathbf{x} \in \mathbb{K}^H$ be a convergent subsequence on which $\limsup_{m \uparrow \infty} \sum_{j=1}^H v_j \mathbb{I}_{i_j}^\theta(x_j^m)$ is attained. Also observe that $\mathbb{I}_l^\theta(x_l) \to \mathbb{I}_{l,\varepsilon_1,\varepsilon_2}(x_l)$ for $1 \leq l \leq N + M$ as $\theta \uparrow \infty$. Then, again by lower semicontinuity,

$$\limsup_{\theta \uparrow \infty} \limsup_{m \uparrow \infty} \limsup_{n \to \infty} \inf_{\mathbf{x} \in D_m(H, \mathbf{v}^n, \Gamma)} \sum_{j=1}^H v_j^n \mathbb{I}_{i_j}^\theta(x_j)$$

$$\geq \limsup_{\theta \uparrow \infty} \sum_{j=1}^H v_j \mathbb{I}_{i_j}^\theta(x_j)$$

$$= \sum_{j=1}^H v_j \mathbb{I}_{i_j, \varepsilon_1, \varepsilon_2}(x_j).$$

To finish the argument, we show that $\mathbf{x} \in D(H, \mathbf{v}, \Gamma) \cap \mathbb{K}^H$. By construction, the diameters of the partitioning cubes $\Delta(m, \cdot)$ uniformly vanish as $m \uparrow \infty$. As $D_m(H, \mathbf{v}^{n_k}, \Gamma)$ is composed of cubes which intersect $D(H, \mathbf{v}^{n_k}, \Gamma)$, we have that any point in $D_m(H, \mathbf{v}^{n_k}, \Gamma)$ is at most a distance $\text{diam}(\Delta(m, \cdot))$ away from $D(H,\mathbf{v}^{n_k},\Gamma)\cap \mathbb{K}^H$. Hence, there are points $\mathbf{y}^{m,n_k}\in D(H,\mathbf{v}^{n_k},\Gamma)\cap \mathbb{K}^H$ such that $|\mathbf{x}^{m,n_k}-\mathbf{y}^{m,n_k}|\leq \mathrm{diam}(\Delta(m,\cdot))$. Let $\mathbf{y}^{m,n'_k}\to \mathbf{y}^m\in \mathbb{K}^H$ be a convergent subsequence. We have then $|\mathbf{x}^m-\mathbf{y}^m|\leq \mathrm{diam}(\Delta(m,\cdot))$. Now since $\Gamma$ is closed and $\sum_{j=1}^H v_j^{n'_k} y_j^{m,n'_k}\in \Gamma$ for all $m,k$, we have

$$\lim_m \lim_k \sum_{j=1}^H v_j^{n'_k} y_j^{m,n'_k} = \lim_m \sum_{j=1}^H v_j y_j^m = \sum_{j=1}^H v_j x_j \in \Gamma$$

and so $\mathbf{x}\in D(H,\mathbf{v},\Gamma)\cap \mathbb{K}^H$.  □

PROOF OF PROPOSITION 4.4.  When $N+M=1$, there is only one irreducible subset $C_1=\Sigma$ and $\overline{\mathcal{P}}_{k,\varepsilon_1,\varepsilon_2}=P(1,\varepsilon_1)$ for $k\geq 2$. So, modulo a first transition (with respect to the constant matrix $\overline{\mathcal{P}}_1$), the measure $\bar{\mu}_\pi$ is a "homogeneous nonnegative process" with respect to $P(1,\varepsilon_1)$. Also, whereas there can no "repeat visits" and $\mathbb{J}_{\mathcal{U}_0,\varepsilon_1,\varepsilon_2}=\mathbb{I}_{1,\varepsilon_1,\varepsilon_2}$ in this case, the proposition follows from the LDP in Proposition 2.1.

We now assume that $N+M\geq 2$. Also, to reduce notation we suppress subscripts $\varepsilon_1$ and $\varepsilon_2$ when there is no confusion in the following text.

STEP 1.  Whereas $Z_n$ takes only values in the set $\mathbb{K}$, we have

(8.1)
$$\bar{\mu}_{\pi,\varepsilon_1,\varepsilon_2}(Z_n\in\overline{\Gamma},\mathbf{X}_n \text{ enters each } C_i \text{ at most once})$$
$$= \sum_{1\leq H\leq N+M} \sum_{\underline{i}_H} \bar{\mu}_{\pi,\varepsilon_1,\varepsilon_2}(Z_n\in\overline{\Gamma}\cap\mathbb{K}, A'_n(H-1), C(\underline{i}_H)),$$

where the sum on $\underline{i}_H$ is over $\binom{N+M}{H}H!$ possibilities.

STEP 2.  We first consider the case when "switching" actually occurs. Let $2\leq H\leq N+M$ and fix indices $\underline{i}_H$. Write, for $n>N+M$ (larger than the number of switches), that

(8.2)
$$\bar{\mu}_\pi(Z_n\in\overline{\Gamma}\cap\mathbb{K}, A'_n(H-1), C(\underline{i}_H))$$
$$= \sum_{\mathbf{k}_H} \bar{\mu}_\pi(Z_n\in\overline{\Gamma}\cap\mathbb{K}, A'_n(H-1), C(\underline{i}_H), S_n(\mathbf{k}_H)),$$

where the sum on $\mathbf{k}_H$ comprises $\binom{n-1}{H-1}$ possibilities.

For convenience, denote $\overline{B}=\overline{\Gamma}\cap\mathbb{K}$ and

$$E_n = A'_n(H-1)\cap C(\underline{i}_H)\cap S_n(\mathbf{k}_H).$$

Let also $\alpha>0$ and let $m\geq m_\alpha$. Recall from part Section 2.4 that $Z_i^j\in\mathbb{K}$ for $i\leq j$, and so we may write the summand in (8.2) equal to

$$\bar{\mu}_\pi(\langle Z_1^{k_1},\ldots,Z_{k_{H-1}+1}^n\rangle\in D(H,\mathbf{v}_{\mathbf{k}_H},\overline{B})\cap\mathbb{K}^H, E_n)$$

$$\leq \bar{\mu}_\pi(\langle Z_1^{k_1},\ldots,Z_{k_{H-1}+1}^n\rangle \in D_m(H,\mathbf{v}_{\mathbf{k}_H},\overline{B}),E_n)$$

(8.3)
$$= \bar{\mu}_\pi\left(\langle Z_1^{k_1},\ldots,Z_{k_{H-1}+1}^n\rangle \in \bigcup_{\mathbf{s}}\Delta(m,\mathbf{s}),E_n\right)$$

$$\leq \sum_{\mathbf{s}} \bar{\mu}_\pi(Z_1^{k_1}\in\Delta_{s_1}^m,\ldots,Z_{k_{H-1}+1}^n\in\Delta_{s_H}^m,E_n),$$

where the union and sum is over $\mathbf{s}\in F(H,m,\mathbf{v}_{\mathbf{k}_H},\overline{B})$

STEP 3. For $1\leq l\leq N+M$, let $\pi_l$ be the uniform distribution on $C_l$ and let $\mathbb{P}_{\pi,\varepsilon_1,\varepsilon_2}^l$ denote the homogeneous nonnegative measure on $C_l$ formed from CON with $U_n\equiv P(l,\varepsilon_1,\varepsilon_2)$ and initial distribution $\pi$. Let also $\theta > \max_{1\leq l\leq N+M}\max_{x\in Q_{l,\varepsilon_1,\varepsilon_2}}\mathbb{I}_{l,\varepsilon_1,\varepsilon_2}(x)$ be a number larger than the maxima of the rate functions on their domains of finiteness (cf. Proposition 2.1).

We now use the Markov property (2.2) and simple estimates to further bound the summand in (8.3) as

$$\bar{\mu}_\pi(Z_1^{k_1}\in\Delta_{s_1}^m,\ldots,Z_{k_{H-1}+1}^n\in\Delta_{s_H}^m,E_n)$$

$$\leq \bar{\mu}_\pi(Z_1^{k_1}\in\Delta_{s_1}^m,\mathbf{X}_1^{k_1}\text{ in }C_{i_1})$$

$$\times \prod_{j=1}^{H-1}|C_{i_j}|\gamma(k_j+1,(i_j,i_{j+1}))$$

(8.4)
$$\times \bar{\mu}_{(\pi_{i_{j+1}},k_j+1)}(Z_{k_j+1}^{k_{j+1}}\in\Delta_{s_{j+1}}^m,\mathbf{X}_{k_j+1}^{k_{j+1}}\text{ in }C_{i_{j+1}})$$

$$\leq \prod_{j=1}^{H-1}\gamma(k_j+1,(i_j,i_{j+1}))\prod_{j=0}^{H-1}|C_{i_{j+1}}|\mathbb{P}_{(\pi_{i_{j+1}},k_j+1)}^{i_{j+1}}(Z_{k_j+1}^{k_{j+1}}\in\Delta_{s_{j+1}}^m).$$

STEP 4. Recall the definition of $\mathbb{I}_l^\theta$ just before Lemma 8.2. Let

$$c(k_{j+1}-k_j;\Delta_{s_{j+1}}^m,\theta,C_{i_{j+1}})$$
$$= \mathbb{P}_{(\pi_{i_{j+1}},k_j+1)}^{i_{j+1}}(Z_{k_j+1}^{k_{j+1}}\in\Delta_{s_{j+1}}^m)\exp((k_{j+1}-k_j)\mathbb{I}_{i_{j+1}}^\theta(\Delta_{s_{j+1}}^m)).$$

From homogeneous nonnegative large deviation upper bounds (cf. Proposition 2.1), uniformly over $H$, $\mathbf{k}_H$, the finite number of cubes $\mathbf{s}$ at level $m$, and $\underline{i}_H$, we have $c(k_{j+1}-k_j;\Delta_{s_{j+1}}^m,\theta,C_{i_{j+1}})\leq e^{o(n)}$.

Also by monotonicity $\gamma(i+1,\ldots)\leq\gamma(i,\ldots)$. Then we have (8.4) is less than

(8.5) $$e^{o(n)}\left[\prod_{j=1}^{H-1}\gamma(k_j,(i_j,i_{j+1}))\right]\exp\left\{-\sum_{j=1}^{H}(k_j-k_{j-1})\mathbb{I}_{i_j}^\theta(\Delta_{s_j}^m)\right\}.$$

STEP 5. At this point, we now bound the terms that correspond to no "switching" in (8.1), that is, when $H = 1$. For $1 \leq i_1 \leq N + M$, we have

$$(8.6) \quad \bar{\mu}_\pi(Z_n \in \overline{\Gamma} \cap \mathbb{K}, A'_n(0), C(i_1)) \leq e^{o(n)} \sum_{s_1 \in F(1,n,1,\overline{B})} \exp\{-n \mathbb{I}^\theta_{i_1}(\Delta^m_{s_1})\}.$$

STEP 6. It is convenient now to define $\gamma(0, (l, l')) = 1$ for distinct $1 \leq l, l' \leq N + M$. We combine (8.5) and (8.6) to bound (8.1) as

$$\bar{\mu}_\pi(Z_n \in \overline{\Gamma}, \mathbf{X}_n \text{ enters each set } C_i \text{ at most once})$$

$$\leq \sum_{1 \leq H \leq N+M} \sum_{\underline{i}_H} \sum_{\mathbf{k}_H} \sum_{\mathbf{s}} \left[ \prod_{j=1}^{H} \gamma(k_{j-1}, (i_{j-1}, i_j)) \right]$$

$$\times \exp\left\{ -\sum_{j=1}^{H} (k_j - k_{j-1}) \mathbb{I}^\theta_{i_j}(\Delta^m_{s_j}) \right\}.$$

[Note that (8.6) corresponds to index $H = 1$.]

Since the sum over $\mathbf{s} \in F(H, m, \mathbf{v}_{\mathbf{k}_H}, \overline{B})$ contains at most $(2^{m-1}(2^d))^H$ terms, we can apply Lemma 5.1 to obtain

$$\limsup \frac{1}{n} \log \bar{\mu}_\pi(Z_n \in \overline{\Gamma}, \mathbf{X}_n \text{ enters each set } C_i \text{ at most once})$$

$$(8.7) \quad \leq \limsup \max_{1 \leq H \leq N+M} \max_{\underline{i}_H} \max_{\mathbf{k}_H} \max_{\mathbf{s}}$$

$$\times \sum_{j=0}^{H-1} \frac{1}{n} \log \gamma(k_j, (i_j, i_{j+1})) - \sum_{j=1}^{H} \frac{(k_j - k_{j-1})}{n} \mathbb{I}^\theta_{i_j}(\Delta^m_{s_j}).$$

STEP 7. Now, by the choice of $\theta$, we have $\mathbb{I}^\theta_l = \mathbb{I}_l$ on $Q_l$ for $1 \leq l \leq N + M$. Also, recall that $\mathbb{I}_l$ is uniformly continuous on $Q_l$ for $1 \leq l \leq N + M$ (Proposition 2.1). Then, for $\mathbf{s} \in F(H, m, \mathbf{v}_{\mathbf{k}_H}, \overline{B})$ such that $\Delta(m, \mathbf{s}) \cap Q_{i_1} \times \cdots \times Q_{i_H} \neq \varnothing$, we have

$$\frac{1}{n} \sum_{j=1}^{H} (k_j - k_{j-1}) \mathbb{I}^\theta_{i_j}(\Delta^m_{s_j})$$

$$= \frac{1}{n} \sum_{j=1}^{H} (k_j - k_{j-1}) \inf_{x_j \in \Delta^m_{s_j} \cap Q_j} \mathbb{I}_{i_j}(x_j)$$

$$(8.8) \quad \geq \inf_{\mathbf{x} \in \Delta(m, \mathbf{s}) \cap \prod_{l=1}^{H} Q_{i_l}} \frac{1}{n} \sum_{j=1}^{H} (k_j - k_{j-1}) \mathbb{I}_{i_j}(x_j) - \alpha$$

$$= \inf_{\mathbf{x} \in \Delta(m,\mathbf{s}) \cap \prod_{l=1}^{H} Q_{i_l}} \frac{1}{n} \sum_{j=1}^{H} (k_j - k_{j-1}) \mathbb{I}_{i_j}^{\theta}(x_j) - \alpha$$

$$\geq \inf_{\mathbf{x} \in \Delta(m,\mathbf{s})} \frac{1}{n} \sum_{j=1}^{H} (k_j - k_{j-1}) \mathbb{I}_{i_j}^{\theta}(x_j) - \alpha.$$

On the other hand, if there exists $G \subset \{1, \ldots, H\}$ such that $\Delta_{s_j}^m \cap Q_{i_j} = \varnothing$ for all $j \in G$, we have that $\mathbb{I}_{i_j}^{\theta}(\Delta_{s_j}^m) = \inf_{x_j \in \Delta_{s_j}^m} \mathbb{I}_{i_j} = \theta$. Then, combining with (8.8), we have

(8.9)
$$\frac{1}{n} \sum_{j=1}^{H} (k_j - k_{j-1}) \mathbb{I}_{i_j}^{\theta}(\Delta_{s_j}^m)$$
$$= \frac{1}{n} \sum_{j \in G^c} (k_j - k_{j-1}) \mathbb{I}_{i_j}^{\theta}(\Delta_{s_j}^m) + \frac{1}{n} \sum_{j \in G} (k_j - k_{j-1}) \mathbb{I}_{i_j}^{\theta}(\Delta_{s_j}^m)$$
$$\geq \inf_{\mathbf{x} \in \Delta(m,\mathbf{s})} \frac{1}{n} \sum_{j=1}^{H} (k_j - k_{j-1}) \mathbb{I}_{i_j}^{\theta}(x_j) - \alpha.$$

With the estimate (8.9), we have that (8.7) is less than

(8.10)
$$\limsup \max_{H} \max_{\underline{i}_H} \max_{\mathbf{k}_H} \sum_{j=0}^{H-1} \frac{1}{n} \log \gamma(k_j, (i_j, i_{j+1}))$$
$$= \inf_{\mathbf{x} \in D_m(H, \mathbf{v}_{\mathbf{k}_H}, \overline{B})} \frac{1}{n} \sum_{j=1}^{H} (k_j - k_{j-1}) \mathbb{I}_{i_j}^{\theta}(x_j) + \alpha.$$

STEP 8. Without loss of generality, we may assume that the lim sup sequence in (8.10) occurs on a subsequence with fixed $1 \leq H \leq N + M$, $\underline{i}_H$ and vectors $\mathbf{k}_H^n$, where

$$\mathbf{v}_{\mathbf{k}_H^n} = \mathbf{v}^n \to \mathbf{v} = \langle v_1, \ldots, v_H \rangle$$

and

$$\lim_{n \to \infty} \frac{1}{n} \log \gamma(k_j^n, (i_j, i_{j+1})) \qquad \text{exists for } 1 \leq j \leq H.$$

Whereas values of $\theta$ and $m$ above a certain range are arbitrary, by Lemma 8.2 we have

$$\limsup_{\theta \uparrow \infty} \limsup_{m \uparrow \infty} \lim_{n \to \infty} \inf_{\mathbf{x} \in D_m(H, \mathbf{v}^n, \overline{B})} \sum_{j=1}^{H} v_j^n \mathbb{I}_{i_j}^{\theta}(x_j)$$
$$\geq \inf_{\mathbf{x} \in D(H, \mathbf{v}, \overline{B}) \cap \mathbb{K}^H} \sum_{j=1}^{H} v_j \mathbb{I}_{i_j, \varepsilon_1, \varepsilon_2}(x_j).$$

STEP 9. We now argue that

$$(8.11) \quad \lim \frac{1}{n} \log \gamma(k_j^n, (i_j, i_{j+1})) \leq \left( \sum_{l=1}^{j} v_l \right) \limsup \frac{1}{n} \log \gamma(n, (i_j, i_{j+1})).$$

Indeed, by definition $\sum_{l=1}^{j} v_l = \lim k_j^n/n$ for $1 \leq j \leq N + M$. Then, whereas

$$\frac{1}{n} \log \gamma(k_j^n, (i_j, i_{j+1})) = \frac{k_j^n}{n} \frac{1}{k_j^n} \log \gamma(k_j^n, (i_j, i_{j+1})),$$

inequality (8.11) follows easily when $\sum_{l=1}^{l} v_l > 0$ or $0 \geq \limsup(\log \gamma(k_j^n, (i_j, i_{j+1})))/n > -\infty$, but in the exceptional case, (8.11) still holds: Whereas $\log \gamma(k_j^n, (i_j, i_{j+1})) \leq 0$, we have by the convention $0 \cdot (-\infty) = 0$ that

$$\lim \frac{1}{n} \log \gamma(k_j^n, (i_j, i_{j+1})) \leq 0 = 0 \cdot (-\infty)$$

$$= \left( \sum_{l=1}^{j} v_l \right) \limsup \frac{1}{n} \log \gamma(n, (i_j, i_{j+1})).$$

Therefore, we have

$$\sum_{j=0}^{H-1} \lim \frac{1}{n} \log \gamma(k_j^n, (i_j, i_{j+1}))$$

$$\leq \mathbb{1}_{[H \geq 2]} \sum_{j=0}^{H-1} \left( \sum_{l=1}^{j} v_l \right) \limsup \frac{1}{n} \log \gamma(n, (i_j, i_{j+1}))$$

$$= \mathbb{1}_{[H \geq 2]} \sum_{j=1}^{H-1} \left( \sum_{l=1}^{j} v_l \right) \mathcal{U}_0(i_j, i_{j+1})$$

from Lemma 8.1, where the indicator reflects that the right-hand side vanishes when $H = 1$. So (8.10) is bounded above by

$$\mathbb{1}_{[H \geq 2]} \sum_{j=0}^{H-1} \left( \sum_{l=1}^{j} v_l \right) \mathcal{U}_0(i_j, i_{j+1}) - \inf_{\mathbf{x} \in D(H, \mathbf{v}, \Gamma) \cap \mathbb{K}^H} \sum_{j=1}^{H} v_j \mathbb{I}_{i_j, \varepsilon_1, \varepsilon_2}(x_j) + \alpha$$

$$\leq -\min_{H} \min_{\underline{i}_H} \min_{\mathbf{v} \in \Omega(H)} -\mathbb{1}_{[H \geq 2]} \sum_{j=0}^{H-1} \left( \sum_{l=1}^{j} v_l \right) \mathcal{U}_0(i_j, i_{j+1})$$

$$+ \inf_{\mathbf{x} \in D(H, \mathbf{v}, \overline{B}) \cap \mathbb{K}^H} \sum_{j=1}^{H} v_j \mathbb{I}_{i_j, \varepsilon_1, \varepsilon_2}(x_j) + \alpha$$

$$\leq -\mathbb{J}_{\mathcal{U}_0, \varepsilon_1, \varepsilon_2}(\overline{B}) + \alpha = -\mathbb{J}_{\mathcal{U}_0, \varepsilon_1, \varepsilon_2}(\overline{\Gamma} \cap \mathbb{K}) + \alpha.$$

Whereas $\alpha$ is arbitrary, the proposition follows. $\square$

**9. Lower coarse graining bounds.** As with the lower surgery estimate, the plan is to restrict the process to conveniently chosen events to derive lower bounds. Recall the notation $A'_n(l)$, $\underline{i}_H$, $C(\underline{i}_H)$, $\mathbf{k}_H$ and $S_n(\mathbf{k}_H)$ from Sections 5 and 8. Also, for $l \in \mathcal{G}$, let $\mathbb{P}^l_\eta$ denote the homogeneous nonnegative measure on $C_l$ with transition matrix $P(l)$ and initial distribution $\eta$.

PROOF OF PROPOSITION 4.9. Whereas $\pi$ is SIE-1 positive, let $e_l \in C_l$ be such that $\pi(e_l) > 0$ for $l \in \mathcal{G}$. Now, when $M = 1$, $\mathcal{G} = \{\zeta_1\}$, $\mathbb{J}_{\mathcal{T}_1} = \mathbb{I}_{\zeta_1}$ and on the set $G_n$, the process never leaves $C_{\zeta_1}$. In this case,

$$\underline{\mu}_\pi(Z_n \in \Gamma_3, G_n) \geq \mathbb{P}^{\zeta_1}_{e_{\zeta_1}}(Z_n \in \Gamma_3)$$

and the desired lower bound follows from Proposition 2.1.

Suppose that $M \geq 2$.

STEP 1. Let now $\underline{i}_M = \langle i_1, \ldots, i_M \rangle$, where $i_j \in \mathcal{G}$ for $1 \leq j \leq M$, be a given ordering of the nondegenerate irreducible sets $\mathcal{G}$. Let also $\Omega^+_M = \{\mathbf{v} \in \Omega_M : v_i > 0 \text{ for } 1 \leq i \leq M\}$ be the set of positive measures and let $\mathbf{v} \in \Omega^+_M$. Define also $v(0) = 0$ and $v(u) = \sum_{j=1}^u v_j$ for $1 \leq u \leq M$ and, in addition, for $n$ large enough so that $\lfloor nv(u) \rfloor < \lfloor nv(u+1) \rfloor$ for $1 \leq u \leq M-1$, that $\mathbf{k}^n = \langle \lfloor nv(1) \rfloor, \ldots, \lfloor nv(M-1) \rfloor \rangle$.

Then, for all large $n$,

$$\begin{aligned}
\underline{\mu}_\pi(Z_n \in \Gamma_3, G_n) \\
\geq \underline{\mu}_\pi(Z_n \in \Gamma_3, A'_n(M-1), C(\underline{i}_M)) \\
= \sum_{\mathbf{k}_M} \underline{\mu}_\pi(Z_n \in \Gamma_3, A'_n(M-1), C(\underline{i}_M), S_n(\mathbf{k}_M)) \\
\geq \underline{\mu}_\pi(Z_n \in \Gamma_3, A'_n(M-1), C(\underline{i}_M), S_n(\mathbf{k}^n)).
\end{aligned} \tag{9.1}$$

STEP 2. Whereas $\Gamma_3$ is open, the set $D(M, \mathbf{v}, \Gamma_3)$ [cf. (2.6)] is also open. Then, for $\mathbf{x} \in D(M, \mathbf{v}, \Gamma_3)$, let $\varepsilon > 0$ be so small so that the open cube $\Delta^\varepsilon(\mathbf{x})$ about $\mathbf{x}$ with side length $\varepsilon$ is contained: $\Delta^\varepsilon(\mathbf{x}) = \prod_{j=1}^M \Delta^\varepsilon(x_j) \subset D(M, \mathbf{v}_M, \Gamma_3)$. Also, for simplicity, let

$$E_n = A'_n(M-1) \cap C(\underline{i}_M) \cap S_n(\mathbf{k}^n)$$

and

$$a_n(u, \mathbf{v}) = \frac{1}{v_u} \frac{\lfloor nv(u) \rfloor - \lfloor nv(u-1) \rfloor}{n}$$

for $1 \leq u \leq M$. Then (9.1) equals

$$\begin{aligned}
\underline{\mu}_\pi(\langle a_n(1, \mathbf{v}) Z_1^{\lfloor nv(1) \rfloor}, \ldots, a_n(M, \mathbf{v}) Z_{\lfloor nv(M-1) \rfloor + 1}^n \rangle \in D(M, \mathbf{v}, \Gamma_3), E_n) \\
\geq \underline{\mu}_\pi(\langle a_n(1, \mathbf{v}) Z_1^{\lfloor nv(1) \rfloor}, \ldots, a_n(M, \mathbf{v}) Z_{\lfloor nv(M-1) \rfloor + 1}^n \rangle \in \Delta^\varepsilon(\mathbf{x}), E_n).
\end{aligned} \tag{9.2}$$

STEP 3. To make notation easier, we now get rid of the $a_n(u, \mathbf{v})$ terms at the cost of a further lower bound. Namely, because $Z_i^u \in \mathbb{K}$ is bounded for all $1 \leq i \leq u$, and $a_n(u, \mathbf{v}) \to 1$ for $1 \leq u \leq M$, we have for all $n$ large enough that

$$\{\langle Z_1^{\lfloor nv(1) \rfloor}, Z_{\lfloor nv(1) \rfloor+2}^{\lfloor nv(2) \rfloor}, \ldots, Z_{\lfloor nv(M-1) \rfloor+2}^n \rangle \in \Delta^{\varepsilon/2}(\mathbf{x})\}$$
$$\subset \{\langle a_n(1, \mathbf{v}) Z_1^{\lfloor nv(1) \rfloor}, \ldots, a_n(M, \mathbf{v}) Z_{\lfloor nv(M-1) \rfloor+1}^n \rangle \in \Delta^{\varepsilon}(\mathbf{x})\}.$$

Therefore, dropping the superscript $\Delta(\mathbf{x}) = \Delta^{\varepsilon/2}(\mathbf{x})$, we have for large $n$ that

$$(9.3) \quad (9.2) \geq \underline{\mu}_\pi(\langle Z_1^{\lfloor nv(1) \rfloor}, Z_{\lfloor nv(1) \rfloor+2}^{\lfloor nv(2) \rfloor}, \ldots, Z_{\lfloor nv(M-1) \rfloor+2}^n \rangle \in \Delta(\mathbf{x}), E_n).$$

STEP 4. We now decompose (9.3) in terms of resting and routing transitions. Recall that the transition probability between states $x \in C_l$ and $y \in C_m$ at time $n$ with respect to $\underline{\mu}_\pi$ equals $\check{\gamma}^0(n+1, (l,m))$ and does not depend on atoms $x$ and $y$.

Bound (9.3) below by

$$\pi(e_{i_1})\underline{\mu}_{e_{i_1}}(\langle Z_1^{\lfloor nv(1) \rfloor}, Z_{\lfloor nv(1) \rfloor+2}^{\lfloor nv(2) \rfloor}, \ldots, Z_{\lfloor nv(M-1) \rfloor+2}^n \rangle \in \Delta(\mathbf{x}),$$
$$X_{\lfloor nv(u-1) \rfloor} = e_{i_u} \text{ and } X_{\lfloor nv(u-1) \rfloor+1} = e_{i_{u+1}} \text{ for } 2 \leq u \leq M-1, E_n)$$

$$(9.4) \quad = \prod_{u=1}^{M-1} \check{\gamma}^0(\lfloor nv(u) \rfloor + 1, (i_u, i_{u+1})) \cdot \mathbb{P}_{e_{i_1}}^{i_1}(Z_1^{\lfloor nv(1) \rfloor} \in \Delta(x_1), X_{\lfloor nv(1) \rfloor} = e_{i_2})$$
$$\times \prod_{u=2}^{M-1} \mathbb{P}_{(\lfloor nv(u-1) \rfloor+1, e_{i_u})}^{i_u}(Z_{\lfloor nv(u-1) \rfloor+2}^{\lfloor nv(u) \rfloor} \in \Delta(x_u), X_{\lfloor nv(u) \rfloor} = e_{i_{u+1}})$$
$$\times \mathbb{P}_{(\lfloor nv(M-1) \rfloor+1, e_{i_M})}^{i_M}(Z_{\lfloor nv(M-1) \rfloor+2}^n \in \Delta(x_M)).$$

STEP 5. Observe, by definition, for distinct $i, j \in \mathcal{G}$, that

$$\liminf \frac{1}{k} \log \check{\gamma}^0(k, (i,j)) = \liminf \frac{1}{k} \log \min_{0 \leq r \leq E(N,M)} \underline{\gamma}^1(k+r, (i,j)) = \mathcal{T}_1(i,j).$$

Then, because large deviations of finite time-homogeneous irreducible chains are independent of the first and last observations, we have

$$\liminf \frac{1}{n} \log(9.4)$$
$$\geq \sum_{u=1}^{M-1} \left( \liminf \frac{\lfloor nv(u) \rfloor + 1}{n} \right)$$

$$(9.5) \qquad \times \left( \liminf \frac{1}{\lfloor nv(u) \rfloor + 1} \log \check{\gamma}^0(\lfloor nv(u) \rfloor + 1, (i_u, i_{u+1})) \right)$$

$$- \sum_{u=1}^{M} v_u \mathbb{I}_{i_u}(\Delta(x_u))$$

$$\geq \sum_{u=1}^{M-1} v(u) \mathcal{T}_1(\zeta_{i_u}, \zeta_{i_{u+1}}) - \sum_{u=1}^{M} v_u \mathbb{I}_{i_u}(x_u).$$

STEP 6. Whereas $\mathbf{v} \in \Omega_M^+$, $\mathbf{x} \in D(M, \mathbf{v}, \Gamma)$ and arrangement $\underline{i}_M$ composed of members in $\mathcal{G}$ are arbitrary, we have from (9.5) that

$$(9.6) \qquad \liminf \frac{1}{n} \log \underline{\mu}_\pi(Z_n \in \Gamma_3, G_n) \geq \sup_{\mathbf{v} \in \Omega_M^+} \max_{\sigma \in \mathbb{S}_M} g(\mathbf{v}, \sigma),$$

where

$$g(\mathbf{v}, \sigma) = \sum_{u=1}^{M-1} v(u) \mathcal{T}_1(\zeta_{\sigma(u)}, \zeta_{\sigma(u+1)}) - \inf_{\mathbf{y} \in D(M, \mathbf{v}, \Gamma_3)} \sum_{u=1}^{M} v_u \mathbb{I}_{\zeta_{\sigma(u)}}(y_u).$$

We now argue that we can replace $\Omega_M^+$ with the larger $\Omega_M$ in (9.6). In Lemma 9.1 below we show, for each $\sigma$, that $g(\cdot, \sigma)$ is lower semicontinuous as a function on $\Omega_M$. In particular, because $\mathbb{S}_M$ is a finite set, $\max_{\sigma \in \mathbb{S}_M} g(\cdot, \sigma)$ is lower semicontinuous. Therefore, by taking limits, we improve the bound in (9.6) to

$$\liminf \frac{1}{n} \log \underline{\mu}_\pi(Z_n \in \Gamma_3, G_n) \geq \sup_{\mathbf{v} \in \Omega_M} \max_{\sigma \in \mathbb{S}_M} g(\mathbf{v}, \sigma),$$

which is identified as $-\inf_{z \in \Gamma_3} \mathbb{J}_{\mathcal{T}_1}(z)$. $\square$

LEMMA 9.1. *Let $B \subset \mathbb{R}^d$ be an open set, and let $M \geq 2$ and $\sigma \in \mathbb{S}_M$. Then $g(\cdot, \sigma) : \Omega_M \to [0, \infty]$ is lower semicontinuous.*

PROOF. Let $\{\mathbf{v}^n\} \subset \Omega_M$ be a sequence which converges, $\mathbf{v}^n \to \mathbf{v}$. Recalling our convention $0 \cdot (-\infty) = 0$, we note that $h_1(\mathbf{v}) = \sum_{u=1}^{M-1} v(u) \mathcal{T}_1(\zeta_{\sigma(u)}, \zeta_{\sigma(u+1)})$ is lower semicontinuous, so we need only to prove $h_2(\mathbf{v}) = \inf_{\mathbf{y} \in D(M, \mathbf{v}, B)} \times \sum_{u=1}^{M} v_u \mathbb{I}_{\zeta_{\sigma(u)}}(y_u)$ is upper semicontinuous.

Let now $\mathbf{w} \in D(M, \mathbf{v}, B)$. Because $B$ is open and $\mathbf{v}^n$ converges to $\mathbf{v}$, we must have $\mathbf{w} \in D(M, \mathbf{v}^n, B)$ for all large $n$. Then,

$$\limsup h_2(\mathbf{v}^n) \leq \limsup \sum_{u=1}^{M} v_u^n \mathbb{I}_{\zeta_{\sigma(u)}}(w_u) = \sum_{u=1}^{M} v_u \mathbb{I}_{\zeta_{\sigma(u)}}(w_u).$$

However, because $\mathbf{w} \in D(M, \mathbf{v}, B)$ is arbitrary, we have in fact that

$$\limsup h_2(\mathbf{v}^n) \leq \inf_{\mathbf{y} \in D(M,\mathbf{v},B)} \sum_{u=1}^M v_u \mathbb{I}_{\zeta_{\sigma(u)}}(y_u) = h_2(\mathbf{v}). \qquad \square$$

**10. Limit estimate on $\mathbb{J}_{\mathcal{U}_0,\varepsilon_1,\varepsilon_2}$.** The proof of Proposition 4.5 follows in two steps (Propositions 10.1 and 10.2). The first step is to take $\varepsilon_1 \downarrow 0$ and estimate in terms of a quantity independent of degenerate transient sets $\mathcal{D}$ in Proposition 10.1. The second step is to let $\varepsilon_2 \downarrow 0$ and recover $\mathbb{J}_{\mathcal{U}_0}$ in the limit in Proposition 10.2.

It will be helpful to reduce the expression $\mathbb{J}_{\mathcal{U}_0,\varepsilon_1,\varepsilon_2}$ for $\varepsilon_1, \varepsilon_2 > 0$ [cf. (4.2)]. Whereas $\mathbb{I}_{i,\varepsilon_1,\varepsilon_2}$ is degenerate around $f(i)$ for $i \in \mathcal{D}$ [cf. (4.1)], we can evaluate $\mathbb{J}_{\mathcal{U}_0,\varepsilon_1,\varepsilon_2}(B)$ for $B \subset \mathbb{R}^d$ and $N + M \geq 2$ as

$$\min_{\sigma \in \mathbb{S}_{N+M}} \inf_{\mathbf{v} \in \Omega_{N+M}} \inf_{\mathbf{x} \in D'(\mathbf{v})} \left\{ -\sum_{i=1}^{N+M-1} \mathcal{U}_0(\sigma(i), \sigma(i+1)) \left[ \sum_{j=1}^i v_j \right] \right.$$
$$\left. - \sum_{\sigma(i) \in \mathcal{D}} v_i \log \varepsilon_1 + \sum_{\sigma(i) \in \mathcal{G}} v_i I_{\sigma(i),\varepsilon_2}(x_i) \right\},$$

where $D'(\mathbf{v}) = \{\mathbf{x} \in D(N+M, \mathbf{v}, B) : x_i = f(\sigma(i)), \text{ for } \sigma(i) \in \mathcal{D}\}$. When $N + M = 1$, the formula collapses to $\mathbb{J}_{\mathcal{U}_0,\varepsilon_1,\varepsilon_2} = \mathbb{I}_{1,\varepsilon_2}$.

We describe now an $\varepsilon_2 \geq 0$ "perturbation" of $\mathbb{J}_{\mathcal{U}_0}$, where we replace rates $\mathbb{I}_i$ with $\mathbb{I}_{i,\varepsilon_2}$ for $i \in \mathcal{G}$. Define, for Borel $B \subset \mathbb{R}^d$ and $M \geq 2$, that

$$\mathbb{J}_{\mathcal{U}_0}^{\varepsilon_2}(B) = \min_{\sigma \in \mathbb{S}_M} \inf_{\mathbf{v} \in \Omega_M} \inf_{\mathbf{x} \in D(M,\mathbf{v},B)}$$
$$- \sum_{i=1}^{M-1} \mathcal{U}_0(\zeta_{\sigma(i)}, \zeta_{\sigma(i+1)}) \left[ \sum_{j=1}^i v_j \right] + \sum_{i=1}^M v_i I_{\zeta_{\sigma(i)},\varepsilon_2}(x_i).$$

When $M = 1$, let $\mathbb{J}_{\mathcal{U}_0}^{\varepsilon_2} = \mathbb{I}_{1,\varepsilon_2}$.

We give now a triangle cost bound useful for the first step.

LEMMA 10.1. *For distinct $i, j, k \in \mathcal{G}$,*

$$\mathcal{U}_0(i,j) + \mathcal{U}_0(j,k) \leq \mathcal{U}_0(i,k).$$

PROOF. By definition, for some $k_1$ and distinct elements $L^1 = \langle l_0^1 = i, l_1^1, \ldots, l_{k_1}^1, l_{k_1+1}^1 = j \rangle$ we have $\mathcal{U}_0(i,j) = \sum_{s=0}^{k_1} \upsilon(l_s^1, l_{s+1}^1)$. Similarly, we have for some $k_2$ and $L^2 = \langle l_0^2 = j, l_1^2, \ldots, l_{k_2}^2, l_{k_2+1}^2 = k \rangle$ that $\mathcal{U}_0(j,k) = \sum_{s=0}^{k_2} \upsilon(l_s^2, l_{s+1}^2)$. Let now $T$ be the first index of an element in $L^1$ which belongs to $L^2$. Clearly, $1 < T \leq k_1 + 1$. Call also $T'$ the index of this element in $L^2$.

Form now $L^3 = \langle l_0^1, l_1^1, \ldots, l_T^1, l_{T'+1}^2, \ldots, l_{k_2+1}^2 \rangle$. From construction, $L^3$ is a list of distinct elements which we relabel as $L^3 = \langle l_0^3, \ldots, l_{k_3}^3 \rangle$ for some $k_3$.

Now, since $\upsilon(a,b) \leq 0$ for all distinct $a,b$, we have

$$\sum_{s=0}^{k_3} \upsilon(l_s^3, l_{s+1}^3) \geq \sum_{s=0}^{k_1} \upsilon(l_s^1, l_{s+1}^1) + \sum_{s=0}^{k_2} \upsilon(l_s^2, l_{s+1}^2).$$

However,

$$\mathcal{U}_0(i,k) = \max_{0 \leq k \leq M-2} \max_{L_k} \sum_{s=0}^{k} \upsilon(l_s, l_{s+1}) \geq \sum_{s=0}^{k_3} \upsilon(l_s^3, l_{s+1}^3)$$

$$\geq \mathcal{U}_0(i,j) + \mathcal{U}_0(j,k). \qquad \square$$

PROPOSITION 10.1. *Let $B \subset \mathbb{K}$ be a compact set and fix $\varepsilon_2 \geq 0$. Then, we have*

$$\liminf_{\varepsilon_1 \downarrow 0} \mathbb{J}_{\mathcal{U}_0, \varepsilon_1, \varepsilon_2}(B) \geq \mathbb{J}_{\mathcal{U}_0}^{\varepsilon_2}(B).$$

PROOF. First, when $N = 0$, we inspect that $\mathbb{J}_{\mathcal{U}_0, \varepsilon_1, \varepsilon_2}(B) = \mathbb{J}_{\mathcal{U}_0}^{\varepsilon_2}(B)$. Therefore, we assume that $N \geq 1$ in the following procedure.

STEP 1. Let $\varepsilon(k) \downarrow 0$, $\mathbf{v}^{\varepsilon(k)}$, $\mathbf{x}^{\varepsilon(k)}$ and $\sigma_{\varepsilon(k)}$ be sequences so that the limit inferior is attained:

$$\liminf_{\varepsilon_1 \downarrow 0} \mathbb{J}_{\mathcal{U}_0, \varepsilon_1, \varepsilon_2}(B)$$

(10.1)
$$= \lim_{k \to \infty} - \sum_{i=1}^{N+M-1} \mathcal{U}_0(\sigma_{\varepsilon(k)}(i), \sigma_{\varepsilon(k)}(i+1)) \left[ \sum_{j=1}^{i} v_j^{\varepsilon(k)} \right]$$

$$- \sum_{\sigma_{\varepsilon(k)}(i) \in \mathcal{D}} v_i^{\varepsilon(k)} \log \varepsilon(k) + \sum_{\sigma_{\varepsilon(k)}(i) \in \mathcal{G}} v_i^{\varepsilon(k)} \mathbb{I}_{\sigma_{\varepsilon(k)}(i), \varepsilon(k)}(x_i^{\varepsilon(k)}).$$

Because $\Omega_{N+M}$ is compact and $\mathbb{S}_{N+M}$ is finite, a further subsequence may be found so that, with the same labels, $\mathbf{v}^{\varepsilon(k)} \to \mathbf{v}$ and $\sigma_{\varepsilon(k)} = \sigma$ for all small $\varepsilon_1$.

STEP 2. When $\sum_{\sigma(i) \in \mathcal{D}} v_i > 0$, we have (10.1) diverges to $\infty$, which is automatically greater than the right-hand side in the proposition. On the other hand, if $\sum_{\sigma(i) \in \mathcal{D}} v_i = 0$, we must have $\sum_{\sigma(i) \in \mathcal{G}} v_i = 1$, because $\mathbf{v}$ is a probability vector. Now, if $(10.1) = \infty$, the proposition bound again holds.

Suppose therefore that (10.1) is finite. Recall that cube $\mathbb{K}$ contains the the domains of finiteness of the rate functions $\{\mathbb{I}_{i,\varepsilon_2} : i \in \mathcal{G}\}$ (cf. Proposition 2.1). Therefore, by taking a subsequence and relabeling, we can take $\mathbf{x}^{\varepsilon(k)} \in D'(\mathbf{v}^{\varepsilon(k)}) \cap \mathbb{K}$ and ensure the sequence is convergent, $\mathbf{x}^{\varepsilon(k)} \to \mathbf{x}$. Moreover, $\mathbf{x} \in D(N+M, \mathbf{v}, B)$ since $\sum_{i=1}^{N+M} v_i^{\varepsilon(k)} x_i^{\varepsilon(k)} \in B$ converges to $\sum_{i=1}^{N+M} v_i x_i$ and $B$ is closed.

Then, because $-\sum_{\sigma(i) \in \mathcal{D}} v_i^{\varepsilon(k)} \log \varepsilon(k) \geq 0$ and the rate functions $\mathbb{I}_{i,\varepsilon_2}$ are lower semicontinuous, we have that

$$(10.1) \geq \liminf_{k \to \infty} - \sum_{i=1}^{N+M-1} \mathcal{U}_0(\sigma(i), \sigma(i+1)) \left[ \sum_{j=1}^{i} v_j^{\varepsilon(k)} \right]$$

$$(10.2) \qquad + \sum_{\sigma(i) \in \mathcal{G}} v_i^{\varepsilon(k)} \mathbb{I}_{\sigma(i),\varepsilon_2}(x_i^{\varepsilon(k)})$$

$$\geq - \sum_{i=1}^{N+M-1} \mathcal{U}_0(\sigma(i), \sigma(i+1)) \left[ \sum_{j=1}^{i} v_j \right] + \sum_{\sigma(i) \in \mathcal{G}} v_i \mathbb{I}_{\sigma(i),\varepsilon_2}(x_i).$$

STEP 3. When $M = 1$ and $N \geq 1$, then $\mathcal{G} = \{\zeta_1\}$ is a singleton and $v_{\zeta_1} = 1$. Moreover, whereas $-\mathcal{U}_0$ is nonnegative, (10.2) is bounded below by $\mathbb{I}_{\zeta_1,\varepsilon_2}(x_{\zeta_1}) \geq \mathbb{J}_{\mathcal{U}_0}^{\varepsilon_2}(B)$ to finish the proof in this case.

STEP 4. Suppose then that $M \geq 2$ and $N \geq 1$. The strategy is to form a permutation $\eta \in \mathbb{S}_M$ and vector $\mathbf{u} \in \Omega_M$ for which (10.2) reduces to an expression that involves only terms that relate to $\mathcal{G}$. Write $\sigma^{-1}(\mathcal{G}) = \{\chi_1, \ldots, \chi_M\}$, where $\chi_i$ is ordered as follows:

$$\chi_1 = \min\{s : \sigma(s) \in \mathcal{G}\} \quad \text{and}$$
$$\chi_i = \min\{s > X_{i-1} : \sigma(s) \in \mathcal{G}\} \qquad \text{when } 2 \leq i \leq M.$$

Now, whereas $v_i = 0$ for $\sigma(i) \notin \mathcal{G}$ and, in particular, $v_i = 0$ for $1 \leq i \leq \chi_1 - 1$ when $\chi_1 \geq 2$, we have

$$- \sum_{i=1}^{N+M-1} \mathcal{U}_0(\sigma(i), \sigma(i+1)) \left[ \sum_{j=1}^{i} v_j \right]$$

$$= - \sum_{i=\chi_1}^{N+M-1} \mathcal{U}_0(\sigma(i), \sigma(i+1)) \sum_{\substack{\chi_1 \leq j \leq i \\ j \in \sigma^{-1}(\mathcal{G})}} v_j$$

$$= - \sum_{k=1}^{M-1} \sum_{i=\chi_k}^{\chi_{k+1}-1} \mathcal{U}_0(\sigma(i), \sigma(i+1)) \sum_{\substack{\chi_1 \leq j \leq i \\ j \in \sigma^{-1}(\mathcal{G})}} v_j + K_0,$$

where

$$K_0 = \begin{cases} -\sum_{i=\chi_M}^{N+M-1} \mathcal{U}_0(\sigma(i),\sigma(i+1)) \left[ \sum_{\substack{1 \leq j \leq i \\ j \in \sigma^{-1}(\mathcal{G})}} v_j \right], & \text{when } \chi_M < N+M, \\ 0, & \text{when } \chi_M = N+M. \end{cases}$$

In any case, because $K_0$ is nonnegative, we have that

(10.3)
$$-\sum_{i=1}^{N+M-1} \mathcal{U}_0(\sigma(i),\sigma(i+1)) \left[ \sum_{j=1}^{i} v_j \right]$$
$$\geq -\sum_{k=1}^{M-1} \left[ \sum_{i=\chi_k}^{\chi_{k+1}-1} \mathcal{U}_0(\sigma(i),\sigma(i+1)) \sum_{\substack{\chi_1 \leq j \leq i \\ j \in \sigma^{-1}(\mathcal{G})}} v_j \right].$$

STEP 5. We now bound individually the terms in large brackets in (10.3). For each $\chi_k \leq i \leq \chi_{k+1}-1$, as $\{v_j : \chi_1 \leq j \leq i \text{ and } j \in \sigma^{-1}(\mathcal{G})\} = \{v_{\chi_s} : 1 \leq s \leq k\}$, we may write

$$\sum_{i=\chi_k}^{\chi_{k+1}-1} \mathcal{U}_0(\sigma(i),\sigma(i+1)) \sum_{\substack{\chi_1 \leq j \leq i \\ j \in \sigma^{-1}(\mathcal{G})}} v_j$$
$$= [\mathcal{U}_0(\sigma(\chi_k),\sigma(\chi_k+1)) + \mathcal{U}_0(\sigma(\chi_k+1),\sigma(\chi_k+2))$$
$$+ \cdots + \mathcal{U}_0(\sigma(\chi_{k+1}-1),\sigma(\chi_{k+1}))] \left[ \sum_{s=1}^{k} v_{\chi_s} \right]$$
$$\leq \mathcal{U}_0(\sigma(\chi_k),\sigma(\chi_{k+1})) \sum_{s=1}^{k} v_{\chi_s}$$

by repeatedly applying the triangle inequality Lemma 10.1.

Hence, pulling together the inequalities, we have

(10.4)
$$-\sum_{i=1}^{N+M-1} \mathcal{U}_0(\sigma(i),\sigma(i+1)) \left[ \sum_{j=1}^{i} v_j \right]$$
$$\geq -\sum_{k=1}^{M-1} \mathcal{U}_0(\sigma(\chi_k),\sigma(\chi_{k+1})) \left[ \sum_{s=1}^{k} v_{\chi_s} \right].$$

STEP 6. Define now $\mathbf{u} \in \Omega_M$ by $u_k = v_{\chi_k}$ for $1 \leq k \leq M$. Then

$$\sum_{i \in \sigma^{-1}(\mathcal{G})} v_i \mathbb{I}_{\sigma(i),\varepsilon_2}(x_i) = \sum_{k=1}^{M} v_{\chi_k} \mathbb{I}_{\sigma(\chi_k),\varepsilon_2}(x_{\chi_k}) = \sum_{k=1}^{M} u_k \mathbb{I}_{\sigma(\chi_k),\varepsilon_2}(x_{\chi_k}).$$

Now let $\eta \in \mathbb{S}_M$ be the permutation where $\zeta_{\eta(i)} = \sigma(\chi_i)$ for $1 \leq i \leq M$. Noting (10.4), we can then bound (10.2) below by

$$
\begin{aligned}
(10.5) \quad & -\sum_{k=1}^{M-1} \mathcal{U}_0(\sigma(\chi_k), \sigma(\chi_{k+1})) \left[\sum_{s=1}^{k} u_s\right] + \sum_{k=1}^{M} u_k \mathbb{I}_{\sigma(\chi_k),\varepsilon_2}(x_{\chi_k}) \\
& = -\sum_{k=1}^{M-1} \mathcal{U}_0(\zeta_{\eta(k)}, \zeta_{\eta(k+1)}) \left[\sum_{s=1}^{k} u_s\right] - \sum_{k=1}^{M} u_k \mathbb{I}_{\zeta_{\eta(k)},\varepsilon_2}(x_{\chi_k}).
\end{aligned}
$$

STEP 7. By construction $\sum_{i=1}^{N+M} v_i x_i \in B$. Then, because $v_j = 0$ when $\sigma(j) \notin \mathcal{G}$, we have

$$
\sum_{j=1}^{N+M} v_j x_j = \sum_{j \in \sigma^{-1}(\mathcal{G})} v_j x_j = \sum_{j \in \sigma^{-1}(\mathcal{G})} v_j x_j = \sum_{s=1}^{M} v_{\chi_s} x_{\chi_s} = \sum_{s=1}^{M} u_s x_{\chi_s}
$$

and so $\langle x_{\chi_1}, \ldots, x_{\chi_M} \rangle \in D(M, \mathbf{u}, B)$. Hence, tracing through the argument,

$$
\begin{aligned}
(10.5) \geq & \inf_{\mathbf{x} \in D(M,\mathbf{u},B)} -\sum_{k=1}^{M-1} \mathcal{U}_0(\zeta_{\eta(k)}, \zeta_{\eta(k+1)}) \left[\sum_{s=1}^{k} u_s\right] - \sum_{k=1}^{M} u_k I_{\zeta_{\eta(k)},\varepsilon_2}(x_k) \\
\geq & \mathbb{J}^{\varepsilon_2}_{\mathcal{U}_0}(B). \qquad \square
\end{aligned}
$$

PROPOSITION 10.2. *Let $\Gamma \subset \mathbb{R}^d$ be compact. Then we have*

$$(10.6) \qquad \liminf_{\varepsilon \downarrow 0} \mathbb{J}^{\varepsilon}_{\mathcal{U}_0}(\Gamma) \geq \mathbb{J}_{\mathcal{U}_0}(\Gamma).$$

PROOF. When $\liminf_{\varepsilon \downarrow 0} \mathbb{J}^{\varepsilon}_{\mathcal{U}_0}(\Gamma) = \infty$, of course (10.6) is immediate.

STEP 1. Suppose then that $\liminf_{\varepsilon \downarrow 0} \mathbb{J}^{\varepsilon}_{\mathcal{U}_0}(\Gamma) < \infty$. As in Step 2 in Proposition 10.1, let $\varepsilon(k) \downarrow 0$, $\sigma_{\varepsilon(k)} = \sigma$ independent of $k$, $\mathbf{v}^{\varepsilon(k)} \to \mathbf{v}$ and $\mathbf{x}^{\varepsilon(k)} \to \mathbf{x} \in D(M, \mathbf{v}, \Gamma)$ be such that

$$
\begin{aligned}
& \liminf_{\varepsilon \downarrow 0} \mathbb{J}^{\varepsilon}_{\mathcal{U}_0}(\Gamma) \\
& = \lim_{k \to \infty} -\sum_{i=1}^{M-1} \left(\sum_{j=1}^{i} v_j^{\varepsilon(k)}\right) \mathcal{U}_0(\zeta_{\sigma(i)}, \zeta_{\sigma(i+1)}) + \sum_{i=1}^{M} v_i^{\varepsilon(k)} \mathbb{I}_{\zeta_{\sigma(i)},\varepsilon(k)}(x_i^{\varepsilon(k)}).
\end{aligned}
$$

STEP 2. We now claim for $i \in \mathcal{G}$ that

$$(10.7) \qquad \liminf_{k \to \infty} \mathbb{I}_{i,\varepsilon(k)}(x_i^{\varepsilon(k)}) \geq \mathbb{I}_i(x_i).$$

For $\lambda \in \mathbb{R}^d$, let $\rho_{i,\varepsilon}(\lambda)$ and $\rho_i(\lambda)$ be the Perron–Frobenius eigenvalues that correspond to the $\lambda$ tilts of $P(i,\varepsilon)$ and $P(i)$ [cf. (2.3)]. From [21], we have that $\lim_{\varepsilon \downarrow 0} \log \rho_{i,\varepsilon}(\lambda) = \log \rho_i(\lambda)$.

Now, for $\lambda' \in \mathbb{R}^d$, observe that

$$\liminf_k \mathbb{I}_{i,\varepsilon(k)}(x_i^{\varepsilon(k)}) = \liminf_k \sup_{\lambda \in \mathbb{R}^d} \langle \lambda, x_i^{\varepsilon(k)} \rangle - \log \rho_{i,\varepsilon(k)}(\lambda)$$

$$\geq \liminf_k \langle \lambda', x_i^{\varepsilon(k)} \rangle - \log \rho_{i,\varepsilon(k)}(\lambda')$$

$$= \langle \lambda', x \rangle - \log \rho_i(\lambda').$$

Hence, because $\lambda'$ is arbitrary, we have $\liminf_k \mathbb{I}_{i,\varepsilon(k)}(x_i^{\varepsilon(k)}) \geq \sup_\lambda \{\langle \lambda, x \rangle - \log \rho_i(\lambda)\} = \mathbb{I}_i(x)$.

STEP 3. In fact, (10.7) proves the proposition when $M = 1$. On the other hand, when $M \geq 2$, we have with (10.7) that

$$\liminf_{\varepsilon \downarrow 0} \mathbb{J}_{\mathcal{U}_0}^\varepsilon(\Gamma) \geq -\sum_{i=1}^{M-1} \left(\sum_{j=1}^i v_j\right) \mathcal{U}_0(\zeta_{\sigma(i)}, \zeta_{\sigma(i+1)}) + \sum_{i=1}^M v_i \mathbb{I}_{\zeta_{\sigma(i)}}(x_i)$$

$$\geq \mathbb{J}_{\mathcal{U}_0}(\Gamma). \qquad \square$$

**11. Routing cost comparisons.** We separate the proof of Proposition 4.10 into two separate results.

PROPOSITION 11.1. *Suppose Assumption* B *holds. Then, for distinct* $i, j \in \mathcal{G}(P)$,

$$\mathcal{T}_1(i,j) \geq \mathcal{T}_0(i,j).$$

PROOF. Recall the definitions of $\underline{\gamma}^1(n, y, z)$ and $\underline{\gamma}^1(n, (i, j))$. It is enough to prove for $y \in C_i$ and $z \in C_j$ that

(11.1) $$\liminf_{n \to \infty} \frac{1}{n} \log \underline{\gamma}^1(n, y, z) \geq \mathcal{T}_0(i,j).$$

Then, clearly

$$\mathcal{T}_1(i,j) = \liminf_{n \to \infty} \frac{1}{n} \log \underline{\gamma}^1(n, (i,j)) \geq \mathcal{T}_0(i,j),$$

finishing the proof.

We now show (11.1). Let $k$ and $L_k = \langle i = l_0, l_1, \ldots, l_k, l_{k+1} = j \rangle$ be such that

$$\mathcal{T}_0(i,j) = \sum_{s=0}^k \tau(l_s, l_{s+1}).$$

To connect with the definition of $\underline{\gamma}^1(n,(i,j))$, form vectors $\mathbf{x}^0 = \langle x_1^0, \ldots, x_{q_0}^0 \rangle, \ldots, \mathbf{x}^{k+1} = \langle x_1^{k+1}, \ldots, x_{q_{k+1}}^{k+1} \rangle$ with distinct elements in $C_{l_0}, \ldots, C_{l_{k+1}}$ such that, for $0 \leq s \leq k$,

$$x_{q_s}^s = a(l_s, l_{s+1}) \quad \text{and} \quad x_1^{s+1} = b_{n+r(s)}(l_s, l_{s+1}),$$

$1 \leq q_s \leq \mathfrak{r} + 1$ and $1 \leq q_0, q_{k+1} \leq \mathfrak{r}$. In addition, because $\{P(i) : i \in \mathcal{G}\}$ are irreducible, we specify that the paths are possible. Namely, for all large $n$,

$$\mathbb{P}_{(n-1,y)}(\mathbf{X}_n^{n+r(0)-1} = \mathbf{x}^0) \geq (p_{\min}/2)^{\mathfrak{r}},$$

$$\mathbb{P}_{(n+r(s-1), x_1^s)}(\mathbf{X}_{n+r(s-1)+1}^{n+r(s)-1} = \mathbf{x}_2^s) \geq (p_{\min}/2)^{\mathfrak{r}}$$

and

$$\mathbb{P}_{(n+r(k+1)-1, x_{q_{k+1}}^{k+1})}(X_{n+r(k+1)} = z) \geq p_{\min}/2$$

when $q_s \geq 2$ and $1 \leq s \leq k+1$. Here, $\mathbf{x}_2^s = \langle x_2^s, \ldots, x_{q_s}^s \rangle$ when $q_s \geq 2$, $r(s) = \sum_{u=0}^{s} q_u$ and $p_{\min}$ is defined in (2.5).

Since the length of the connecting path from $y$ to $z$ is at most $E_0(N, M)$ [cf. near (2.10)], we have

$$\liminf_{n \to \infty} \frac{1}{n} \log \underline{\gamma}^1(n, y, z)$$

$$\geq \liminf \left[ \frac{\log(p_{\min}/2)^{E_0(N,M)}}{n} \right.$$

$$\left. + \frac{1}{n} \sum_{s=0}^{k+1} \log p_{n+r(s)}(a(l_s, l_{s+1}), b_{n+r(s)}(l_s, l_{s+1})) \right]$$

$$= \sum_{s=0}^{k+1} \tau(l_s, l_{s+1}) = \mathcal{T}_0(i, j)$$

from Assumption B. □

PROPOSITION 11.2. *Suppose that Assumption C holds. Then, for distinct $i, j \in \mathcal{G}$,*

$$\mathcal{T}_1(i,j) \geq \mathcal{T}_0(i,j).$$

PROOF. The proof is similar to that of Proposition 11.1. As before, it is enough to show (11.1). Let $k$ and $L_k = \langle i = l_0, l_1, \ldots, l_k, l_{k+1} = j \rangle$ be such that $\mathcal{T}_0(i,j) = \sum_{s=0}^{k} \tau(l_s, l_{s+1})$. Form the path vector $\mathbf{x}^0 = \langle x_1^0, \ldots, x_{q_0}^0 \rangle$ with $1 \leq q_0 \leq \mathfrak{r}$ of distinct elements in $C_i$ and state $x_1^1 \in C_{l_1}$ such that

$$p_{n-1+(q_0+1)}(x_{q_0}^0, x_1^1) = t(n + q_0, (i, l_1))$$

and
$$\liminf_{n\to\infty} \frac{1}{n} \log \mathbb{P}_{(n-1,y)}(\mathbf{X}_n^{n+q_0-1} = \mathbf{x}^0)$$
$$= \liminf \frac{1}{n}[\log p_n(y, x_1^0) + \log p_n(x_1^0, x_2^0) + \cdots + \log p_n(x_{q_0-1}^0, x_{q_0}^0)] = 0.$$

Such a vector $\mathbf{x}^0$ exists from the primitivity of $P^*(i)$.

Similarly, form vectors $\mathbf{x}^s = \langle x_1^s, \ldots, x_{q_s}^s \rangle$ in $C_{l_s}$, where $1 \leq q_s \leq \mathfrak{r}+1$ for $1 \leq s \leq k$ and $1 \leq q_{k+1} \leq \mathfrak{r}$. Also specify that
$$p_{n-1+r(s)+1}(x_{q_s}^s, x_1^{s+1}) = t(n - 1 + r(s) + 1, (l_s, l_{s+1}))$$
for $1 \leq s \leq k$. In addition, the paths are chosen so
$$\liminf \frac{1}{n} \log \mathbb{P}_{(n-1+r(s-1)+1, x_1^s)}(\mathbf{X}_{n+r(s-1)+1}^{n+r(s)-1} = \mathbf{x}_2^s) = 0$$
and
$$\liminf \frac{1}{n} \log \mathbb{P}_{(n+r(k+1)-1, x_{q_{k+1}}^{k+1})}(X_{n+r(k+1)} = z) = 0$$
when $l_s \in \mathcal{G}$ and $q_s \geq 2$, and $\mathbf{x}_2^s$ and $r(s)$ are as before. Then
$$\liminf_{n\to\infty} \frac{1}{n} \log \underline{\gamma}^1(n, y, z) \geq \liminf \frac{1}{n} \sum_{s=0}^{k+1} \log t(n - 1 + r(s), (l_s, l_{s+1}))$$
$$\geq \sum_{s=0}^{k+1} \tau(l_s, l_{s+1}) = \mathcal{T}_0(i, j). \qquad \square$$

**12. Examples.** In this section, we present three examples that concern possible LD behaviors of $\{Z_n(f)\}$ under $\mathbb{P}_\pi \in \mathbb{A}(P)$. The first shows that even if Assumption A is violated, an LDP may still hold with respect to some processes and functions $f$. The second example shows that the bounds in Theorems 3.1 and 3.2(ii) may be achieved. The third example shows that it is possible that an LDP is nonexistent under Assumption A when one of the submatrices $\{P(i) : i \in \mathcal{G}\}$ is periodic and Assumptions B and C do not hold.

12.1. *Assumption* A *is not necessary for LDP.* The point is that if the connecting transition probabilities oscillate so that Assumption A fails, but not too wildly, then the process on the large deviation scale can wait an $o(n)$ time to select optimal connections. Let $\Sigma = \{0, 1\}$ and initial distribution $\pi = \langle 1/2, 1/2 \rangle$. Let also $f : \Sigma \to \mathbb{R}$ be given by $f(0) = 1$ and $f(1) = 0$, and for $k \geq 1$, define transition matrices
$$A_k = \begin{bmatrix} 1 - (\frac{1}{2})^k & (\frac{1}{2})^k \\ 0 & 1 \end{bmatrix} \quad \text{and} \quad B_k = \begin{bmatrix} 1 - (\frac{1}{3})^k & (\frac{1}{3})^k \\ 0 & 1 \end{bmatrix}.$$

Then, for $n \geq 1$, let
$$P_n = \begin{cases} A_n, & \text{for } n \text{ even,} \\ B_n, & \text{for } n \text{ odd.} \end{cases}$$
The limit matrix $P$ is the $2 \times 2$ identity matrix $I_2$, with two irreducible sets, $C_0 = \{0\}$ and $C_1 = \{1\}$. Both sets correspond to degenerate rate functions, for $i = 0, 1$, $\mathbb{I}_i(x) = 0$ for $x = 1 - i$ and $= \infty$ otherwise. Also, one sees that $\tau(0,1) = -\log 3 < -\log 2 = v(0,1)$, so Assumption A is not satisfied here. Of course, $\tau(1,0) = v(1,0) = -\infty$. Also, the process satisfies Condition SIE-1.

To identify the large deviations of $\{Z_n(f)\}$ under $\mathbb{P}_\pi^{\{P_n\}}$, we focus on sets $\Gamma = (a, b]$ for $0 < a < b < 1$, because the analysis on other types of sets is similar.

As before, $A(0)$ and $A(1)$ are the events that $\mathbf{X}_n$ does not switch and switches exactly once between sets $C_0$ and $C_1$. Since $\Gamma$ is such that $\mathbb{P}_\pi(Z_n \in \Gamma, A(0)) = 0$ and also since the chain cannot switch from state 1 to 0, we have
$$\mathbb{P}_\pi(Z_n \in \Gamma) = \mathbb{P}_\pi(Z_n \in \Gamma, A(1)) = \mathbb{P}_\pi(Z_n \in \Gamma, A(1), X_1 = 0, X_n = 1).$$
The event $\{A(1), X_1 = 0, X_n = 1\} \subset \Sigma^n$ consists exactly of $n-1$ paths $\mathbf{x}_{n,i}$ that start at 0 but switch to 1 at time $1 \leq i \leq n-1$. Now compute that
$$\mathbb{P}_\pi(\mathbf{X}_n = \mathbf{x}_{n,i}) = \pi(0) \prod_{k=1}^{i} (1 - \alpha(k)^k)(\alpha(i+1))^{i+1} \prod_{l=i+2}^{n} (1 - \alpha(l)^l)$$
$$= e^{o(n)} (\alpha(i+1))^{i+1},$$
where $\alpha(k) = 1/2$ for $k$ even and $= 1/3$ for $k$ odd. Also, on the path $\mathbf{x}_{n,i}$, we have that $Z_n = i/n$.

Let $G_n^o = \{1 \leq i \leq n : i/n \in \Gamma^o\}$. Then, by Lemma 5.1, we have
$$\liminf \frac{1}{n} \log \mathbb{P}_\pi(Z_n \in \Gamma^o, A(1), X_1 = 0, X_n = 1)$$
$$= \liminf \max_{i \in G_n^o} \frac{1}{n} \log \mathbb{P}_\pi(\mathbf{X}_n = \mathbf{x}_{n,i})$$
$$= \liminf \max\left\{\frac{\lceil an \rceil}{n} \log(\alpha(\lceil an \rceil)), \frac{\lceil an \rceil + 1}{n} \log(\alpha(\lceil an \rceil + 1))\right\}$$
$$= a \log\left(\frac{1}{2}\right) = -a \log(2).$$
Similarly, $\limsup (1/n) \log \mathbb{P}_\pi(Z_n \in \overline{\Gamma}, A(1), X_1 = 0, X_n = 1) = -a \log(2)$.

A related analysis works for more general $\Gamma$ and we have that $\{Z_n(f)\}$ satisfies an LDP with rate function
$$\mathbb{I}(z) = \begin{cases} z \log 2, & z \in [0, 1), \\ 0, & z = 1, \\ \infty, & \text{otherwise.} \end{cases}$$

12.2. *Bounds may be sharp in Theorems* 3.1 *and* 3.2. The key in this example is that the connection probabilities oscillate "unboundedly," so picking out the optimal strategy is time-dependent. As before, let $\Sigma = \{0, 1\}$, $\pi = \langle 1/2, 1/2 \rangle$ and let $f : \Sigma \to \mathbb{R}$ be given by $f(0) = 1$ and $f(1) = 0$. Let $\{g(n)\}$ be a fast divergent sequence of integers, $g(n) \uparrow \infty$, $g(n) < g(n+1)$ and $g(n-1)/g(n) \to 0$. Also, for $k \geq 1$, let

$$P_i = \begin{cases} I_2, & \text{for } 1 \leq i \leq g(2), \\ A_i, & \text{for } g(2k) < i \leq g(2k+1), \\ B_i, & \text{for } g(2k+1) < i \leq g(2k+2), \end{cases}$$

where $A_i$ and $B_i$ are defined in Section 12.1.

To compute the large deviations of $\{Z_n(f)\}$, we focus now on sets $\Gamma = (a, b) \subset [0, 1]$, where $0 \leq a < b < 1$. Calculations for other sets are analogous. Then, in the notation of the previous example,

$$\liminf \frac{1}{n} \log \mathbb{P}_\pi(Z_n \in \Gamma) = \liminf \frac{1}{n} \log \mathbb{P}_\pi(Z_n \in \Gamma, A(1), X_1 = 0, X_n = 1).$$

Let now $n_k = g(2k+2)$ for $k \geq 1$. Then $i/n_k \in \Gamma$ exactly when $\lceil g(2k+2)a \rceil \leq i \leq \lfloor g(2k+2)b \rfloor$. Also, whereas

$$\lim_n \frac{\lceil g(2k+2)a \rceil}{g(2k+2)} = a > 0 = \lim \frac{g(2k+1)}{g(2k+2)},$$

we have for all large $k$ that $g(2k+1) + 1 \leq \lceil g(2k+2)a \rceil \leq \lfloor g(2k+2)b \rfloor \leq g(2k+2)$. Note also that $P_i = B_i$ for $g(2k+1) + 1 \leq i \leq g(2k+2)$. Hence,

$$\lim \frac{1}{n_k} \log \mathbb{P}_\pi(Z_{n_k} \in \Gamma, A(1), X_1 = 0, X_{n_k} = 1)$$

$$= \liminf \max_{i : i/n_k \in \Gamma} \frac{1}{n_k} \log \mathbb{P}_\pi(\mathbf{X}_{n_k} = \mathbf{x}_{n_k, i})$$

$$= \lim \frac{\lceil g(2n+2)a \rceil}{g(2n+2)} \log\left(\frac{1}{3}\right) = -a \log(3).$$

Moreover, in fact $\liminf (1/n) \log \mathbb{P}_\pi(Z_n \in \Gamma) = -a \log(3)$.

Similarly, by considering subsequence $n_k = g(2k+1)$, we get

$$\limsup \frac{1}{n} \log \mathbb{P}_\pi(Z_n \in \overline{\Gamma}, A(1), X_1 = 0, X_n = 1) = -a \log(2).$$

These calculations, and analogous ideas give, for any $\Gamma$, that

$$\limsup \frac{1}{n} \log \mu_\pi(Z_n \in \overline{\Gamma}) = -\inf_{z \in \overline{\Gamma}} \overline{\mathbb{J}}(z)$$

and

$$\liminf \frac{1}{n} \log \mu_\pi(Z_n \in \Gamma^o) = -\inf_{z \in \Gamma^o} \underline{\mathbb{J}}(z),$$

where

$$\overline{\mathbb{J}}(z) = \begin{cases} z \log 2, & \text{for } z \in [0,1), \\ 0, & z = 1, \\ \infty, & \text{otherwise,} \end{cases}$$

and

$$\underline{\mathbb{J}}(z) = \begin{cases} z \log 3, & \text{for } z \in [0,1), \\ 0, & z = 1, \\ \infty, & \text{otherwise.} \end{cases}$$

On the other hand, these lower and upper rate functions match those in Theorems 3.1 and 3.2(i). Whereas $\mathcal{T}_0(0,1) = -\log 3$, $\mathcal{U}_0(0,1) = -\log 2$ and $t(k,(1,0)) = 0$ for all $k \geq 1$, we have

$$\mathbb{J}_{\mathcal{T}_0}(z) = -\inf_{\delta \in [0,1]} \inf_{\langle x,y \rangle \in D(2,\langle \delta, 1-\delta \rangle, z)} \min\{\mathbb{I}_0(y), \delta \log(3) + \delta \mathbb{I}_0(x) + (1-\delta)\mathbb{I}_1(y)\}$$
$$= \underline{\mathbb{J}}(z)$$

and analogously $\mathbb{J}_{\mathcal{U}_0} = \overline{\mathbb{J}}$.

12.3. *Periodicity and nonexistence of LDP.* We consider a process which satisfies Assumptions A but not Assumptions B or C for which an LDP cannot hold through an explicit contradiction. Also, we show that the lower bound with respect to $\mathcal{T}_0$ in Theorem 3.2 does not work for this example.

Let $\Sigma = \{1, \ldots, 9\}$ and let $\pi$ be the uniform distribution on $\Sigma$. For $n$ of the form $n = 1 + 3j$ for $j \geq 0$, except when $n = 3^{2^j} + 1$ for $j \geq 5$, let

$$\overline{P}_n = \begin{bmatrix} 1/3 & 1/3 & 1/3 & 0 & 0 & 0 & 0 & 0 & 0 \\ 1/3 & 1/3 & 1/3 & 0 & 0 & t(n,(1,2)) & 0 & 0 & 0 \\ 1/3 & 1/3 & 1/3 & 0 & 0 & 0 & 0 & 0 & 0 \\ 0 & 0 & 0 & 0 & 1 & 0 & 0 & 0 & t(n,(2,3)) \\ 0 & 0 & 0 & 0 & 0 & 1 & 0 & 0 & 0 \\ 0 & 0 & 0 & 1 & 0 & 0 & 0 & 0 & \underline{t}(n,(2,3)) \\ 0 & 0 & 0 & 0 & 0 & 0 & 1/3 & 1/3 & 1/3 \\ 0 & 0 & 0 & 0 & 0 & 0 & 1/3 & 1/3 & 1/3 \\ 0 & 0 & 0 & 0 & 0 & 0 & 1/3 & 1/3 & 1/3 \end{bmatrix},$$

$$\overline{P}_{n+1} = \begin{bmatrix} 1/3 & 1/3 & 1/3 & 0 & 0 & 0 & 0 & 0 & 0 \\ 1/3 & 1/3 & 1/3 & 0 & 0 & 0 & 0 & 0 & 0 \\ 1/3 & 1/3 & 1/3 & t(n+1,(1,2)) & 0 & 0 & 0 & 0 & 0 \\ 0 & 0 & 0 & 0 & 1 & 0 & \underline{t}(n+1,(2,3)) & 0 & 0 \\ 0 & 0 & 0 & 0 & 0 & 1 & t(n+1,(2,3)) & 0 & 0 \\ 0 & 0 & 0 & 1 & 0 & 0 & 0 & 0 & 0 \\ 0 & 0 & 0 & 0 & 0 & 0 & 1/3 & 1/3 & 1/3 \\ 0 & 0 & 0 & 0 & 0 & 0 & 1/3 & 1/3 & 1/3 \\ 0 & 0 & 0 & 0 & 0 & 0 & 1/3 & 1/3 & 1/3 \end{bmatrix},$$

$$\overline{P}_{n+2} = \begin{bmatrix} 1/3 & 1/3 & 1/3 & 0 & t(n+2,(1,2)) & 0 & 0 & 0 & 0 \\ 1/3 & 1/3 & 1/3 & 0 & 0 & 0 & 0 & 0 & 0 \\ 1/3 & 1/3 & 1/3 & 0 & 0 & 0 & 0 & 0 & 0 \\ 0 & 0 & 0 & 0 & 1 & 0 & 0 & 0 & 0 \\ 0 & 0 & 0 & 0 & 0 & 1 & 0 & \underline{t}(n+2,(2,3)) & 0 \\ 0 & 0 & 0 & 1 & 0 & 0 & 0 & t(n+2,(2,3)) & 0 \\ 0 & 0 & 0 & 0 & 0 & 0 & 1/3 & 1/3 & 1/3 \\ 0 & 0 & 0 & 0 & 0 & 0 & 1/3 & 1/3 & 1/3 \\ 0 & 0 & 0 & 0 & 0 & 0 & 1/3 & 1/3 & 1/3 \end{bmatrix}.$$

For $n = 3^{2^j} + 1$ for $j \geq 5$, let $\widehat{P}_{n+1}$ and $\widehat{P}_{n+2}$ be defined as before, but now let

$$\overline{P}_n = \begin{bmatrix} 1/3 & 1/3 & 1/3 & 0 & 0 & 0 & 0 & 0 & 0 \\ 1/3 & 1/3 & 1/3 & 0 & 0 & t(n,(1,2)) & 0 & 0 & 0 \\ 1/3 & 1/3 & 1/3 & 0 & 0 & 0 & 0 & 0 & 0 \\ 0 & 0 & 0 & 0 & 1 & 0 & 0 & 0 & t(n,(2,3)) \\ 0 & 0 & 0 & 0 & 0 & 1 & 0 & 0 & \underline{t}(n,(2,3)) \\ 0 & 0 & 0 & 1 & 0 & 0 & 0 & 0 & 0 \\ 0 & 0 & 0 & 0 & 0 & 0 & 1/3 & 1/3 & 1/3 \\ 0 & 0 & 0 & 0 & 0 & 0 & 1/3 & 1/3 & 1/3 \\ 0 & 0 & 0 & 0 & 0 & 0 & 1/3 & 1/3 & 1/3 \end{bmatrix}.$$

Suppose now that $t(n,(1,2))$, $t(n,(2,3))$ and $\underline{t}(n,(2,3))$ vanish as $n$ tends to infinity and limits

$$\lim \frac{1}{n} \log t(n(1,2)), \quad \lim \frac{1}{n} \log t(n,(2,3)) \quad \text{and} \quad \lim \frac{1}{n} \log \underline{t}(n,(2,3))$$

exist and equal, respectively,

$$v(1,2) = \tau(1,2) = 0, \quad v(2,3) = \tau(2,3) = A$$

and

$$\lim \frac{1}{n} \log \underline{t}(n,(2,3)) = 2A + \varepsilon,$$

where $A < 0$ and $\varepsilon > 0$ is chosen small enough so that $2A + \varepsilon < A$.

Define the diagonal matrix $\Delta_n = \text{diag}\{\lambda_1^{-1}, \ldots, \lambda_9^{-1}\}$, where $\lambda_i$ is the $i$th row sum of $\overline{P}_n$. Then $\lim \Delta_n = I_9$. Let $P_n = \Delta_n \overline{P}_n$ for $n \geq 1$. The limit matrix $P = \lim P_n = \lim \overline{P}_n$ corresponds to three sets: $C_1 = \{1,2,3\}$, $C_2 = \{4,5,6\}$ and $C_3 = \{7,8,9\}$.

Let also $f$ be a one-dimensional function on the state space such that $f(1) = f(2) = f(3) = 1$, $f(4) = f(5) = f(6) = 2$ and $f(7) = f(8) = f(9) = 3$. We now concentrate the sequence $\{Z_n(f)\}$ with respect to the process $\mathbb{P}_\pi^{\{P_n\}}$.

*Assumptions.* By inspection, it is clear that Condition SIE-1 and Assumption A hold, but Assumptions B and C do not hold.

*Nonexistence of LDP.* First, let $\mu_\pi$ be the measure constructed from $\{\overline{P}_n\}$ and $\pi$ through CON. It is not difficult to see that the large deviation of $Z_n$ under $\mathbb{P}_\pi$ is the same as with respect to $\mu_\pi$, that is, for Borel $\Gamma \subset \mathbb{R}^d$,

$$\limsup \frac{1}{n} \log \mathbb{P}_\pi(Z_n \in \overline{\Gamma}) = \limsup \frac{1}{n} \log \mu_\pi(Z_n \in \overline{\Gamma})$$

and

$$\liminf \frac{1}{n} \log \mathbb{P}_\pi(Z_n \in \Gamma^o) = \liminf \frac{1}{n} \log \mu_\pi(Z_n \in \Gamma^o)$$

(cf. Proposition 7.1). Second, the rate functions on the three sets are degenerate:

$$\mathbb{I}_i(z) = \begin{cases} 0, & \text{if } z = i, \\ \infty, & \text{otherwise,} \end{cases} \text{ for } i = 1, 2, 3.$$

Consider now the following two lemmas, which are proved later.

LEMMA 12.1. *For $0 < \varepsilon < 1/2$, let $\Gamma = [2+\varepsilon, 2+2\varepsilon]$. Then*

$$\limsup \frac{1}{n} \log \mu_\pi(Z_n \in \overline{\Gamma}) > (1-2\varepsilon)A.$$

LEMMA 12.2. *For $0 < \varepsilon < 1/2$ and $\theta > 0$, let $\Gamma(\theta) = (2+\varepsilon-\theta, 2+2\varepsilon+\theta)$. Then*

$$\liminf_{\theta \downarrow 0} \liminf_{n \to \infty} \frac{1}{n} \log \mu_\pi(Z_n \in \Gamma(\theta)) \leq (1-2\varepsilon)A.$$

These results show that no LDP is possible. If an LDP were to hold with rate function $I$, say, then

$$(1-2\varepsilon)A \geq \liminf_{\theta \downarrow 0} \liminf_{n \to \infty} \frac{1}{n} \log \mu_\pi(Z_n \in \Gamma(\theta))$$

$$\geq \liminf_{\theta \downarrow 0} - \inf_{x \in \Gamma(\theta)} I(x) \geq - \inf_{x \in \overline{\Gamma}} I(x)$$

$$\geq \limsup_{n \to \infty} \frac{1}{n} \log \mu_\pi(Z_n \in \overline{\Gamma}) > (1-2\varepsilon)A,$$

leading to a contradiction.

*Lower bound in Theorem* 3.2(ii) *does not hold.* Consider the following lemma proved at this end of this section.

LEMMA 12.3. *With respect to $\Gamma(\theta)$ as in Lemma 12.2, we have*

$$- \inf_{z \in \Gamma(\theta)} \mathbb{J}_{\mathcal{T}_0}(z) = (1-2\varepsilon-\theta)\frac{A}{2}.$$

Then a clear contradiction with Lemma 12.2 would arise if the lower bound in Theorem 3.2(ii) were valid.

PROOF OF LEMMA 12.1. No $n$-word $\mathbf{x}_n$ that remains solely in a single closed set can have an average in $\overline{\Gamma}$. Also, by construction, no $n$-word can pass from $C_i$ to $C_j$ for $i > j$, or in one step from $C_1$ to $C_3$. Therefore, the only $n$-words such that $(1/n)\sum_{i=1}^n f(x_i) \in [2+\varepsilon, 2+2\varepsilon]$ are those which visit succesively $C_1$, $C_2$ and $C_3$ or those which visit first $C_2$ and then $C_3$.

We now examine $(1/n)\log\mu_\pi(Z_n \in \overline{\Gamma})$ along the sequence

$$n_k = \left\lceil 3^{2^k} / \left(\frac{1-2\varepsilon}{2}\right) \right\rceil$$

for $k \geq 1$. Let now $A(n_k)$ be the set of $n_k$-words $\mathbf{x}_{n_k}$ which stays in $C_1$ until time $3^{2^k}$, spends one time unit in $C_2$ and then switches to $C_3$. By definition, for $\mathbf{x}_{n_k} \in A(n_k)$ and $k$ large enough, we have

$$\frac{1}{n_k}\sum_{i=1}^{n_k} f(x_i) = \frac{3^{2^k}}{n_k} + \frac{2}{n_k} + \left(1 - \frac{3^{2^k}+1}{n_k}\right) 3 \in [2+\varepsilon, 2+2\varepsilon].$$

Then, with $\delta(\varepsilon) = (1-2\varepsilon)/2$, we have

$$\limsup_{n\to\infty} \frac{1}{n}\log\mu_\pi^n(Z_n \in \overline{\Gamma})$$

$$\geq \liminf_{k\to\infty} \frac{1}{n_k}\log\mu_\pi^{n_k}(Z_{n_k} \in \overline{\Gamma}, A(n_k))$$

$$\geq - \inf_{\langle x,y\rangle \in D(2,\langle\delta(\varepsilon),1-\delta(\varepsilon)\rangle,\overline{\Gamma})}$$

$$- \delta(\varepsilon)\left\{\liminf \frac{1}{k}\log t(k,(1,2)) + \liminf \frac{1}{k}\log \underline{t}(k,(2,3))\right\}$$

$$+ \delta(\varepsilon)\mathbb{I}_1(x) + (1-\delta(\varepsilon))\mathbb{I}_3(y)$$

$$= \delta(\varepsilon)\left\{\lim \frac{1}{k}\log t(k,(1,2)) + \lim \frac{1}{k}\log \underline{t}(k,(2,3))\right\}$$

$$- \delta(\varepsilon)\mathbb{I}_1(1) - (1-\delta(\varepsilon))\mathbb{I}_3(3)$$

$$= \frac{1-2\varepsilon}{2}(2A+\varepsilon) > (1-2\varepsilon)A. \qquad \square$$

PROOF OF LEMMA 12.2. Let now $n_k = 3^{2^k}$ for $k \geq 1$. We first show that $\mathbf{x}_{n_k}$ cannot visit $C_1$, $C_2$ and $C_3$ in succession and satisfy $\frac{1}{n_k}\sum_{i=1}^{n_k} f(x_i) \in \Gamma(\theta)$ for all small $\theta$. Indeed, by construction, a path $\mathbf{x}_{n_k}$ which visits $C_1, C_2$ and $C_3$ must switch from $C_2$ to $C_3$ at a time less than or equal to $3^{2^{k-1}}$.

However, then, because $f(\cdot) \geq 1$ and $3^{2^{k-1}}/n_k \to 0$, we have for large $k$ and $\theta$ sufficiently small that

$$\frac{1}{n_k}\sum_{i=1}^{n_k} f(x_i) \geq \frac{3^{2^{k-1}}}{n_k} + \left(1 - \frac{3^{2^{k-1}}+1}{n_k}\right)3 > 2 + 2\varepsilon + \theta.$$

Thus, if $\mathbf{x}_{n_k} \in \Gamma(\theta)$, we deduce $\mathbf{x}_{n_k}$ begins in $C_2$ and then switches to $C_3$. Now let $\tau(\varepsilon) = 1 - 2\varepsilon$. We have

$$\liminf_{\theta \downarrow 0} \liminf_{n \to \infty} \frac{1}{n} \log P_\pi(Z_n \in \Gamma(\theta))$$

$$\leq \liminf_{\theta \downarrow 0} \liminf_{k \to \infty} \frac{1}{n_k} \log P_\pi(Z_{n_k} \in \Gamma(\theta))$$

$$= \liminf_{\theta \downarrow 0} \sup_{0 \leq \delta \leq 1} \sup_{\langle x,y \rangle \in D(2,\langle \delta, 1-\delta \rangle, \Gamma(\theta))} \delta \limsup \frac{1}{k} \log t(k,(2,3))$$
$$- \delta \mathbb{I}_2(x) - (1-\delta)\mathbb{I}_3(y)$$

$$= \tau(\varepsilon) \limsup \frac{1}{k} \log t(k,(2,3)) = (1 - 2\varepsilon)A,$$

because $\tau(\varepsilon)$ is the smallest $\delta$ such that $(2,3) \in D(2,\langle \delta, 1-\delta \rangle, [2-\varepsilon, 2+2\varepsilon])$. □

PROOF OF LEMMA 12.3. Since motion is possible only from $C_1$ to $C_2$ to $C_3$, and the corresponding rate functions are degenerate at $x_1 = 1$, $x_2 = 2$ and $x_3 = 3$, we have

$$\mathbb{J}_{\mathcal{T}_0}(\Gamma(\theta)) = \sup_{\substack{v_1+v_2+v_3=1 \\ 0 \leq v_1,v_2,v_3 \leq 1}} \sup_{\mathbf{x} \in D(3,\mathbf{v},\Gamma(\theta))} v_1 \tau(1,2) + (v_1+v_2)\tau(2,3) - \sum_{i=1}^{3} v_i \mathbb{I}_i(x_i)$$

$$= \sup_{\substack{v_1+2v_2+3(1-v_1-v_2) \in \Gamma(\theta) \\ 0 \leq v_1,v_2 \leq 1}} (v_1+v_2)A$$

$$= (1 - 2\varepsilon - \theta)(A/2). \qquad \Box$$

## APPENDIX

**A.1.** *Proof of Lemma* 4.1. We consider separately the situations when $0 < \delta < 1$ and $\delta = 0$.

CASE $\delta > 0$. Let $\hat{t}_n = \sup_{s \geq n} t_s$. Then $t_n \leq \hat{t}_n$, $0 < \hat{t}_n \leq 1$ and $\hat{t}_n \downarrow \delta$. Also

$$\lim \frac{1}{n} \log \hat{t}_n \to 0 = \limsup \frac{1}{n} \log t_n.$$

CASE $\delta = 0$. The proof is split into two subcases.

*Subcase* 1. $\limsup (1/n) \log t_n = t < 0$. If $t_n$ vanishes eventually, that is, $t_n = 0$ for $n \geq N_0$, some $N_0 \geq 1$, then we may take

$$\hat{t}_n = \begin{cases} 1, & \text{for } 1 \leq n < N_0, \\ e^{-n^2}, & \text{for } n \geq N_0. \end{cases}$$

Otherwise, let $a_n = \sup_{j \geq n}(1/j) \log t_j$ and $\hat{t}_n = \exp\{\sup_{l \geq n} l a_l\}$. Note $a_n \downarrow t$ and $\hat{t}_n \geq \exp\{n a_n\} \geq \exp\{n(1/n) \log t_n\} = t_n$, and also that $1 \geq \hat{t}_n > 0$. In addition, $(1/n) \log \hat{t}_n \geq a_n \to t$.

Let now $1 > \varepsilon > 0$ and let $N_1$ be such that $a_n < (1-\varepsilon)t$ for $n \geq N_1$. Then

$$\frac{1}{n} \log \hat{t}_n \leq \frac{1}{n} \sup_{l \geq n} lt(1-\varepsilon) = t(1-\varepsilon)$$

for $n \geq N_1$. Whereas $\varepsilon$ is arbitrary, we then have $(1/n) \log \hat{t}_n \to t$.

*Subcase* 2. $\limsup(1/n) \log t_n = t = 0$. As $t_n \to 0$, we have $t_n < 1$ for $n \geq N_2$, say. Let $b_j = \max_{N_2 \leq l \leq j}(1/l) \log t_l$ for $j \geq N_2$ and let

$$\hat{t}_n = \begin{cases} 1, & \text{for } n < N_2, \\ \exp\left\{\sup_{j \geq n} j b_j\right\}, & \text{for } n \geq N_2. \end{cases}$$

Note that $t_n \leq \hat{t}_n$ and $1 \geq \hat{t}_n > 0$, and as $\sup_{j \geq n} j b_j$ decreases with $n$, that $\hat{t}_n$ is a decreasing sequence.

We now identify the limit. Note that $b_j \leq 0$ for all $j \geq N_2$ and $(1/l) \log t_l \to 0$. Then, for each $K \geq N_2$, there is an index $J_K \geq K$ such that

$$b_j = \max_{K \leq l \leq j}(1/l) \log t_l \text{ for } j \geq J_K.$$

Hence, for large $n$ and given $K \geq N_2$,

$$\hat{t}_n = \exp\left\{\sup_{j \geq n} j \max_{K \leq l \leq j}(1/l) \log t_l\right\} \leq \exp\left\{\sup_{j \geq n} \max_{K \leq l \leq j} \log t_l\right\} = \sup_{j \geq n} \max_{K \leq l \leq j} t_l.$$

Whereas $K$ is arbitrary, we have that $\hat{t}_n \downarrow 0$.

Finally, as $b_j \to 0$, we have for $\varepsilon > 0$ and large $n$ that

$$0 \geq (1/n) \log \hat{t}_n = (1/n) \sup_{j \geq n} j b_j \geq (1/n) \sup_{j \geq n} j(-\varepsilon) = -\varepsilon.$$

Whereas $\varepsilon$ is arbitrary, we have $(1/n) \log \hat{t}_n \to 0$. □

A.2. *An extended Gärtner–Ellis theorem.* We give here a minor extension of the Gärtner–Ellis theorem and state some general conditions under which a sequence of bounded nonnegative measures $\{\mu_n\}$ on $\mathbb{R}^d$ satisfies an LDP.

For $\lambda \in \mathbb{R}^d$, define the extended real sequence $\Lambda_n(\lambda) = \log \int_{\mathbb{R}^d} e^{\langle \lambda, x \rangle} d\mu_n(x)$ and also $\Lambda(\lambda) = \lim_{n \to \infty} (1/n) \Lambda_n(n\lambda)$, provided the extended limit exists. We now recall when $\Lambda$ is *essential smoothness* (cf. [10]).

ASSUMPTION E.

1. For all $\lambda \in \mathbb{R}^d$, $\Lambda(\lambda)$ exists as an extended real number in $(-\infty, \infty]$.
2. Let $D_\Lambda = \{\lambda : -\infty < \Lambda(\lambda) < \infty\}$. Suppose $0 \in D_\Lambda^o$.
3. The function $\Lambda(\cdot)$ is differentiable throughout $D_\Lambda^o$.
4. When $\{\lambda_n\} \subset D_\Lambda^o$ converges to a boundary point of $D_\Lambda$, we have $|\nabla \Lambda(\lambda_n)| \to \infty$.
5. The function $\Lambda(\lambda)$ is a lower semicontinuous function.

We now state the standard Gärtner–Ellis theorem (cf. [10]).

PROPOSITION A.1. *Let $\{\nu_n\}$ be a sequence of probability measures which satisfy Assumption* E. *Let $\mathbb{I}$ be the Legendre transform of $\Lambda$. Then $\mathbb{I}$ is a rate function and $\{\nu_n\}$ satisfies LDP* (2.1).

The main result of this section is the following proposition.

PROPOSITION A.2. *Let $\{\mu_n\}$ be a sequence of bounded nonnegative measures on $\mathbb{R}^d$ that satisfy Assumption* E. *Let $\mathbb{I}$ be the Legendre transform of $\Lambda$. Then $\mathbb{I}$ is an extended rate function and the LDP* (2.1) *holds. Moreover, $\mathbb{I}$ can be decomposed as the difference of a rate function of a probability sequence and a constant, $\mathbb{I} = \mathbb{I}^1 - \Lambda(0)$.*

PROOF. By Assumption E, with $\lambda = 0$, we have that $(1/n) \log \mu_n(\mathbb{R}^d) \to \Lambda(0) \in \mathbb{R}$. Consider now the probability measures $\nu_n(\cdot) = \mu_n(\cdot)/\mu_n(\mathbb{R}^d)$. The pressure of the sequence $\{\nu_n\}$ is calculated as $\Lambda(\cdot) - \Lambda(0)$. Since Assumption E holds for $\Lambda(\cdot)$, it also holds for the shifted function $\Lambda(\cdot) - \Lambda(0)$. Therefore, by Proposition A.1, we have that $\{\nu_n\}$ satisfies (2.1) with rate function $\mathbb{I}^1$ given by

$$\mathbb{I}^1(x) = \sup_\lambda \{\langle \lambda, x \rangle - (\Lambda(\lambda) - \Lambda(0))\}$$
$$= \sup_\lambda \{\langle \lambda, x \rangle - \Lambda(\lambda)\} + \Lambda(0).$$

Let now $\mathbb{I}(x) = \sup_\lambda \{\langle \lambda, x \rangle - \Lambda(\lambda)\}$, so that $\mathbb{I} = \mathbb{I}^1 - \Lambda(0)$. Whereas $\mu_n(\cdot) = \mu_n(\mathbb{R}^d) \nu_n(\cdot)$, by translating we obtain that (2.1) holds for the $\{\mu_n\}$ sequence with rate function $\mathbb{I}$. $\square$

A.3. *Proof of Proposition* 2.1.

*Extended pressure* $\Lambda$. We follow the method in [10] to identify the extended pressure of the sequence $\{\mu_n\}$:

$$\begin{aligned}\Lambda(\lambda) &= \lim \frac{1}{n} \log \Lambda_n(n\lambda) \\ &= \lim \frac{1}{n} \log \int_{\mathbf{X}_n \in C^n} \exp\left(\left\langle \lambda, \sum f(X_i) \right\rangle\right) d\mathbb{U}_\pi \\ &= \lim \frac{1}{n} \log(\pi^t (\Pi_{C,\lambda})^n \mathbf{1}).\end{aligned}$$

Since $\Pi_{C,\lambda}$ is an irreducible matrix, the Perron–Frobenius eigenvalue $\rho(C,\lambda)$ possesses a right Perron–Frobenius eigenvector $\mathbf{v}(\lambda)$ with positive entries. Let $a$ and $b$ be the smallest and largest entries. Then

$$\log(\pi^t(\Pi_{C,\lambda})^n \mathbf{1}) \leq \log((1/a)\pi^t(\Pi_{C,\lambda})^n \mathbf{v}) = \frac{1}{n}\log\left(\frac{1}{a}\pi^t \mathbf{v}\right) + \log \rho(C,\lambda)$$

and, similarly, $\log(\pi^t(\Pi_{C,\lambda})^n \mathbf{1}) \geq \log \rho(C,\lambda) + o(1)$. Hence,

$$\Lambda(\lambda) = \lim \frac{1}{n} \log \Lambda_n(n\lambda) = \log \rho(C,\lambda).$$

*Analyticity, convexity and essential smoothness of* $\Lambda$. Perron–Frobenius theory guarantees that $\rho(\lambda)$ has multiplicity 1 and is positive for all $\lambda \in \mathbb{R}^d$. Then, by Theorem 7.7.1 in [21], $\rho(\cdot)$ is analytic and so $\Lambda(\cdot)$ is analytic. Now, because $\Lambda(\lambda)$ is the limit of a sequence of convex functions, it is convex. Finally, by the comments of Section 3.1 in [10], we have that $\Lambda$ is essentially smooth.

$\mathbb{I}$ *is an extended rate function and* $\{\mu_n\}$ *satisfies an LDP.* Recall now that $\mathbb{I} = \mathbb{I}_C$ is the Legendre transform $\mathbb{I}(x) = \sup_{\lambda \in \mathbb{R}^d} \langle \lambda, x \rangle - \Lambda(\lambda)$. By Proposition A.2, we have that $\mathbb{I}$ is an extended rate function and $\{\mu_n\}$ satisfies an LDP with respect to $\mathbb{I}$.

$\mathbb{I}$ *is a rate function when* $U_C$ *is substochastic.* When $U_C$ is substochastic, we have $\Lambda(0) \leq 0$. Hence, by Proposition A.2, $\mathbb{I} = \mathbb{I}^1 - \Lambda(0) \geq 0$ and so is a rate function.

$\mathbb{I}$ *is not identically* $\infty$. Let $\hat{x} = \nabla \Lambda(0)$. Then, by Theorem 23.5 in [25],

$$\mathbb{I}(\hat{x}) = \sup_{\lambda \in \mathbb{R}^d} \langle \lambda, \nabla \Lambda(0) \rangle - \Lambda(\lambda) = \langle 0, \nabla \Lambda(0) \rangle - \Lambda(0) = -\Lambda(0) < \infty.$$

*Convexity of $\mathbb{I}$ and strict convexity on the relative interior of $Q_C$.* Whereas $\Lambda$ is convex, the Legendre transform $\mathbb{I}$ is convex. Also, because $\Lambda(\cdot)$ is real-valued and lower semicontinuous, by Lemma 4.5.8 of [10], $\Lambda$ is the conjugate of $\mathbb{I}$. Whereas $\mathbb{I}$ is not identically $\infty$, it is a proper convex function. Moreover, since $\mathbb{I}$ is lower semicontinuous, it is a closed convex function as well (cf. [25], page 52). Then, since $\Lambda$ is essentially smooth, we have from Theorem 26.3 of [25] that $\mathbb{I}$ is strictly convex on the relative interior of its domain of finiteness $Q_C$.

*$Q_C$ is convex and $Q_C \subset \mathbb{K}$.* Let $x, y \in Q_C$. The convexity of $\mathbb{I}$ implies that $\mathbb{I}((x+y)/2) \leq (\mathbb{I}(x) + \mathbb{I}(y))/2 < \infty$. Hence, $Q_C$ is convex.

For $\lambda \in \mathbb{R}^d$, let $\bar{\lambda} = \langle |\lambda_1|, \ldots, |\lambda_d| \rangle$. Then

$$\exp\left\langle -\bar{\lambda}, \left(\max_i |f(i)|\right) \mathbf{1}_d \right\rangle P_C \leq \Pi_{C,\lambda} \leq \exp\left\langle \bar{\lambda}, \left(\max_i |f(i)|\right) \mathbf{1}_d \right\rangle P_C.$$

Whereas the Perron–Frobenius value of $P_C$ is 1, we have

$$\exp\left\langle -\bar{\lambda}, \left(\max_i |f(i)|\right) \mathbf{1}_d \right\rangle \leq \rho(\lambda) \leq \exp\left\langle \bar{\lambda}, \left(\max_i |f(i)|\right) \mathbf{1}_d \right\rangle.$$

Now let $x$ be such that $x_j > \max_i |f(i)|$ for some $1 \leq j \leq d$. Then, for $\alpha \in \mathbb{R}$, let $\lambda^{j,\alpha} \in \mathbb{R}^d$ be such that $\lambda_i^{j,\alpha} = 0$ for $i \neq j$ and $\lambda_j^{j,\alpha} = \alpha$. We have then

$$\mathbb{I}(x) \geq \sup_{\lambda \in \mathbb{R}^d} \langle \lambda, x \rangle - \left\langle \bar{\lambda}, \left(\max_i |f_i|\right) \mathbf{1}_d \right\rangle$$

$$\geq \langle \lambda^{j,\alpha}, x \rangle - \left\langle \overline{\lambda^{j,\alpha}}, \left(\max_i |f_i|\right) \mathbf{1}_d \right\rangle \geq \alpha x_j - |\alpha| \max_i |f_i|.$$

By taking $\alpha \uparrow \infty$, we have that $\mathbb{I}(x) = \infty$. Similarly, if $x_j < -\max_i |f_i|$, then $\mathbb{I}(x) = \infty$. Thus, $\mathbb{I}(x) < \infty$ implies $\max_i |x_i| < \max_i |f_i|$ and so $Q_C \subset \mathbb{K}$.

*$Q_C$ is compact.* If $\mathbb{I}$ can be shown to be uniformly bounded on $Q_C$, then the lower semicontinuity of $\mathbb{I}$ will imply that $Q_C$ is closed. Also, since it was shown above that $Q_C$ is bounded, $Q_C$ will then be compact.

Let $\underline{p}$ be the smallest positive entry in $P_C$ and let $G = \{x : \mathbb{I}(x) \leq -\log \underline{p}\}$. By the lower semicontinuity of $\mathbb{I}$, $G$ is a closed set. Let $x_0 \in G^c$. We show that $\mathbb{I}(x_0) = \infty$ and hence $Q_C \subset G$.

Since $G^c$ is open, let $B = \overline{B}(x_0; \delta) \subset G^c$ be a closed ball around $x_0$ with some radius $\delta > 0$. If now $\limsup (1/n) \log \mu_n(Z_n \in B) > -\infty$, then there exists a sequence $\{\mathbf{x}_{n_k}\}$ such that $\sum_{i=1}^{n_k} f(x_i)/n_k \in B$ and $\mu_n(\mathbf{X}_n = \mathbf{x}_{n_k}) > 0$. However, we have $\mu_n(\mathbf{X}_n = \mathbf{x}_{n_k}) \geq \underline{p}^{n_k}$, and so $\limsup (1/n_k) \log \mu_n(Z_{n_k} \in B) \geq \log \underline{p}$. Hence, using the LD upper bound,

$$-\mathbb{I}(B) \geq \limsup \frac{1}{n} \log \mu_n(Z_n \in B) \geq \log \underline{p}.$$

However, since $\mathbb{I}$ is lower semicontinuous, $\mathbb{I}(B) = \mathbb{I}(x_1)$ on some point $x_1$ in the compact set $B \subset G^c$. Hence, $\mathbb{I}(B) > -\log \underline{p}$, giving a contradiction.

Therefore, we must have $\mathbb{I}(x_0) = \infty$ because

$$-\infty = \limsup \frac{1}{n} \log \mu_n(Z_n \in B)$$
$$\geq \liminf \frac{1}{n} \log \mu_n(Z_n \in B) \geq -\mathbb{I}(B^o) \geq -\mathbb{I}(x_0).$$

$\mathbb{I}$ *is uniformly continuous on* $Q_C$. Whereas $\mathbb{I}$ is convex, $\mathbb{I}$ restricted to $Q_C$ is continuous. Since $Q_C$ is compact, $\mathbb{I}$ is in fact uniformly continuous on $Q_C$.

$\mathbb{I}$ *is a good rate function.* Whereas $\mathbb{I}$ is lower semicontinuous, the level set $\{x : \mathbb{I}(x) \leq a\}$ for $a \in \mathbb{R}$ is a closed subset of $Q_C$ and hence compact.

A.4. *Proof of Proposition* 4.1. When $M = 1$, $P(\zeta_1)$ is stochastic and $\mathbb{J}_U = \mathbb{I}_{\zeta_1}$, and the proof follows from Proposition 2.1. Suppose now that $M \geq 2$. Consider that $\mathbb{J}_U \leq \min\{\mathbb{I}_i : i \in \mathcal{G}\}$ and so $Q_{\mathbb{J}_U} \supset \bigcup_{i \in \mathcal{G}} Q_i$ is nonempty. Also, $Q_{\mathbb{J}_U} \subset \mathbb{K}$: Indeed, for $z \notin \mathbb{K}$, and $\mathbf{v} \in \Omega_M$ and $\mathbf{x} \in D(M, \mathbf{v}, z)$ we must have that $v_i > 0$ and $x_i \notin \mathbb{K}$ for some $1 \leq i \leq M$. Then $C_{\mathbf{v},U}(\sigma, \mathbf{x}) = \infty$ and so $\mathbb{J}_U(z) = \infty$.

In addition, $\mathbb{J}_U$ is lower semicontinuous and nonnegative because $\{\mathbb{I}_i\}$ that correspond to substochastic matrices $\{P(i) : i \in \mathcal{G}\}$ are rate functions with compact domains of finiteness. Finally, $\mathbb{J}_U$ is a good rate function from the same argument given for Proposition 2.1.

**Acknowledgments.** We thank the referee for careful remarks, which led to improvements in the article, and also for pointing out useful references. Also, we thank J. Schmalian for some helpful conversations.

DEPARTMENT OF MATHEMATICS
TULANE UNIVERSITY
6823 ST. CHARLES AVENUE
NEW ORLEANS, LOUISIANA 70118
USA
E-MAIL: zdietz@math.tulane.edu

DEPARTMENT OF MATHEMATICS
IOWA STATE UNIVERSITY
400 CARVER HALL
AMES, IOWA 50011
USA
E-MAIL: sethuram@iastate.edu